\documentclass[11pt,a4paper]{article}
\usepackage[margin=1in]{geometry}
\usepackage{lmodern}

\usepackage[skip=0.5em,indent=0pt]{ parskip}

\usepackage[utf8]{inputenc} %
\usepackage[T1]{fontenc}    %
\usepackage{hyperref}       %
\usepackage{url}            %
\usepackage{booktabs}       %
\usepackage{amsfonts}       %
\usepackage{nicefrac}       %
\usepackage{microtype}      %
\usepackage{xcolor}         %

\usepackage{amssymb,amsfonts}
\usepackage{caption}
\usepackage{subcaption}
\usepackage{graphicx}
\usepackage{enumitem}

\usepackage{url}
\usepackage{multicol}
\usepackage{multirow}
\usepackage{mathtools}
\usepackage{tabularray}
\usepackage{cals}
\usepackage{algorithm}
\usepackage[noend]{algpseudocode}
\usepackage{tikz}
\usepackage{color}
\usepackage{array}
\usepackage{multirow}
\usepackage{comment}
\usepackage{tikz}
\usepackage{tikz-cd}
\usetikzlibrary{arrows.meta, positioning}
\usetikzlibrary{shapes}
\usetikzlibrary{calc}
\usepackage{enumitem}
\usepackage{cleveref}

\usepackage{bm}

\usepackage[normalem]{ulem}

\DeclareMathOperator*{\argmin}{arg\,min}

\DeclareMathOperator{\tr}{tr}
\DeclareMathOperator{\KL}{KL}
\DeclareMathOperator{\Ent}{Ent}

\DeclareMathOperator{\proj}{proj}
\DeclareMathOperator{\E}{\mathbb{E}}

\newcommand{\Tan}{\mathrm{Tan}}

\DeclareMathOperator{\Unif}{Unif}

\DeclareMathOperator{\diag}{diag}

\DeclareMathOperator{\dom}{dom}

\DeclareMathSymbol{\mhyphen}{\mathord}{AMSa}{"39}

\definecolor{Blue}{HTML}{000000}
\definecolor{newBlue}{HTML}{000000}
\definecolor{newTeal}{HTML}{000000}%

\usepackage{natbib}
 \bibpunct[, ]{(}{)}{,}{a}{}{,}%

\usepackage{amsthm}
\theoremstyle{plain}
\newtheorem{theorem}{Theorem}
\newtheorem{lemma}{Lemma}
\newtheorem{definition}{Definition}
\newtheorem{corollary}{Corollary}
\newtheorem{proposition}{Proposition}
\newtheorem{remark}{Remark}
\newtheorem{assumption}{Assumption}

\usepackage{etoc}

\newcommand{\dd}{\, \mathrm{d}}
\newcommand{\R}{\mathbb{R}}

\newcommand{\eps}{\varepsilon}

\newcommand{\defeq}{\coloneqq}
\newcommand{\BW}{\mathsf{BW}}
\newcommand{\MF}{\mathsf{MF}}
\newcommand{\Rad}{\mathsf{Rad}}
\newcommand{\PSD}{\mathbb{S}^d_+}

\newcommand{\Sym}{\mathbb{S}^d}
\newcommand{\Orth}{\mathbb{O}^d}
\newcommand{\id}{\mathrm{Id}}
\newcommand{\Id}{\mathrm{Id}}

\newcommand{\cA}{\mathcal{A}}
\newcommand{\cB}{\mathcal{B}}

\newcommand{\cF}{\mathcal{F}}
\newcommand{\cG}{\mathcal{G}}
\newcommand{\cH}{\mathcal{H}}

\newcommand{\cM}{\mathcal{M}}
\newcommand{\cN}{\mathcal{N}}
\newcommand{\cO}{\mathcal{O}}
\newcommand{\cP}{\mathcal{P}}

\newcommand{\cR}{\mathcal{R}}
\newcommand{\cS}{\mathcal{S}}
\newcommand{\cT}{\mathcal{T}}

\newcommand{\EE}{\mathbb{E}}

\newcommand{\HH}{\mathbb{H}}

\newcommand{\LL}{\mathbb{L}}
\newcommand{\MM}{\mathbb{M}}
\newcommand{\NN}{\mathbb{N}}

\newcommand{\PP}{\mathbb{P}}

\newcommand{\RR}{\mathbb{R}}
\renewcommand{\SS}{\mathbb{S}}
\newcommand{\TT}{\mathbb{T}}

\newcommand{\sQ}{\mathsf{Q}}

\let\originalleft\left
\let\originalright\right
\renewcommand{\left}{\mathopen{}\mathclose\bgroup\originalleft}
\renewcommand{\right}{\aftergroup\egroup\originalright}

\title{An Accelerated Proximal Gradient Method for Wasserstein Space}

\usepackage{authblk}

\author[1]{Jakob Nylöf\thanks{Correspondence:
\href{mailto:jakob.nylof@epfl.ch}{\texttt{jakob.nylof@epfl.ch}}.}}
\author[1]{Sloan Nietert}
\author[1]{Julien Pallage}
\author[1]{\authorcr Giancarlo Ferrari-Trecate}
\author[1]{Daniel Kuhn}
\author[2]{John Lygeros}

\affil[1]{EPFL, Lausanne, Switzerland}
\affil[2]{ETH Zürich, Zürich, Switzerland}

\date{}

\begin{document}

\maketitle

\begin{abstract}
We study composite optimization over the 2-Wasserstein space of probability measures. 
We identify conditions under which the objective lifts to a Hilbert space of transport maps as the sum of a smooth convex functional and a convex proximable functional.
This leads us to \texttt{W-FISTA}, an accelerated proximal gradient method that achieves an objective gap of $O(k^{-2})$ after $k$ iterations and accelerated linear convergence when the smooth lift is strongly convex.
Under typical assumptions, its deterministic gradient oracle complexity is minimax optimal up to constant factors over the general composite class. 
We also derive stochastic convergence bounds under a multiplicative noise model, showing that acceleration governs the initial convergence regime, while noise dominates in later iterations. 
Our principal application is variational inference with a log-concave, log-smooth target over the classes of Gaussian, product, and radially symmetric measures. 
We handle negative entropy through tractable proximal steps and obtain deterministic outer iteration bounds with square-root dependence on the target condition number in the strongly log-concave case.
Synthetic experiments for Gaussian variational inference illustrate the predicted improvements.
\end{abstract}

\etocdepthtag.toc{mainmatter}
\etocsettagdepth{mainmatter}{subsection}
\etocsettagdepth{appendix}{none}

\section{Introduction}
\label{sec:intro}

Optimization over probability measures is a prominent paradigm in modern machine learning and statistics, with applications including generative modeling, distributionally robust optimization, and variational inference \citep{Arjovsky_2017_Wasserstein, Esfahani_2017_DRO,Blei_2017_Variational}.
In this work, we study composite problems of the form
\begin{equation}
\label{eq:composite_problem}
\inf_{\mu \in \sQ} \Bigl[ \cF(\mu) \defeq \cG(\mu) + \cR(\mu)\Bigr],
\end{equation}
where $\sQ$ is a subset of the 2-Wasserstein space $\cP_2(\R^d)$ of probability measures with finite second moment.
We consider the full space $\cP_2(\R^d)$, the Bures--Wasserstein (BW) family of Gaussian measures, the mean-field family of product
measures, and the family of radially symmetric measures. Mirroring the typical assumptions for Euclidean proximal gradient methods, we require that $\cG$ satisfies a suitable notion of convexity and smoothness, and that $\cR$ is appropriately convex and proximable (as defined shortly). Our main contribution is identifying an interesting class of problems of the form \eqref{eq:composite_problem} that satisfy these assumptions and can thus be solved efficiently at accelerated rates, via a black-box application of proximal gradient methods for Hilbert spaces.

Our principal application is variational inference (VI),
which approximates a target distribution $\pi$ with density
proportional to $\exp(-V)$, for a smooth potential
$V:\R^d\to\R$, by solving
\begin{equation}
\label{eq:VI}
\inf_{\mu \in \sQ} \KL(\mu \,\|\, \pi).
\end{equation}
Up to an additive constant independent of $\mu$, this is \eqref{eq:composite_problem} with potential energy $\cG(\mu) = \int V\dd\mu$ and negative entropy $\cR = -\Ent$. Specifically, $\Ent(\mu)\defeq-\int \log \mu(x) \dd \mu(x)$ is the differential entropy, identifying $\mu$ with its Lebesgue density. Under the common assumption that $V$ is convex and smooth, this problem falls under our framework. For the BW, mean-field, and radial families, we refer to \eqref{eq:VI} as Gaussian, mean-field, and radial VI, respectively.

We represent each $\mu \in \sQ$ as the pushforward of a reference measure $\rho \in \cP_2(\R^d)$ through a map $T \in \cT$, where $\cT$ is a closed convex subset of $L^2(\rho;\R^d)$, and we maintain momentum iterates lying in a closed linear subspace $\cH$ containing $\cT$. We choose a lifted objective $g:\cH\to\R$ satisfying
$g(T)=\cG(T_\sharp\rho)$ on $\cT$, and define
$r(T) \defeq \cR(T_\sharp\rho)$ on $\cT$ and $r(T)=+\infty$
otherwise.
Then \eqref{eq:composite_problem} is equivalent to
minimizing $g+r$ over $\cH$.
Our accelerated proximal gradient methods apply when:
\begin{enumerate}[label=(\roman*)]
    \item $g$ is $\ell$-strongly convex and $L$-smooth, with $0\le\ell\le L < \infty$, and $\nabla g$ can be evaluated on $\cH$;
    \item $r$ is proper, convex, and lower semicontinuous, and its proximal operator can be evaluated on $\cH$.
\end{enumerate}

Under these conditions, our algorithm \texttt{W-FISTA} applies a refinement of the fast iterative shrinkage-thresholding algorithm (FISTA) of \cite{beck2009fast} (namely, the strongly convex variant of \citealp{chambolle2016introduction}) to the composite objective $g + r$, and it returns $k$-th iterate $\mu_k=(T_k)_\sharp \rho$ with objective gap
\begin{equation*}
    \cF(\mu_k) - \cF(\mu^\star) \lesssim \frac{L R^2}{(k+1)^2}
\exp\left(-\frac{k+1}{2}\sqrt{\frac{\ell}{L}}\,\right).
\end{equation*}
Here, $\mu^\star=(T^\star)_\sharp\rho$ is a minimizer of $\cF$ and $R = \|T^\star-T_0\|_{L^2(\rho)}$ is the initial distance, equal to the 2-Wasserstein distance $W_2(\mu^\star,\mu_0)$ in many of our settings.
When $L = R = O(1)$, \texttt{W-FISTA} reaches objective gap $\eps$ within $O(\eps^{-1/2})$ iterations and, for the strongly convex regime where $\ell>0$ and $\eps \lesssim \ell$, within $O(\sqrt{L/\ell}\,\log(\ell/\eps))$ iterations.
Each iteration requires one evaluation of $\nabla g$, which for the canonical lift $g(T) = \cG(T_\sharp \rho)$ can be computed using a single Wasserstein gradient of $\cG$. For the full $\cP_2(\R^d)$, BW, and mean-field settings, the resulting oracle complexities are minimax optimal for general composite objectives up to absolute multiplicative constants (comparing to deterministic methods in sufficiently high dimension).
\texttt{W-FISTA} also supports stochastic gradients, which we analyze under a suitable multiplicative noise model. This version enjoys accelerated convergence up until an accuracy threshold depending on the noise level.

Further, we establish the required lift properties for concrete
Wasserstein objectives. When $\cG$ is a potential or interaction
energy, we show how convexity and smoothness of the
potential function transfer to the lift $g$, express $\nabla g$ using projections of ambient Wasserstein gradients, and
construct unbiased stochastic gradient fields with
explicit variance bounds. When $\cT$ consists of optimal
transport maps from $\rho$, we show that $r$ is convex
exactly when $\cR$ is convex along generalized geodesics
with pivot $\rho$, a property that holds for negative entropy.

For VI, we derive exact entropy proximal steps for the Gaussian case and develop finite-dimensional formulations for
the mean-field and radial families, with guarantees on
their approximation error. Treating entropy proximally
keeps the smoothness constant governing the outer iterations
equal to that of $V$, independently of the resolution
of the approximation. Our approach improves existing outer iteration bounds for VI,
achieving square-root dependence on the target condition number
and logarithmic dependence on inverse accuracy for strongly
log-concave targets.

After discussing related work, Section~\ref{sec:prelims} introduces the necessary Wasserstein
preliminaries. Section~\ref{sec:W-FISTA} develops our
transport map framework, establishes deterministic and
stochastic convergence guarantees for \texttt{W-FISTA},
and applies them to VI. Section~\ref{sec:tractable} gives
tractable specializations for the Gaussian, mean-field, and
radial families, and Section~\ref{sec:experiments} presents
synthetic experiments illustrating our improved rates for Gaussian VI. Proofs are collected in
the appendices.

\subsection{Related work}

Continuous-time gradient flows and discrete-time optimization methods are both well-studied over Wasserstein and BW space \citep{Jordan_1998_JKO, malago2018wasserstein,han2021riemannian}. For mean-field VI, \cite{du2026particle} analyze a particle discretization of the Wasserstein gradient flow over product measures. In the BW setting, \cite{Chewi_2020_Bures,altschuler2021averaging} develop fast first-order methods for computing barycenters.
\cite{lambert2022variational} study the Gaussian VI problem through the lens of BW optimization, taking explicit gradient steps on the nonsmooth entropy term. \cite{diao2023forward}, a direct inspiration for this work, obtain faster rates by efficiently solving a BW proximal step between gradient updates, without acceleration. 
\citet{domke2020provable,Kim_2023_VI} formulate Gaussian VI as a composite Euclidean problem in mean and Cholesky parameters, admitting an inexpensive entropy
proximal step. \cite{salim2020wasserstein} analyze proximal gradient methods over full $W_2$ space, and \cite{yao2024wasserstein} develop block-wise proximal methods over products of Wasserstein spaces, including mean-field VI. Geodesic extrapolation has also been used to construct
a variational BDF2 scheme for Wasserstein gradient flows
\citep{gallouet2024geodesic}.

Finite-dimensional transport approximations provide another approach to VI.
\citet{jiang2025algorithms} develop accelerated projected gradient methods for mean-field VI over polyhedral families of transport maps, while \citet{ghafourpour2026variational} develop projected gradient methods for radial VI.
Both differentiate the full KL objective over finite-dimensional transport families (approximating infinite-dimensional classes) and impose positive lower bounds on the derivatives of the admissible maps. Their smoothness bounds include an entropy contribution that depends on the resolution of the transport approximation, absent from our deterministic outer iteration bounds.

Accelerated methods over probability measures have been
developed in continuous time through accelerated and
Hamiltonian flows
\citep{taghvaei2019accelerated,wang2022accelerated,
chen2025accelerating}. In discrete time,
\citet{tanaka2023accelerated,tang2026nesterov} develop
Nesterov-type methods using transport map or random variable
lifts. 
Building on this lifted perspective, we combine acceleration with proximal steps to obtain convergence guarantees for composite Wasserstein objectives, including VI with negative entropy handled proximally.

Relatedly, Wasserstein space carries a formal Riemannian structure \citep{otto2001geometry} which suggests using Riemannian accelerated methods. Such methods exist, including proximal variants \citep{feng2025riemannian}, but their guarantees are local, only eventually accelerated, or tied to curvature assumptions \citep{zhang2018estimate,ahn2020nesterov,alimisis2021momentum,kim2022accelerated,martinezrubio2022global,martinezrubio2023accelerated}. These assumptions bound the sectional curvature from above and below on a bounded region assumed to contain all iterates, or require it to be nonpositive or constant.
In some negatively curved geometries, oracle lower bounds rule out Euclidean accelerated rates \citep{hamilton2021nogo,criscitiello2022negative}.
The available curvature-dependent Riemannian guarantees do not directly extend to the full Wasserstein space. Its Alexandrov curvature is nonnegative \citep{ambrosio2008gradient}, but sectional curvature is unbounded above even on the Gaussian family when $d\ge2$ \citep{takatsu2011wasserstein,massart2019curvature},
and geodesics need not admit extensions beyond their endpoints \citep{gallouet2024geodesic}.
Riemannian accelerated methods have also been adapted
to particle-based VI \citep{liu2019understanding}, and dynamic step-size schedulers, such as silver step sizes, have been used to partially accelerate gradient descent for smooth objectives \citep{park2025acceleration}. Our lift avoids curvature assumptions by working in a Hilbert space of transport maps, where global convergence
follows by the stated convexity and smoothness assumptions on the lifted objective.

\textbf{Concurrent work.}
A preliminary version of this work \citep{nylof2026accelerated} focused on accelerated composite optimization over
Bures--Wasserstein space, corresponding to
Section~\ref{sec:BW-FISTA} of the present version.
That version was
submitted on September 2, 2026, and accepted on September 30, 2026, to the Geometric Distributional Deep Learning (GDDL) workshop at NeurIPS 2026.\footnote{\texttt{https://neurips.cc/virtual/2026/loc/paris/workshop/137622}} Moreover, we observed in \citep[Section~6]{nylof2026accelerated} that the approach therein extends directly to 2-Wasserstein space, that is, the setting of the present work. By that time, we had already conceived the main ideas of the present work. That includes developing variants of \texttt{W-FISTA} with deterministic and stochastic gradient oracles, viewing \texttt{BW-FISTA} as a specialization of \texttt{W-FISTA} to the
Bures--Wasserstein space (Section~\ref{sec:BW-FISTA}), and
exploiting finite-dimensional approximations of
mean-field transport maps to obtain numerically
tractable proximal steps (Section~\ref{sec:MF-FISTA}).

The following concurrent works appeared on arXiv after the results of Section~\ref{sec:BW-FISTA} had been accepted to the GDDL workshop and during the development of our extended manuscript.
\begin{itemize}
\item \citet{park2026acceleration} study accelerated
proximal gradient methods based on affine transport
map lifts for Gaussian VI with strongly log-concave
targets, overlapping with our results for $\texttt{BW-FISTA}$ in Section~\ref{sec:BW-FISTA}.
They further investigate the limitations of acceleration
in the stochastic setting.
Our framework additionally accommodates merely
log-concave targets and more general composite
objectives, including nonsmooth components other
than the negative entropy.
    \item \citet{ushiyama2026wasserstein} 
share the underlying transport map lifting principle and
establish similar deterministic convergence
guarantees. Beyond the deterministic setting, our work establishes
convergence guarantees for stochastic variants of
\texttt{W-FISTA}. Moreover, whereas their numerical
approach approximates the proximal steps using
input-convex neural networks, we exploit the structure
of specific distribution families to obtain numerically tractable convex
formulations of the proximal steps. These can be solved exactly as in \texttt{BW-FISTA} or approximated efficiently as in \texttt{MF-FISTA} and \texttt{Rad-FISTA} for VI. 
\end{itemize}

\subsection{AI Usage Statement}

During the preparation of this manuscript, the authors used ChatGPT (OpenAI) to explore and refine their proposed approach of applying composite optimization methods to Wasserstein optimization problems through transport map lifts. They also used ChatGPT to discuss mathematical arguments, investigate alternative proof strategies, identify relevant literature, and improve the clarity and presentation of the manuscript. Furthermore, the authors used Claude (Anthropic) to assist in extending and validating the numerical experiments. The authors take full responsibility for the mathematical correctness, originality, and content of the work.

\section{Preliminaries on Wasserstein Space}
\label{sec:prelims}

We first recall some basics of optimal transport, identify important subfamilies of Wasserstein space, and establish some definitions and results for Wasserstein differentiation.

\textbf{Notation.}
Let $\NN$ and $\NN_+$ denote the sets of nonnegative and positive integers, respectively, and $[d] \coloneq \{1,\dots,d\}$. We write $\langle \cdot,\cdot \rangle$ and $\|\cdot\|$ for the Euclidean inner product and norm, respectively.
For an arbitrary set $S$, let $\iota_S$ denote the indicator function with $\iota_S(x) = 0$ for $x \in S$ and $+\infty$ otherwise.
We denote by $\mathbb{S}^d$, $\mathbb{S}^d_{+}$ and $\mathbb{S}^d_{++}$ the sets
of real symmetric, symmetric positive semidefinite and symmetric positive
definite $d\times d$ matrices, respectively, by $I_d$ the identity matrix, by
$\mathrm{Id}$ the identity map on $\mathbb{R}^d$, by $\|\cdot\|_{\mathrm{F}}$
the Frobenius norm, and by $\tr(\cdot)$ the trace. 
For $x \in \mathbb{R}^d$, let $\mathrm{diag}(x)$ denote the diagonal matrix with diagonal $x$. For $A \in \mathbb{S}^{d}$, denote by $\lambda(A)\in \mathbb{R}^d$ the vector of eigenvalues of $A$ in descending order, and $\lambda_{\max}(A)$ the largest eigenvalue.  For
$A\in\mathbb{S}^d_{+}$, $A^{1/2}\in\mathbb{S}^d_{+}$ is the unique positive
semidefinite square root. We also write
$\Orth\defeq\{U\in\R^{d\times d}:U^\top U=I_d\}$ for the group of orthogonal matrices,  $\Unif(\Orth)$ for its normalized Haar measure, and $\cS^{d-1}$ for the unit sphere in $\R^d$.

Write $\cP_2(\mathbb{R}^d)$ for the set of Borel probability measures on $\mathbb{R}^d$ with finite second moment and $\cP_{2,\mathrm{ac}}(\R^d)$ for those which are absolutely continuous with respect to the Lebesgue measure. Let $\mathcal{N}(m,\Sigma)$ denote the Gaussian with mean $m\in\mathbb{R}^d$ and covariance $\Sigma\in\mathbb{S}^d_{+}$. For general $\mu \in \mathcal{P}_2(\mathbb{R}^d)$, denote by $m_\mu$ and $\Sigma_\mu$ its mean and covariance matrix, respectively. For Borel $T:\mathbb{R}^d\to\mathbb{R}^d$ and 
$\mu \in \mathcal{P}_2(\mathbb{R}^d)$, the pushforward measure $T_\sharp\mu$ is defined by
$(T_\sharp\mu)(\cdot)=\mu(T^{-1}(\cdot))$.
Write $\E_\mu[T] = \E_{X \sim \mu}[T(X)] \defeq \int T(x) \dd \mu(x)$. 
If $\cT$ is a set of Borel maps, then $\cT_\sharp\rho \defeq \{ T_\sharp \rho : T \in \cT \}$. Write $C^k(\mathbb{R}^d)$ for the set of $k$-times continuously differentiable real-valued functions on $\mathbb{R}^d$ and $C^\infty_c(\R^d)$ for those which are  infinitely differentiable and compactly supported.

For a real-valued function $f$ defined on an open subset of a Hilbert space $\cH$, its (Fr\'echet) derivative $Df(x)$ is the bounded linear functional satisfying $f(x+h) = f(x)+Df(x)[h]+o(\|h\|_{\mathcal H})$ as $h \to 0$, unique if it exists. Its gradient $\nabla f(x)\in\mathcal H$ is characterized by $Df(x)[h]=\langle\nabla f(x),h\rangle_{\mathcal H}$. When $\cH = \RR^d$, we write $\partial_i f(x) \coloneq (\nabla f(x))_i$, for all $i \in [d]$. We say that $f$ is $L$-smooth if it is continuously differentiable and $\nabla f$ is $L$-Lipschitz. A function $\varphi:\mathcal H \to (-\infty,\infty]$ is
proper if 
$\operatorname{dom}\varphi
\defeq\{x:\varphi(x)<+\infty\}$ is nonempty.
For proper convex $\varphi$, its subdifferential is given by
\[
    \partial\varphi(x)
    \defeq
    \left\{
        v\in\mathcal H:
        \varphi(y)\ge
        \varphi(x)+\langle v,y-x\rangle_{\mathcal H}
        \text{ for every }y\in\mathcal H
    \right\},
    \qquad x\in\operatorname{dom}\varphi,
\]
and $\partial\varphi(x)=\emptyset$ otherwise.

We write $a\lesssim b$ if $a\leq Cb$ for an absolute constant $C>0$, and $a\asymp b$ if both $a\lesssim b$
and $b\lesssim a$ hold. We use $O$ for big-$O$ notation and $\widetilde{O}$ when suppressing polylogarithmic factors.

\subsection{Wasserstein Space and Subfamilies of Interest}
\label{ssec:w2-space}
Given $\mu,\nu \in \cP_2(\R^d)$, their \emph{2-Wasserstein distance} is defined by
\begin{equation}
\label{eq:W2}
    W_2(\mu,\nu)^2 \defeq \inf_{\pi \in \Pi(\mu,\nu)} \int_{\R^d \times \R^d} \|x - y\|^2 \dd \pi(x,y),
\end{equation}
where $\Pi(\mu,\nu) \subseteq \cP_2(\R^d \times \R^d)$ is the set of couplings of $\mu$ and $\nu$. We recall the following important facts (see \cite{villani2009optimal} for a comprehensive survey): 
\begin{enumerate}[label=(\roman*), nosep]
    \item $(\cP_2(\R^d),W_2)$ is a metric space, which we often call $W_2$ space;
    \item The infimum in \eqref{eq:W2} is always achieved.  We write $\Pi_{\mathrm{opt}}(\mu,\nu) \subseteq \Pi(\mu,\nu)$ for the non-empty set of optimal couplings, and say that a map $T \in L^2(\mu;\R^d)$ is an \emph{optimal transport map} for $W_2(\mu,\nu)$ if $(\id,T)_\sharp \mu \in \Pi_{\mathrm{opt}}(\mu,\nu)$;
    \item If $\mu \in \cP_{2,\mathrm{ac}}(\R^d)$, then $\Pi_{\mathrm{opt}}(\mu,\nu) = \{(\id,T^{\mu \to \nu})_\sharp \mu\}$ for a unique optimal transport map  $T^{\mu \to \nu}$ belonging to the cyclically monotone family
\begin{equation*}
    \mathcal M_{\mathrm{cm}}(\mu)
    =
    \left\{
    T\in L^2(\mu;\mathbb R^d):
    \begin{array}{l}
    \exists\,\varphi:\mathbb R^d\to\mathbb R\cup\{+\infty\}
    \text{ proper and convex},\\
    T(x)\in\partial\varphi(x)\quad\mu\text{-almost everywhere (a.e.)}
    \end{array}
    \right\}.
\end{equation*}
Likewise, we have $W_2(\mu,T_\sharp \mu) = \|T - \id\|_{L^2(\mu)}$ for each map $T \in \cM_{\mathrm{cm}}(\mu)$.
\end{enumerate}

\textbf{Main distribution families.} We consider
the BW space of Gaussian measures, the
mean-field space of product measures, and the family of
radially symmetric measures about the origin:
\begin{align*}
    \mathsf{BW}
    &\defeq
    \left\{
        \cN(m,\Sigma):
        m\in\R^d,\ \Sigma\in\PSD
    \right\}\\
    \mathsf{MF}
    &\defeq
    \left\{
        \mu_1\otimes\cdots\otimes\mu_d:
        \mu_i\in\cP_2(\R),\ i=1,\ldots,d
    \right\}\\
    \mathsf{Rad}
    &\defeq
    \left\{
        \mu\in\cP_2(\R^d):
        U_\sharp\mu=\mu
        \text{ for every }U\in\Orth
    \right\}.
\end{align*}
We refer to
\cite{bhatia2019bures} for a survey on BW geometry, and we employ the term ``mean-field'' for this context in line with
\cite{jiang2025algorithms}, derived from physics applications.

\textbf{Transport map spaces.} 
In addition to the class of cyclically monotone maps, consider
the spaces of affine, coordinate-wise separable, and radial transport maps defined by
\begin{align*}
    \cM_{\mathrm{aff}}(\mathbb M)
    &\defeq
    \left\{
        x \in \R^d \mapsto Ax+b:
        A\in\mathbb M,\ b\in\R^d
    \right\} \qquad \MM \subseteq \R^{d \times d}\\
    \cM_{\mathrm{sep}}
    &\defeq
    \bigl\{
        x \mapsto  (T_1(x_1),\ldots,T_d(x_d)) \,\big|\, \text{each $T_i:\R \to \R$ Borel} \bigr\} \\
    \cM_{\mathrm{rad}}
    &\defeq
    \bigl\{
        T :\R^d \to \R^d \text{ Borel} \,\big|\,  T \circ U = U \circ T \text{ for every } U \in \Orth \bigr\}.
\end{align*}
Each $T\in\cM_{\mathrm{rad}}$ is uniquely determined by
its signed radial profile
\[
    \Psi_T(r)\defeq\langle T(ru),u\rangle,
    \qquad r\ge0,
\]
for any unit vector $u \in \R^d$. Equivariance implies that $\Psi_T$ is well-defined and $T(ru)=\Psi_T(r)u$.

We prove the following in Appendix~\ref{prf:pushforward-expressivity}. Here, for a set $\cA$ of $L^2$ maps and $\cB$ of Borel maps, we abuse notation and write $\cA \cap \cB$ for the set of maps in $\cA$ admitting a representative in $\cB$.

\begin{proposition}
\label{prop:pushforward-expressivity}
Let $\sQ \in \{\BW, \mathsf{MF}, \mathsf{Rad}\}$ and fix reference measure $\rho \in \sQ \cap \cP_{2,\mathrm{ac}}(\R^d)$. Setting
\begin{align*}
\cM_\sQ \defeq \begin{cases}
    \cM_{\mathrm{aff}}(\Sym), & \sQ = \mathsf{BW}\\    \cM_{\mathrm{sep}}, & \sQ = \mathsf{MF}\\
    \cM_{\mathrm{rad}}, & \sQ = \mathsf{Rad},
\end{cases}
\end{align*}
we have 
\begin{equation*}
    (L^2(\rho;\R^d) \cap \cM_\sQ)_\sharp \rho = \sQ \qquad \text{and} \qquad  T^{\rho \to \mu} \in \cM_{\mathrm{cm}}(\rho) \cap \cM_\sQ, \quad \forall \mu \in \sQ
\end{equation*}
An element of $L^2(\rho;\R^d)$ belongs to $\cM_{\mathrm{cm}}(\rho) \cap \cM_\sQ$ if and only if it admits a Borel representative $T:\R^d \to \R^d$ such that
\begin{itemize}[nosep]
    \item ($\sQ = \BW$) $T \in \cM_{\mathrm{aff}}(\PSD)$;
    \item ($\sQ = \MF$) $T \in \cM_{\mathrm{sep}}$ and each $T_i$ is non-decreasing $\rho_i$-a.e.;
    \item ($\sQ = \Rad$) $T \in \cM_{\mathrm{rad}}$ and $\Psi_T$ is non-negative and non-decreasing a.e.\ under $\|\cdot\|_\sharp \rho$.
\end{itemize}
\end{proposition}
We often take $\rho = \cN(0,I_d)$, since it is a valid reference measure for all of our families.

\subsection{Wasserstein Derivatives}

We now establish some basics of $W_2$ calculus, deferring proofs to Appendices~\ref{prf:main-tangent-spaces}--\ref{prf:interaction-energy-Q-grads}. These definitions and results are standard when $\sQ = \cP_2(\R^d)$ and differentiation occurs at $\mu \in \cP_{2,\mathrm{ac}}(\R^d)$, with many appearing in \cite{ambrosio2008gradient}, for example. However, we implement some small tweaks so that all of our subfamilies and singular measures are covered in a unified way.

\begin{definition}[$W_2$ tangent cone]
For a coupling $\pi \in \Pi(\mu,\nu)$, define the barycentric displacement vector field $v_\pi \in L^2(\mu;\R^d)$ by $v_\pi(x) \defeq \int(y-x) \dd \pi(y|x)$, where $\pi(\cdot |x)$ denotes a disintegration of $\pi$ with respect to its first marginal $\mu$. For a family $\sQ\subseteq\cP_2(\R^d)$ and $\mu,\nu \in \sQ$, let $\Pi_{\mathrm{opt}}^{\sQ}(\mu,\nu)$ denote the set of
couplings $\pi\in\Pi_{\mathrm{opt}}(\mu,\nu)$ such that the \emph{constant-speed $W_2$ geodesic} defined by
\[
    \mu_t \defeq ((1-t)x+ty)_\sharp\pi,\qquad t\in[0,1],
\]
lies entirely in $\sQ$. We then define the \emph{barycentric $W_2$ tangent cone to $\mathsf{Q}$ at $\mu \in \mathsf{Q}$} by
\[
    \Tan_\mu \mathsf{Q} \defeq \overline{\left\{\lambda v_\pi : \lambda \geq 0, \nu \in \mathsf{Q}, \pi \in \Pi_{\mathrm{opt}}^\mathsf{Q}(\mu,\nu)\right\}}^{L^2(\mu)}.
\]
\end{definition}

For the distribution families considered above, tangent cones are always closed linear subspaces.

\begin{proposition}[Main tangent spaces]
\label{prop:main-tangent-spaces}
For $\sQ \in \{ \cP_2(\R^d),\BW, \mathsf{MF}, \mathsf{Rad}\}$ and $\mu \in \sQ$, we have
\begin{align*}
    \Tan_\mu \sQ = \begin{cases}
        \overline{\{ \nabla \varphi : \varphi \in C_c^\infty(\R^d) \}}^{L^2(\mu)}, & \sQ = \cP_2(\R^d)\\
        \cM_{\mathrm{aff}}(\Sym) \cap L^2(\mu;\R^d), & \sQ = \BW\\
        \cM_{\mathrm{sep}} \cap L^2(\mu;\R^d), & \sQ = \mathsf{MF}\\
        \cM_{\mathrm{rad}} \cap L^2(\mu;\R^d), & \sQ = \mathsf{Rad}.
    \end{cases}
\end{align*}
\end{proposition}

We then define $W_2$-gradients as follows.

\begin{definition}[$W_2$-differentiability]
\label{def:W2-diff}
Let $\sQ \subseteq \cP_2(\R^d)$. We say that $\cG:\mathsf{Q} \to \R$ is \emph{$\mathsf{Q}$-differentiable at $\mu \in \mathsf{Q}$} if there exists a gradient vector field $\xi \in \Tan_\mu \mathsf{Q}$ such that, for $\nu \in \mathsf{Q}$ and $\pi \in \Pi_{\mathrm{opt}}^\mathsf{Q}(\mu,\nu)$,
\begin{equation*}
    \cG(\nu) = \cG(\mu) + \langle \xi, v_\pi \rangle_{L^2(\mu)} + o(W_2(\mu,\nu))
\end{equation*}
as $\nu \to \mu$ within $\mathsf{Q}$, uniformly over $\pi$. This gradient is unique if it exists, in which case we call it the \emph{$\sQ$-gradient} of $\cG$ at $\mu$ and denote it by $\nabla_\mathsf{Q} \cG(\mu)$. When $\sQ = \cP_2(\R^d)$, we call $\sQ$-differentiability $W_2$-differentiability and denote the $\cP_2(\R^d)$-gradient, termed the $W_2$-gradient, by $\nabla_{W_2} \cG(\mu)$.
\end{definition}

Interaction energies serve as a primary example of $W_2$-differentiable functions.

\begin{proposition}[Interaction energy gradients]
\label{prop:interaction-gradient}
Fix integer $n\ge1$ and $\cG:\cP_2(\R^d) \to \R$ given by $\cG(\mu) \defeq \int V \dd\mu^{\otimes n}$ for a potential $V \!\in\! C^2((\R^d)^n)$ with bounded Hessian. For $\mu \in \cP_2(\R^d)$, define
\[
    w_\mu(x)
    \defeq
    \sum_{j=1}^n
    \E\left[
        \partial_j V(X_1,\ldots,X_n)\mid X_j=x
    \right] \in \R^d,
    \qquad
    X_1,\ldots,X_n\overset{\mathrm{iid}}{\sim}\mu.
\]
When $n= 1$, note that $\cG(\mu) = \EE_\mu[V]$ and $w_\mu = \nabla V$. Then $\nabla_{W_2}\cG(\mu) = w_\mu$.
\end{proposition}

We refer to $\cG$ above as an \emph{$n$-th order interaction energy}. A first order interaction energy is also called a \emph{potential energy}. We typically compute $\sQ$-gradients as projections of $W_2$-gradients.

\begin{proposition}
\label{prop:Q-grad-projection}
Let $\sQ \subseteq \cP_2(\R^d)$ and $\mu \in \sQ$. Assume that $\Tan_\mu \sQ$ is a closed linear subspace. If $\cG:\cP_2(\R^d)\to\R$ is $W_2$-differentiable at $\mu$, then its restriction $\cG|_{\sQ}$ is $\sQ$-differentiable at $\mu$ and
\begin{equation}
    \nabla_{\sQ}(\cG|_{\sQ})(\mu)
    =
    \proj_{\,\Tan_\mu\sQ}\nabla_{W_2}\cG(\mu),
\end{equation}
where $\proj_{\,\Tan_\mu\sQ}$ denotes the orthogonal projection in $L^2(\mu;\R^d)$ onto $\Tan_\mu\sQ$. When there is no ambiguity, we write $\nabla_{\sQ}\cG(\mu)$ for $\nabla_{\sQ}(\cG|_{\sQ})(\mu)$.
\end{proposition}

We now make these projections explicit for the primary tangent spaces.

\begin{proposition}
\label{prop:primary-projections}
Let $\sQ \in \{\BW, \mathsf{MF}, \mathsf{Rad}\}$ and fix $\mu \in \sQ$. For $\xi \in L^2(\mu;\R^d)$, the following formulas hold,
with all vector-field identities understood
$\mu$-almost everywhere.
\begin{enumerate}[label=(\roman*),nosep]
    \item ($\sQ = \BW$) Setting $a \defeq\int \xi \dd\mu$ and $C\defeq \int \xi(x)(x-m_\mu)^\top \dd\mu(x)$, we have
    \begin{equation*}
        [\proj_{\Tan_\mu \BW}(\xi)](x) = a + A(x - m_\mu),
    \end{equation*}
    where $A\in\Sym$ is any solution of $A\Sigma_\mu+\Sigma_\mu A=C+C^\top$. A symmetric solution always exists, and all such solutions induce the same $L^2(\mu)$ vector field.\vspace{1mm}

    \item ($\sQ = \MF$) For each $i = 1, \dots, d$, we have $[\proj_{\Tan_\mu \mathsf{MF}}(\xi)]_i(x_i)
        =
        \E_\mu[\xi_i(X)\mid X_i=x_i]$.
    \item ($\sQ = \Rad$) We have $[\proj_{\Tan_\mu\mathsf{Rad}}\xi](x)
    =
    \E_{U \sim \Unif(\Orth)}[ U^\top\xi(Ux)]$.
\end{enumerate}
\end{proposition}

Combining Proposition~\ref{prop:primary-projections} with
Proposition~\ref{prop:Q-grad-projection} yields the
following gradients for interaction energies.

\begin{corollary}
\label{cor:interaction-energy-Q-grads}
Under the setting of Proposition~\ref{prop:interaction-gradient}, with $\cG$, $V$ and $w_\mu$ as therein, we have
\begin{align*}
    \nabla_{\mathsf{BW}}\cG(\mu)(x)
    &=
    \E_\mu[w_\mu(X)]
    +\E_\mu[D w_\mu(X)](x-m_\mu),\\
    [\nabla_{\mathsf{MF}}\cG(\mu)(x)]_i
    &=
    \E_\mu\left[w_{\mu,i}(X)\mid X_i=x_i\right],
    \quad i=1,\ldots,d,\\
    \nabla_{\mathsf{Rad}}\cG(\mu)(x)
    &=
    \E_{\Unif(\Orth)}
    \!\left[U^\top w_\mu(Ux)\right]\\
    &=\left\{\begin{aligned}&\EE_{\mathrm{Unif}(\cS^{d-1})}[\langle w_\mu(\|x\|_2 X), X\rangle]\dfrac{x}{\|x\|_2}, &&\text{if }x \neq 0,\\[-3mm]
    &0, &&\text{ otherwise,}\end{aligned}\right.
\end{align*}
assuming $\mu$ belongs to the corresponding family. Here, identities hold $\mu$-a.e.\ and $D w_\mu$ is the Jacobian of $w_\mu$. Note when $n=1$ that $w_\mu = \nabla V$ and $Dw_\mu = \nabla^2 V$.
\end{corollary}

\section{Accelerated Optimization over Transport Maps}
\label{sec:W-FISTA}

Existing approaches to optimization over probability space exploit ambient linear structures for $\cP_2(\R^d)$, including spaces of signed measures in mirror descent \citep{aubin2022mirror} and spaces of transport maps or random variables in lifted formulations \citep{gangbo2019differentiability,bonet2024mirror}. Building on the transport map perspective, we specify conditions under which guarantees for accelerated composite optimization over Hilbert spaces transfer to \eqref{eq:composite_problem}.

Throughout, we fix a domain $\sQ \subseteq \cP_2(\R^d)$ and a composite objective $\cF = \cG + \cR$, where $\cG:\sQ \to \R$ and $\cR:\sQ \to \R \cup \{+\infty\}$. All proofs appear in Appendix~\ref{app:W-FISTA}.

\subsection{The Transport Map Lift}

We adopt the following lifted setup.

\begin{assumption} \label{assm:lift}
There exists a reference measure $\rho \in \cP_{2,\mathrm{ac}}(\R^d)$, a closed linear subspace $\cH \subseteq L^2(\rho;\R^d)$, and a closed convex subset $\cT \subseteq \cH$ such that the following hold:
\begin{enumerate}[label=(\roman*)]
        \item (Transport representability.) 
        $\mathsf Q = \cT_\sharp\rho$ and $\Id \in \cT$;
        \item (Smooth lift.) There exists an $\ell$-strongly convex, $L$-smooth function $g:\cH \to \mathbb{R}$ such that $g(T) = \cG(T_\sharp\rho)$ for all $T \in \cT$, where $0 \leq \ell \leq L < \infty$ and $L > 0$;
        \item (Nonsmooth lift.) The function $r:\cH \to \mathbb{R}\cup\{+\infty\}$ given by $r(T) \coloneq \cR(T_\sharp\rho)$ if $T \in \cT$, and $r(T) \coloneq +\infty$ otherwise, is proper, convex, and lower semicontinuous.
\end{enumerate}
\end{assumption}

Here, the values of $g$ are prescribed by $\cG$ only on $\cT$, whereas $r$ is fully determined by $\cR$.

Under Assumption~\ref{assm:lift}, our optimization problem over probability measures \eqref{eq:composite_problem} is equivalent to the following composite optimization problem over the lifted Hilbert space $\cH$:
\begin{equation}
    \label{eq:lifted_composite_problem}
    \min_{T \in \cH} \ g(T) + r(T).
\end{equation}

For cleaner statements, we sometimes impose additional structure.

\begin{assumption}
\label{assm:extras}
For $\rho$, $\cH$, and $\cT$ as in Assumption~\ref{assm:lift}, the following hold:
\begin{enumerate}[label=(\roman*)]
    \item (Weak geodesic convexity.) For all $\mu,\nu \in \sQ$, we have $\Pi_{\mathrm{opt}}^\sQ(\mu,\nu) \neq \emptyset$;
    \item (Cyclic monotonicity.) Each $T \in \cT$ is cyclically monotone, i.e., $\cT \subseteq \cM_{\mathrm{cm}}(\rho)$;
    \item (Canonical smooth lift.) $\cH_\sharp\rho=\sQ$ and $g(T) = \cG(T_\sharp \rho)$ for all $T \in \cH$.
\end{enumerate}
\end{assumption}

Condition \emph{(i)} holds exactly when each pair of measures in $\sQ$ is connected by a constant-speed $W_2$ geodesic that lies fully in $\sQ$. This holds for all of our running example spaces.

\begin{proposition}
\label{prop:weak-geodesic-convexity}
Each family $\sQ \in \{ \cP_2(\R^d), \BW, \MF, \Rad \}$ is weakly geodesically convex.
\end{proposition}

Condition \emph{(ii)} makes each $T\in\cT$ the optimal
transport map from $\rho$ to $T_\sharp\rho$, so
$\|T-\Id\|_{L^2(\rho)}=W_2(\rho,T_\sharp\rho)$. This lets us express the initial radius in our convergence
bounds using the intrinsic $W_2$ distance. Together, conditions \emph{(i)} and \emph{(iii)} allow us to translate differentiability and gradient evaluations from $\sQ$ to $\cH$.

\begin{proposition}[$\cH$-gradients from $\sQ$-gradients]
\label{prop:gradient_connection}
Suppose Assumption~\ref{assm:lift}(i) and Assumption~\ref{assm:extras}(i,iii) hold. Fix $T \in \cH$ and let $\cG$ be $\mathsf{Q}$-differentiable at $T_\sharp\rho \in \mathsf{Q}$. Then, $g$ is Fréchet differentiable at $T$ with gradient $$\nabla_\mathcal{H} g(T) = \mathrm{proj}_\cH\left[\nabla_\mathsf{Q}\cG(T_\sharp \rho) \circ T\right].$$
\end{proposition}

The Assumption~\ref{assm:extras} conditions hold for our main applications, but relaxing them can be useful.
Dropping condition \emph{(ii)} can simplify computation
(see Remark~\ref{rem:lower-triangular}), while allowing a noncanonical smooth lift can accommodate
objectives whose canonical lifts are not differentiable.

\begin{remark}[Non-canonical smooth lifts]
For $\cG(\mu)=W_2(\mu,\rho)^2$, the canonical lift on
$L^2(\rho;\R^d)$ is not differentiable at constant maps.
Yet, if $\cT\subseteq\cM_{\mathrm{cm}}(\rho)$,
then $g(T)=\|T-\id\|_{L^2(\rho)}^2$ agrees with $\cG(T_\sharp\rho)$ on $\cT$ and is
$2$-strongly convex and $2$-smooth over $L^2(\rho;\R^d)$. We consider additional examples for the BW setting in Section~\ref{sec:BW-FISTA}.
\end{remark}

We now give standard choices for the transport map spaces.

\begin{proposition}[Standard choices for $\cH$ and $\cT$]
\label{prop:standard_choices}
Fix $\sQ \in \{ \cP_2(\R^d),\BW, \mathsf{MF}, \mathsf{Rad}\}$ and $\rho \in \sQ \cap \cP_{2,\mathrm{ac}}(\R^d)$. Take $\cH = L^2(\rho;\R^d)$ if $\sQ = \cP_2(\R^d)$ and otherwise set $\cH = L^2(\rho;\R^d) \cap \cM_\sQ$, for $\cM_\sQ$ as in Proposition~\ref{prop:pushforward-expressivity}. Let $\cT = \cH \cap \cM_{\mathrm{cm}}(\rho)$. Then $\cH$ is a closed linear subspace with $\cH_\sharp\rho=\sQ$, $\cT$ is closed and convex, and Assumptions ~\ref{assm:lift}(i) and ~\ref{assm:extras}(i--ii) hold.
\end{proposition}

We next give conditions ensuring the required convexity
and regularity of the lifted objective components. Call $\cG:\cP_2(\R^d)\to\R$ \emph{geodesically convex} if $t\mapsto\cG(\mu_t)$ is convex along every constant-speed $W_2$ geodesic $\{\mu_t\}_{t\in[0,1]}$. The following result relates this property to convexity of the canonical lift. Here, we say that $\rho \in \cP_2(\RR^d)$ \emph{atomless} if $\rho(\{x\})=0$ for every $x\in\RR^d$.

\begin{proposition}[Convexity of canonical lift]
\label{prop:convexity_canonical_lift_on_W2}
Let $\rho\in\cP_2(\R^d)$,
$\cG:\cP_2(\R^d)\to\R$, and
$g(T)\defeq\cG(T_\sharp\rho)$ on $L^2(\rho;\R^d)$.
If $\rho$ is atomless, convexity of $g$ implies
geodesic convexity of $\cG$.
Conversely, geodesic convexity of $\cG$ implies
convexity of $g$ if $\cG$ is everywhere
$W_2$-differentiable \citep[Theorem~1.1]{parker2024some},
or if $d\geq 2$ and $\cG$ is $W_2$-continuous
\citep[Theorem~5.6]{Cavagnari_2026_Lagrangian}.
\end{proposition}

We often consider interaction energies $\cG(\mu) = \int V \dd\mu^{\otimes n}$, in which case differentiability, convexity and smoothness are inherited from the integrand $V$.

\begin{proposition}[Interaction energy lift]
\label{prop:interaction_energy_lift}
Let $\rho \in \cP_{2,\mathrm{ac}}(\R^d)$ and
$\cH\subseteq L^2(\rho;\RR^d)$ be a closed linear
subspace. Fix an integer $n\geq1$ and consider the
interaction energy $\cG:\cP_2(\RR^d)\to\RR$ defined by $\cG(\mu)\defeq\int V\dd\mu^{\otimes n}$, for an $\ell_V$-strongly convex and $L_V$-smooth
potential $V\in C^2((\RR^d)^n)$, with
$0\leq\ell_V\leq L_V<\infty$ and $L_V>0$.
Then the lift $g:\cH\to\RR$ defined by
$g(T)\defeq\cG(T_\sharp\rho)$ is Fr\'echet differentiable,
$n\ell_V$-strongly convex, and $nL_V$-smooth, with gradient $\nabla_\cH g(T)=\mathrm{proj}_\cH[w_{T_\sharp\rho}\circ T]$, where $w_\mu$ is defined in
Proposition~\ref{prop:interaction-gradient}.
\end{proposition}

Our main choice of the nonsmooth component $\cR$ is
the negative entropy functional $-\mathrm{Ent}$.
The following result characterizes the regularity of
the nonsmooth lift $r$ in terms of the corresponding
properties of $\cR$. For negative entropy, it suffices that $-\Ent$ is finite
at some measure in $\sQ$.

\begin{proposition}[Nonsmooth lift regularity]
\label{prop:nonsmooth_lift_regular}
Suppose that Assumptions~\ref{assm:lift}{(i)} and \ref{assm:extras}(ii) hold.
Then Assumption~\ref{assm:lift}(iii) holds if and only if $\cR$ is proper, $W_2$-lower semicontinuous on $\sQ$, and generalized geodesically convex with pivot $\rho$, i.e., for all $\mu,\nu \in \sQ$, the map $t \mapsto \cR(((1-t)T^{\rho \to \mu} + tT^{\rho \to \nu})_\sharp \rho)$ is convex over $t \in [0,1]$. In particular, this implies that Assumption~\ref{assm:lift}(iii) holds for $\cR = -\Ent\!|_\sQ$ so long as $\sQ$ contains a measure $\mu$ with $-\Ent(\mu) < \infty$.
\end{proposition}

\subsection{Gradient Oracle Models}

Our optimization methods will access $g$ through its gradients or, in the stochastic case, through the following stochastic first-order oracle. Here and throughout, (Borel) measurability conditions are included for completeness, but we encourage the reader to ignore these during an initial read.

\begin{definition}[Stochastic $\cH$-gradient oracle for $g$]
\label{def:oracle-Hilbert}
Suppose Assumption~\ref{assm:lift} holds. Let $(\Xi,\PP)$ be a probability space, set $\sigma \geq 0$, and let $\zeta_\cH:\mathcal{H}\to[0,+\infty)$ be 1-Lipschitz. A \emph{stochastic $\cH$-gradient oracle} for $g$ with noise parameters $\sigma,\zeta_\cH$ is a jointly measurable map $\cO_\cH:\cH\times\Xi \to \cH$ satisfying
\[
\E_{\xi\sim\PP}[\cO_\cH(T,\xi)] = \nabla_{\cH}g(T), \quad \E_{\xi\sim\PP}
\left[\left\|\cO_\cH(T,\xi)-\nabla_{\cH}g(T)\right\|_{L^2(\rho)}^2\right]\leq
\sigma^2 \zeta_\cH(T)^2 \qquad \forall T \in \cH.
\]
\end{definition}

By Proposition~\ref{prop:gradient_connection}, deterministic $\cH$-gradients of $g$ can be simulated using $\sQ$-gradients of $\cG$.
A similar reduction holds for the stochastic case.

\begin{definition}[Stochastic $\sQ$-gradient oracle for $\cG$]
\label{def:oracle-Q}
Suppose that $\cG:\sQ\to\RR$ is $\sQ$-differentiable.
Let $\Xi_\mathsf{Q}$ be a measurable space and let
$\{\PP_\mu\}_{\mu\in\sQ}$ be a family of
probability measures on $\Xi_\mathsf{Q}$.
Set $\sigma\geq0$ and let
$\zeta_\sQ:\sQ\to[0,+\infty)$ be $1$-Lipschitz
with respect to $W_2$.
A \emph{stochastic $\sQ$-gradient oracle} for $\cG$
with noise parameters $\sigma,\zeta_\sQ$ is a mapping
from each $(\mu,\eta)\in\sQ\times\Xi_\mathsf{Q}$ to a vector field
$\cO_\sQ(\mu,\eta)\in\Tan_\mu\sQ$, such that, for every $\mu \in \mathsf{Q}$,
$\cO_\sQ(\mu,\cdot)$ is $\PP_\mu$-measurable and
\[
\E_{\eta\sim\PP_\mu}[\cO_\sQ(\mu,\eta)]
=\nabla_\sQ\cG(\mu),
\quad
\E_{\eta\sim\PP_\mu}
\left[
\|\cO_\sQ(\mu,\eta)-\nabla_\sQ\cG(\mu)\|_{L^2(\mu)}^2
\right]
\leq \sigma^2\zeta_\sQ(\mu)^2.
\]
\end{definition}

\begin{proposition}[Stochastic $\cH$-gradient oracles
from $\sQ$-gradient oracles]
\label{prop:stochastic_gradient_connection}
Suppose that Assumption~\ref{assm:lift} and
Assumption~\ref{assm:extras}(i,iii) hold and that
$\cG$ is $\sQ$-differentiable.
Let $\cO_\sQ$ be a stochastic $\sQ$-gradient oracle
for $\cG$ with the associated family of probability spaces $\{(\Xi_\mathsf{Q}, \PP_\mu)\}_{\mu \in \mathsf{Q}}$ and noise parameters $\sigma,\zeta_\sQ$. Let $(\Xi,\PP)$ be a probability space and let
$\eta:\cH\times\Xi \to \Xi_\mathsf{Q}$ be a map satisfying $\eta(T,\cdot)_\sharp\PP=\PP_{T_\sharp\rho}$, for all $T\in\cH$.
If the map $\cO_\cH: \cH \times \Xi \to \cH$ given by
\[
\cO_\cH(T,\xi)\coloneq
\mathrm{proj}_\cH
\left[
\cO_\sQ(T_\sharp\rho,\eta(T,\xi))\circ T
\right]
\]
is jointly measurable, then it is a stochastic $\cH$-gradient oracle for $g$
with noise parameters $\sigma$ and
$\zeta_\cH(T)\coloneq\zeta_\sQ(T_\sharp\rho)$.
\end{proposition}

We allow the sampling measure $\PP_\mu$ to depend
on the query measure $\mu$, since stochastic $\sQ$-gradient
oracles often draw samples $\eta\sim\mu$.
This is the case for the BW and mean-field oracles introduced
in the following proposition, which gives  stochastic $\sQ$-gradients for interaction energies.

\begin{proposition}
\label{prop:stochastic-Q-examples}
Consider the setting of
Proposition~\ref{prop:interaction-gradient} with $n$, $\cG$, $V$, and $w_\mu$ as therein, where $V$ is further $L_V$-smooth for $L_V \geq 0$. 
Consider the %
probability spaces $\{(\Xi_\mathsf{Q}, \PP_\mu)\}_{\mu \in \mathsf{Q}}$, where
\[\Xi_\mathsf{Q} =
\begin{cases}
\RR^d,
&\sQ\in\{\mathsf{BW},\mathsf{MF}\},\\
\cS^{d-1},
&\sQ=\mathsf{Rad},
\end{cases}\qquad
\PP_\mu =
\begin{cases}
\mu,
&\sQ\in\{\mathsf{BW},\mathsf{MF}\},\\
\operatorname{Unif}(\cS^{d-1}),
&\sQ=\mathsf{Rad}.
\end{cases}
\]
Define the map $ (\mu,\eta) \in \sQ \times \Xi_\sQ \mapsto \cO_\mathsf{Q}(\mu,\eta) \in \mathrm{Tan}_\mu \mathsf{Q}$ by
\[
\cO_\sQ(\mu,\eta)(x)\coloneq
\begin{cases}
w_\mu(\eta)+Dw_\mu(\eta)(x-m_\mu),
&\sQ=\mathsf{BW},\\[1mm]
\bigl(w_{\mu,i}(x_i,\eta_{-i})\bigr)_{i=1}^d,
&\sQ=\mathsf{MF},\\[1mm]
\displaystyle
\left\langle
w_\mu(\|x\|_2\eta)-w_\mu(0),\eta
\right\rangle
\frac{x}{\|x\|_2},
&\sQ=\mathsf{Rad},\quad x\ne0,
\end{cases}
\]
with $\cO_\mathsf{Rad}(\mu,\eta)(0)\coloneq0$.
Here, $(x_i,\eta_{-i})$ denotes the vector obtained
by replacing coordinate $i$ of $\eta$ by $x_i$. Then $\cO_\sQ$ is a stochastic $\sQ$-gradient oracle
for $\cG$, with
\[
\zeta_\sQ(\mu)
=
\left(\EE_\mu\|X-m_\mu\|_2^2\right)^{1/2},
\qquad
\sigma^2=
\begin{cases}
2n^2L_V^2,&\sQ=\mathsf{BW},\\
n^2L_V^2,&\sQ\in\{\mathsf{MF},\mathsf{Rad}\}.
\end{cases}
\]
For $n=1$, recall that $w_\mu=\nabla V$ and
$Dw_\mu=\nabla^2V$.
\end{proposition}

Combining Propositions~\ref{prop:stochastic-Q-examples} and ~\ref{prop:stochastic_gradient_connection}, one obtains corresponding stochastic $\cH$-gradient oracles for interaction energies in each of the cases $\mathsf{Q} \in \{\mathsf{BW}, \mathsf{MF}, \mathsf{Rad}\}$. We use this in Section~\ref{sec:vi} for the potential energy.

\subsection{Wasserstein FISTA: An Accelerated Proximal Gradient Method for Wasserstein Space}

We now introduce Algorithm \ref{alg:W-FISTA}, the Wasserstein Fast Iterative Shrinkage Thresholding Algorithm (\texttt{W-FISTA}). We arrive at \texttt{W-FISTA} by applying an accelerated proximal gradient method for composite optimization in Hilbert spaces to Problem \eqref{eq:lifted_composite_problem}. 
The following theorem gives its deterministic
convergence rate and iteration complexity.

\begin{theorem}[Deterministic \texttt{W-FISTA} guarantee]
\label{thm:w-fista_rate}
Suppose Assumption~\ref{assm:lift} holds. Let $\{T_k\}_{k\geq0}$
be generated by Algorithm~\ref{alg:W-FISTA} with deterministic
gradients, step size $0<h\leq1/L$, and
$q\defeq h\ell\in[0,1]$.
Assume that there exists a minimizer
$\mu^\star\in\argmin_{\mu\in\sQ}\cF(\mu)$, and define the
initial radius
$R\defeq\inf\{\|T^\star-T_0\|_{L^2(\rho)}\!:\!
T^\star\!\in\!\cT, (T^\star)_\sharp\rho\!=\!\mu^\star\}$.
Writing $\mu_k\!\defeq\!(T_k)_\sharp\rho$, we have for each
$k\ge1$ that
\begin{equation}\label{eq:chambolle_bound}
\cF(\mu_k)-\cF(\mu^\star)
\leq
\frac{8R^2}{h(k+1)^2}
\exp\left(-\frac{k+1}{2}\sqrt q\right).
\end{equation}
Hence, for each accuracy $\eps>0$, there exists an iteration count $N_{\mathrm{det}}(\eps;\ell,L,R)$ satisfying
\begin{equation}\label{eq:iteration-complexity-det}
N_{\mathrm{det}}(\eps;\ell,L,R)
\asymp
1+
\begin{cases}
\displaystyle
\sqrt{\frac{LR^2}{\eps}},
& \displaystyle \eps>\frac{\ell R^2}{4},
\\[4mm]
\displaystyle
\sqrt{\frac{L}{\ell}}\,
\log\frac{\ell R^2}{\eps},
& \displaystyle 0<\eps\leq\frac{\ell R^2}{4},
\end{cases}
\end{equation}
such that
$\cF(\mu_k)-\cF(\mu^\star)\leq\eps$ for all
$k\geq N_{\mathrm{det}}(\eps;\ell,L,R)$, taking $h = 1/L$.
Under Assumption~\ref{assm:extras}(ii), we further have $R=W_2(\mu_0,\mu^\star)=W_2(\rho,\mu^\star)$.
\end{theorem}

When $\ell=0$, \texttt{W-FISTA} is the standard FISTA procedure
of \cite{beck2009fast} applied to $g+r$.
For $\ell>0$, it coincides with the strongly convex variant of
\cite{chambolle2016introduction}.
Our proof refines their analysis to obtain
\eqref{eq:chambolle_bound}.
The following proposition shows that these iteration
complexities are optimal up to absolute constants for
deterministic $\sQ$-gradient queries.

\begin{proposition}[$\sQ$-gradient oracle complexity]
\label{prop:Q_gradient_oracle_complexity}
Suppose that $\cG$ is $\sQ$-differentiable, that Assumptions ~\ref{assm:lift} and \ref{assm:extras} hold, and that $W_2(\rho,\mu^\star) \leq R$ for $\mu^\star \in \argmin_{\mu \in \sQ} \cF(\mu)$. Then Proposition~\ref{prop:gradient_connection} and Theorem~\ref{thm:w-fista_rate} imply that \texttt{W-FISTA} can achieve objective gap $\eps > 0$ using $N_{\mathrm{det}}(\eps;\ell,L,R) $ evaluations of $\nabla_\sQ \cG$. Supposing further that $\sQ$ is closed under translations, $R > 0$, $d > 2N_{\mathrm{det}}(\eps;\ell,L,R)$, and $0 \leq \ell \leq L/4$, there exists an absolute constant $c > 0$ such that any deterministic procedure which accesses $\cG$ through only fewer than $c N_{\mathrm{det}}(\eps;\ell,L,R)$ queries to $\nabla_\sQ \cG$ must suffer objective gap at least $\eps$ for a valid instance.
\end{proposition}

The stochastic case requires a customized proof due to our slightly non-standard noise model.

\begin{algorithm}[t]
\caption{W-FISTA}
\label{alg:W-FISTA}
\begin{algorithmic}[1]
\Statex \textbf{Parameters:} Step size $h>0$, momentum tuning constant $q \in [0,1]$
\State Initialize $T_{-1}=T_0=\mathrm{Id}, \ t_0 = 0$
\For{$k=0,1,2,\ldots$}
    \State $t_{k+1}
    =
    \frac{1}{2}(1-qt_k^2+\sqrt{(1-qt_k^2)^2+4t_k^2}), \ \beta_k=\frac{t_{k}-1}{t_{k+1} + q t_k^2}$\label{step:weight-recurrence}
    \State $Y_k = T_k+\beta_k(T_k-T_{k-1})$
    \If{deterministic gradients} \label{line:det}
    \State $Z_k =Y_k - h \nabla_{\cH}g(Y_k)$
    \ElsIf{stochastic gradients} 
    \State $Z_k = Y_k-h\,\cO_{\mathcal{H}}(Y_k,\xi_k)$, where $\xi_k \sim \PP$ is drawn independently 
    \EndIf \label{line:stochastic}
    \State
    $T_{k+1}
    =
    \argmin_{T\in\cT}
    \left\{
    r(T)
    +
    \frac{1}{2h}
    \|T-Z_k\|_{L^2(\rho)}^2
    \right\}
    $ \label{line:proximal_step}
\EndFor
\end{algorithmic}
\end{algorithm}

\begin{theorem}[Stochastic \texttt{W-FISTA} guarantee]
\label{thm:w-fista_rate_stochastic}
Suppose Assumption~\ref{assm:lift} holds, and let $\cO_\cH$
be a stochastic gradient oracle for $g$ with noise parameters
$\sigma$ and $\zeta_\cH$.
Let $\{T_k\}_{k\ge0}$ be generated by Algorithm~\ref{alg:W-FISTA} using $\cO_\cH$, step size $0<h\le1/(2L)$, and $q\defeq h\ell\in[0,1/2]$. Take iterates $\{\mu_k\}_{k \geq 0}$, minimizer $\mu^\star$, and radius $R$ as in Theorem~\ref{thm:w-fista_rate}, and set variance parameter $\bar{\sigma}^2 \defeq 2 \sigma^2 \left[(\zeta_\cH(T_0) + R)^2 + R^2 \right]$. Suppose $R > 0$; otherwise $\mu_0$ is optimal. Then, for $k \geq 1$, we have
\begin{equation}\label{eq:WSfista_rate}
\begin{aligned}
\E\bigl[\cF(\mu_k)-\cF(\mu^\star)\bigr]
&\le
\frac{8R^2}{h(k+1)^2}
\exp\left(-\frac{k+1}{2}\sqrt q\right) 
+2\bar\sigma^2h\min\{k,q^{-1/2}\},
\end{aligned}
\end{equation}
provided that either
$h\le1/(2\sigma k^{3/2})$, or
$\ell>0$ and $h\le\ell^3/(64\sigma^4)$. 

These restrictions are vacuous when $\sigma=0$, and we take $q^{-1/2}=+\infty$ when $q=0$. Consequently, for each accuracy $\eps>0$, there exists an iteration count $N_{\mathrm{stoc}}(\eps;\ell,R,\bar\sigma)$ satisfying
\begin{equation}\label{eq:iteration-complexity-stoc}
N_{\mathrm{stoc}}(\eps;\ell,R,\bar\sigma)
\asymp 1 + 
\begin{cases}
\displaystyle
\frac{\bar\sigma^2R^2}{\eps^2},
& \displaystyle \eps>\frac{\ell R^2}{4},
\\[4mm]
\displaystyle
\frac{\bar\sigma^2}{\ell\eps}
\log\frac{\ell R^2}{\eps},
& \displaystyle 0<\eps\le\frac{\ell R^2}{4},
\end{cases}
\end{equation}
such that $\E[\cF(\mu_k)-\cF(\mu^\star)]\le\eps$ whenever
$k\ge
N_{\mathrm{det}}(\eps;\ell,L,R)
+
N_{\mathrm{stoc}}(\eps;\ell,R,\bar\sigma)$ and the step size $h$ is properly tuned as a function of $k$.
\end{theorem}

\begin{remark}[Noise-dominated regime]
One can show that $N_{\mathrm{det}}$ dominates $N_{\mathrm{stoc}}$ when
\[
\eps\gtrsim \eps_{\mathrm{noise}}
\defeq
\min\left\{
\left(\frac{\bar\sigma^4R^2}{L}\right)^{1/3},
\frac{\bar\sigma^2}{\sqrt{\ell L}}
\right\},
\]
where the second threshold is interpreted as $+\infty$ when
$\ell=0$. Below this threshold, the stochastic contribution dominates the iteration complexity.
\end{remark}

We view our general reduction to Hilbert space optimization, rather than the exact rates and guarantees above, as our primary conceptual contribution. This enables faster rates under additional assumptions, for example if $r$ is also strongly convex as in \cite{chambolle2016introduction}, as well as refined guarantees using more complex algorithms. 

\begin{remark}[Horizon dependence of step size]
In the convex stochastic regime, horizon-dependent step sizes are standard for constant-step Euclidean accelerated methods (see, e.g., Corollary~15 of \citealp{kulunchakov2020estimate}).
This can be avoided by modifying the
algorithm, for example through restarting and decreasing step
sizes in the strongly convex setting
\citep[Corollary~14]{kulunchakov2020estimate}, or through
anytime accelerated constructions in the convex setting
\citep{cutkosky2019anytime}. 
Above, we prioritized a shared simple algorithm at the cost of horizon dependence.
\end{remark}

Finally, we common on an apparent oversight of the smooth setting.

\begin{remark}[Reduction to optimization over $\R^d$ when $\cR=0$]
Our work immediately jumps to the nonsmooth setting without
explicitly treating the case where $\cR=0$ and $\cF=\cG$.
However, suppose that $\sQ$ is closed under translations and
that the conditions of Proposition~\ref{prop:gradient_connection}
hold. If $T^\star\in\cT$ is a minimizing map, set
$m^\star\defeq\E_\rho[T^\star]$ and let
$\bar T^\star\equiv m^\star$. Then, translation closure and transport
representability imply that $\bar T^\star\in\cT$, while
Proposition~\ref{prop:gradient_connection} implies that
$\nabla_\cH g(\bar T^\star)$ is constant. Hence
\[
\langle\nabla_\cH g(\bar T^\star),
T^\star-\bar T^\star\rangle_{L^2(\rho)}=0,
\]
and convexity gives
$g(\bar T^\star)\le g(T^\star)$. Thus $\delta_{m^\star}$ is a minimizer, and $\inf_{\mu\in\sQ}\cG(\mu)
=
\inf_{m\in\R^d}\cG(\delta_m)$. The smooth problem therefore reduces to ordinary convex
optimization of $m\mapsto\cG(\delta_m)$ over $\R^d$.
\end{remark}

\subsection{Application to Variational Inference}
\label{sec:vi}

As a principal application of \texttt{W-FISTA}, we consider the VI problem over $\sQ$,
\begin{equation}
    \label{eq:VI_in_sec}
    \min_{\mu\in\sQ}\KL(\mu\Vert\pi).
\end{equation}
Here, $\pi\in\cP_2(\RR^d)$ is the target distribution whose Lebesgue density is proportional to $\exp(-V)$, for an $\ell_V$-strongly convex and
$L_V$-smooth potential $V:\R^d \to \R$, and the Kullback-Leibler (KL) divergence is defined as
$$
\mathrm{KL}(\mu\Vert\pi)
\coloneq
\begin{cases}
\displaystyle\int \log\left(\frac{\mathrm d\mu}{\mathrm d\pi}\right)\mathrm d\mu,
& \mu\ll\pi,\\[3mm]
+\infty,&\text{otherwise},
\end{cases}
$$
where $\frac{\dd\mu}{\dd\pi}$ denotes the Radon–Nikodym derivative when $\mu$ is absolutely continuous with respect to $\pi$, denoted by $\mu\ll\pi$. Solving \eqref{eq:VI_in_sec} finds the closest distribution $\mu \in \sQ$ to $\pi$ in terms of the KL divergence. VI restricts the search to a prescribed family \(\sQ \subseteq \cP_2(\RR^d)\), as the unique minimizer over \(\cP_2(\RR^d)\) is \(\pi\). %

The KL divergence decomposes into the potential energy and the negative entropy up to an additive constant (see, e.g., \citealp{salim2020wasserstein}),
\[
    \KL(\mu\Vert\pi)
    =
    \cG(\mu) + \cR(\mu)
    + \mathrm{const}., \quad \cG(\mu) = \int V\dd\mu, \quad \cR(\mu) = -\mathrm{Ent}(\mu).
\]

For each $\sQ\in\{\BW,\MF,\Rad\}$, representing Gaussian, mean-field and Radial VI, respectively, consider the standard
choices of $\cT$ and $\cH$ in
Proposition~\ref{prop:standard_choices} and define the lifts
\begin{equation}
\label{eq:VI_lifts}
\begin{aligned}
    &g:\cH\to\RR,\quad
    g(T)=\cG(T_\sharp\rho)=\int V(T(x))\dd\rho(x),\\
    &r:\cH\to\RR\cup\{+\infty\},\quad
    r(T)=-\Ent(T_\sharp\rho)+\iota_\cT(T).
\end{aligned}
\end{equation}
Each choice of $\sQ$ yields one specialization of \texttt{W-FISTA} for the VI problem in \eqref{eq:VI_in_sec}. Table~\ref{tab:vi_summary} lists the
deterministic gradients and stochastic gradient oracles corresponding to each specialization. The following proposition establishes that
Assumptions~\ref{assm:lift} and~\ref{assm:extras}
hold in all cases and justifies the entries in the table.

\begin{proposition}[Regularity and gradients for VI]
\label{prop:vi_summary}
Fix $\rho=\cN(0,I_d)$. For $\sQ\in\{\BW,\MF,\Rad\}$,
let $\cH$ and $\cT$ be as in
Proposition~\ref{prop:standard_choices}, and let $g$ and $r$
be the lifts in \eqref{eq:VI_lifts}.
Suppose that $V \in C^2(\RR^d)$ is $\ell_V$-strongly convex and
$L_V$-smooth. Then, Assumptions~\ref{assm:lift} and \ref{assm:extras} hold with $\ell = \ell_V$ and $L = L_V$. Moreover, for every $T\in\cH$, the $\cH$-gradient
$\nabla_\cH g(T)$ is given in Table~\ref{tab:vi_summary},
and the maps $\cO_\cH$ therein are stochastic
$\cH$-gradient oracles in the sense of
Definition~\ref{def:oracle-Hilbert}, where
$\xi\sim\rho$ for $\sQ\in\{\BW,\MF\}$ and
$\xi\sim\Unif(\cS^{d-1})$ for $\sQ=\Rad$,
with noise parameters
\[
    \sigma^2=
    \begin{cases}
        2L_V^2,&\sQ=\BW,\\
        L_V^2,&\sQ\in\{\MF,\Rad\},
    \end{cases}
    \quad
    \zeta_\cH(T)=
    \|T-\EE_\rho[T(X)]\|_{L^2(\rho)}.
\]
\end{proposition}

Proposition~\ref{prop:vi_summary} implies that the convergence rates of
Theorems~\ref{thm:w-fista_rate}
and~\ref{thm:w-fista_rate_stochastic} apply to VI.

\begin{table}[t]
\centering
\small
\caption{Deterministic gradients and stochastic gradient
oracles for the VI problem in \eqref{eq:VI_in_sec} for $\sQ \in \{\BW, \MF, \Rad\}$, given that $\rho = \cN(0,I_d)$, the spaces $\cH$ and $\cT$ are as in Proposition~\ref{prop:standard_choices}, and the lifts $g$ and $r$ are as in \eqref{eq:VI_lifts}. For a square matrix $M$, we write $\mathrm{sym}(M) \defeq (M + M^\top)/2$. The radial gradient formulas hold for $T(x) \neq 0$ and are defined to be zero when $T(x) = 0$.}
\label{tab:vi_summary}
\begin{tabular}{lll}
\toprule
$\sQ$ & $\nabla_\cH g(T)(x)$ & $\cO_\cH(T,\xi)(x)$ \\
\midrule
$\BW$
& $\EE_\rho[\nabla V(T(X))]+\mathrm{sym}(\EE_\rho[\nabla^2V(T(X))]A)x$
& $\nabla V(T(\xi))+\mathrm{sym}(\nabla^2V(T(\xi))A)x$
\\[1ex]
$\MF$
& $\big(\EE_{\rho}[\partial_i V(T_i(x_i),T_{-i}(X_{-i}))]\big)_{i=1}^d$
& $\big(\partial_i V(T_i(x_i),T_{-i}(\xi_{-i}))\big)_{i=1}^d$
\\[1ex]
$\Rad$
& $\EE_{\mathrm{Unif(\cS^{d-1})}}\left[
    \left\langle\nabla V(\|T(x)\|_2X),X\right\rangle
    \right]\frac{T(x)}{\|T(x)\|_2}$
& $\left\langle
    \nabla V(\|T(x)\|_2\xi)-\nabla V(0),\xi
    \right\rangle\frac{T(x)}{\|T(x)\|_2}$
\\
\bottomrule
\end{tabular}
\end{table}

\section{Numerically Tractable Specializations of W-FISTA}
\label{sec:tractable}

Algorithm~\ref{alg:W-FISTA} leaves the implementation of
the proximal step at line~\ref{line:proximal_step}
unspecified. Thus, the iteration bounds in
Section~\ref{sec:W-FISTA} do not account for the cost of
solving a potentially infinite-dimensional optimization
problem over $\cT\subseteq L^2(\rho;\RR^d)$ at each iteration.

We focus on two settings in which the proximal operator
in \texttt{W-FISTA} can be computed or approximated by solving a numerically
tractable convex optimization problem.
In the first setting, $\cH$ is finite-dimensional,
allowing the proximal problem to be formulated in terms
of the coefficients of a basis expansion.
In the second setting, $\cH$ is infinite-dimensional
but admits a low-dimensional structure suitable for finite-dimensional approximations. The same coefficient representation can then be used
to solve the resulting finite-dimensional problems.

    \subsection{Bures--Wasserstein FISTA}
    \label{sec:BW-FISTA}

    We first specialize \texttt{W-FISTA} to the Bures--Wasserstein space $\BW$ of Gaussian measures, which admits a natural finite-dimensional parametrization through affine transport maps. This setting illustrates the more general finite-dimensional formulation developed in
    Section~\ref{sec:finite-W-FISTA}. Moreover, this setting admits proximal steps that reduce to numerically tractable convex programs.

    Specifically, we refer to \texttt{BW-FISTA} as Algorithm~\ref{alg:W-FISTA} with $\sQ=\BW$ and\footnote{
    The construction extends to other nondegenerate
    Gaussian reference measures $\rho$. We choose
    $\rho=\cN(0,I_d)$ because it yields the simple
    inner product
    \(
        \langle T_{b,A},T_{d,C}\rangle_{L^2(\rho)}
        = b^\top d+\tr(A^\top C).
    \)}
    \begin{align}
    \label{eq:BW-choices}
    &\rho = \cN(0,I_d),\quad  \cT = L^2(\rho;\RR^d) \cap \cM_{\mathrm{aff}}(\TT), \quad \cH = L^2(\rho;\RR^d) \cap \cM_{\mathrm{aff}}(\HH) ,\end{align}
    where $\TT$ is a closed convex subset of $\RR^{d\times d}$ and $\HH$ is a linear subspace of $\RR^{d\times d}$ containing $\TT$. In other words, we consider a more general affine parametrization than the standard choice in Proposition~\ref{prop:standard_choices}, allowing for transport maps $\cT$ that are not necessarily cyclically monotone. The subspace condition on $\HH$ ensures that $\cH$ is a Hilbert space. Assumption~\ref{assm:lift}\emph{(i)} holds if and only if $I_d \in \TT$ and the map $\mathbb{T} \ni A \mapsto AA^\top \in \mathbb{S}_+^d$ is surjective. Consequently, any $\HH$ such that $\TT \subseteq \HH \subseteq \RR^{d \times d}$ makes $\cH_\sharp\rho = \BW$.
    
    For example, $\TT$ can be taken as $\SS^d_+$ (by the matrix square root) or the set of lower triangular matrices with non-negative diagonals $\LL^d_+$ (by the Cholesky factorization). In either case, a natural choice for the Hilbert space $\HH$ is the linear span of $\TT$, that is, $\SS^d$ or the set of lower triangular matrices $\LL^d$, respectively. If $\TT = \SS_+^d$, then $\cT$ is the set of optimal transport maps from $\rho$ to any Gaussian distribution, while $\TT = \LL_+^d$ includes suboptimal transport maps.

    \textbf{Affine coefficient representation.} \texttt{BW-FISTA} is a special case of finite-dimensional \texttt{W-FISTA}, discussed in Section~\ref{sec:finite-W-FISTA}, as $\cH$ is finite-dimensional. However, rather than introducing an explicit basis for $\cH$, we exploit the isometric isomorphism
    $$\begin{aligned}&\RR^d \times \HH \to \cM_\mathrm{aff}(\HH), \ (b,A) \mapsto T_{b,A},\end{aligned}$$
    where $T_{b,A}$ denotes a map in $\cM_\mathrm{aff}(\HH)$ with affine coefficients $b$ and $A$, and $\RR^d\times\HH$ is equipped with the inner product $\langle (b,A), (d,C)\rangle = b^\top d + \mathrm{tr}(A^\top C)$.
    Let $g$ and $r$ be as in Assumption~\ref{assm:lift} and define the functions
    \begin{align*}
        &G:\RR^d \times \HH \to \RR, \ G(b,A) \coloneq g(T_{b,A}),\\
        &\bar{R}:\RR^d \times \HH \to \RR\cup\{+\infty\}, \ \bar{R}(b,A) \coloneq r(T_{b,A}).
    \end{align*}
    Since $(b,A)\mapsto T_{b,A}$ is a linear isometric isomorphism, $G$ is differentiable at $(b,A)$ if and only if $g$ is differentiable at $T_{b,A}$.
In this case, the chain rule yields
\[
    \nabla G(b,A)=(p,P)
    \quad\Longleftrightarrow\quad
    \nabla_\cH g(T_{b,A}) = T_{p, P}.
\]
Similarly, for the stochastic gradient setting,
define $\cO:\RR^d\times\HH\times\Xi\to\RR^d\times\HH$
through the correspondence
\[
    \cO(b,A,\xi)=(p,P)
    \quad\Longleftrightarrow\quad
    \cO_\cH(T_{b,A},\xi)=T_{p,P}.
\]
We further have by the isometric isomorphism that $\bar{R}$ is proper, convex,
and lower semicontinuous if and only if $r$ is. Additionally, if $T_{b,A} \in \cH$, then the proximal problem \begin{equation*}
    \label{eq:affine_prox}
    \min_{T \in \cT} r(T) + \frac{1}{2h}\|T-T_{b,A}\|_{L^2(\rho)}^2,\end{equation*}
    has the unique minimizer $T_{b^\star,A^\star} \in \cT$, where $(b^\star, A^\star)$ is the unique solution to the convex problem
    \begin{equation}
    \label{eq:affine_coeff_prox}
    \min_{(b',A') \in \RR^{d} \times \TT}
\bar{R}(b',A')
+
\frac{1}{2h}
(
\|b'-b\|_2^2+\|A'-A\|_{\mathrm{F}}^2
).\end{equation}

    We are now ready to state \texttt{BW-FISTA} in terms of the affine coefficients of the iterates in Algorithm~\ref{alg:W-FISTA}. Let $G$ and $\bar{R}$ be as above and let $(b_k, A_k), (b_k^Y, A_k^Y)$ and $(b_k^Z, A_k^Z)$ denote the affine coefficients of $T_k$, $Y_k$ and $Z_k$, respectively.
    The affine coefficients of $\nabla_{\cH} g(Y_k)$ and $\cO_\cH(Y_k,\xi_k)$ in lines~\ref{line:det}--\ref{line:stochastic} in Algorithm \ref{alg:W-FISTA} are given by $\nabla G(b_k^Y, A_k^Y)$ and $\cO(b_k^Y, A_k^Y,\xi_k)$, respectively.
Finally, the proximal problem in line~\ref{line:proximal_step} can be replaced by Problem \eqref{eq:affine_coeff_prox}.
Altogether, \texttt{BW-FISTA} amounts to setting $T_k = T_{b_k, A_k}$ in Algorithm~\ref{alg:W-FISTA}, where the affine coefficients $(b_k,A_k)$ are defined as follows: initialize $(b_{0}, A_{0}) = (b_{-1}, A_{-1}) = (0, I_d)$ and, for $k \geq 0$, set
\begin{equation}
\label{eq:BW_FISTA}
\begin{aligned}
    &(b^Y_k, A^Y_k)=(b_k, A_k)+\beta_k((b_k, A_k)-(b_{k-1}, A_{k-1})),\\
    &(b^Z_k, A^Z_k)
= \begin{cases}
(b^Y_k, A^Y_k)
-
h\nabla G(b_k^Y, A_k^Y) & \text{if deterministic gradients},\\
(b_k^Y, A_k^Y) - h \cO(b_k^Y, A_k^Y, \xi_k) & \text{if stochastic gradients, where $\xi_k \sim \mathbb{P}$,} 
\end{cases}\\
&(b_{k+1}, A_{k+1})
\in 
\argmin_{(b,A) \in \RR^d \times\TT}
\left\{
\bar{R}(b,A) +
\frac{1}{2h}
\left(\|b - b^Z_k\|_2^2 + \|A - A_k^Z\|_\mathrm{F}^2\right)\right\}.
\end{aligned}
\end{equation}

\textbf{Examples of $\cG$.} 
We present a selection of spaces $\TT$ and $\HH$ and smooth components $\cG$ in Table~\ref{tab:G_examples} that satisfy Assumption~\ref{assm:lift}\emph{(i-ii)}. The table lists smooth lifts $g$ satisfying Assumption~\ref{assm:lift}\emph{(ii)} that are both canonical (Assumption~\ref{assm:extras}\emph{(iii)}) and non-canonical. The next proposition justifies the table entries.%

\begin{table}[t]
\centering
\caption{Examples of functions $\cG$ and spaces $\TT,\HH$
satisfying Assumption~\ref{assm:lift}\emph{(i--ii)}
with appropriate lifts $g$; see
Proposition~\ref{prop:sufficient_g} for details. In the first two rows, $g$ is the canonical lift in
Assumption~\ref{assm:extras}\emph{(iii)}, while the
remaining rows list non-canonical smooth lifts $g$.}
\label{tab:G_examples}
\resizebox{\textwidth}{!}{
\begin{tabular}{ll|lll}
\toprule

$\cG(\mu)$ & $\TT,\HH$ & $g(T_{b,A})$
& $\nabla G(b,A)$ & $\ell,L$\\
\midrule

$\int V\dd\mu^{\otimes n}$
& $\TT,\HH$ (general)
& $\int V\dd((T_{b,A})_\sharp\rho)^{\otimes n}$ (canonical)
& $\left(\EE_\mu[w_\mu(X)],
    \mathrm{proj}_\HH[\EE_\mu[Dw_\mu(X)]A]\right)$
& $n\ell_V,nL_V$ \\[2mm]

$G_1(m_\mu,\Sigma_\mu)$
& $\TT,\HH$ (general)
& $G_1(b,AA^\top)$ (canonical)
& $\left(\nabla_1 G_1(b,AA^\top),
    \mathrm{proj}_\HH[2\nabla_2 G_1(b,AA^\top)A]\right)$
& $\ell_1,L_1$ \\[2mm]

$G_2(m_\mu,\Sigma_\mu^{1/2})$
& $\SS_+^d,\SS^d$
& $G_2(b,A)$
& $\left(\nabla_1 G_2(b,A),\nabla_2 G_2(b,A)\right)$
& $\ell_2,L_2$ \\[2mm]

$\frac12 W_2^2(\mu,\cN(m_\nu,\Sigma_\nu))$
& $\SS_+^d\Sigma_\nu^{1/2},\SS^d\Sigma_\nu^{1/2}$
& $\frac12\|b-m_\nu\|_2^2
    +\frac12\|A-\Sigma_\nu^{1/2}\|_\mathrm{F}^2$
& $(b-m_\nu,A-\Sigma_\nu^{1/2})$
& $1,1$ \\[1mm]

\bottomrule
\end{tabular}}
\end{table}

\begin{proposition}
\label{prop:sufficient_g}
Consider the functions $\cG$, spaces
$\TT,\HH$ such that $\TT \subseteq \HH \subseteq \RR^{d\times d}$, $\TT$ is closed and convex and $\HH$ is a linear subspace, and lifts $g$ in
Table~\ref{tab:G_examples}, under the respective assumptions: (i) $V\in C^2((\RR^d)^n)$ is $\ell_V$-strongly
convex and $L_V$-smooth for $0 \leq \ell_V  \leq L_V < \infty$, $L_V > 0$;
(ii) $G_1:\RR^d\times\SS_+^d\to\RR$ is Fréchet differentiable relative to its domain, and $(b,A)\mapsto G_1(b,AA^\top)$ is $\ell_1$-strongly convex and $L_1$-smooth on $\RR^d\times\HH$; (iii) $G_2:\RR^d\times\SS^d\to\RR$ is differentiable, $\ell_2$-strongly convex, and $L_2$-smooth; (iv) $\Sigma_\nu\succ0$. Here, $0\leq\ell_i\leq L_i<\infty$, $L_i>0$, for $i \in \{1,2\}$. Whenever $\TT$ and $\HH$ are not specified, they are chosen such that Assumption~\ref{assm:lift}(i) holds. Then each choice in Table~\ref{tab:G_examples} satisfies
Assumption~\ref{assm:lift}(i-ii), with $\nabla_\cH g(T_{b,A})=T_{\nabla G(b,A)}$, where $\nabla G(b,A)$ and the constants $\ell,L$ are given in Table~\ref{tab:G_examples}.
\end{proposition}

For interaction energies, a stochastic $\cH$-gradient oracle can be constructed by projecting the $\mathsf{BW}$-gradient oracle in Proposition~\ref{prop:stochastic-Q-examples} onto $\cH$, which the following corollary makes precise.

\begin{corollary}
\label{cor:BW_variance}
Consider the setting of
Proposition~\ref{prop:stochastic-Q-examples}
with $\sQ=\BW$, and let $\rho$, $\cT$, $\cH$,
$\TT$, and $\HH$ be as in \eqref{eq:BW-choices},
where $\TT$ is chosen such that
Assumption~\ref{assm:lift}(i) holds.
Let $g(T)=\cG(T_\sharp\rho)$ be the canonical lift
on $\cH$.
Taking $(\Xi,\PP)=(\RR^d,\rho)$, the map
$\cO_\cH:\cH\times\Xi\to\cH$, defined by
\[
\cO_\cH(T_{b,A},\xi)(x)
\coloneq
w_\mu(T_{b,A}(\xi))
+
\mathrm{proj}_\HH[Dw_\mu(T_{b,A}(\xi))A]x,
\quad \mu=(T_{b,A})_\sharp\rho,
\]
is a stochastic $\cH$-gradient oracle for $g$
in the sense of Definition~\ref{def:oracle-Hilbert},
with $\sigma^2=2n^2L_V^2$ and
$\zeta_\cH(T_{b,A})=\|A\|_{\mathrm F}$.
\end{corollary}

Note that the results above extend those for
potential energies in the BW setting of Section~\ref{sec:vi} to higher-order interaction energies and more general affine parametrizations.

\textbf{Examples of $\cR$.} We restrict attention to a special case when $\cR$ decomposes into a function of the mean of $\mu$ and a spectral function of the covariance of $\mu$. More precisely, \begin{equation}\label{eq:spectral_R}\cR(\mu) = R_1(m_\mu)+ R_2(\lambda(\Sigma_\mu^{1/2})), \quad \mu \in \mathsf{BW},\end{equation}
where $R_1:\RR^d \to \RR \cup \{+\infty\}$ and $R_2:\mathbb{R}^d_+\to \mathbb{R}\cup\{+\infty\}$ are proper, lower semicontinuous and convex, and $R_2$ is permutation invariant, that is, $R_2 \circ \sigma = R_2$ for all permutations $\sigma:\mathbb{R}^d \to \mathbb{R}^d$. The next proposition shows that, if $\TT = \SS_+^d$ and $\HH = \SS^d$, then this decomposition leads to a simple form of the proximal step in~\eqref{eq:BW_FISTA}. The result is a direct consequence of convex analysis for spectral functions; see \cite[Corollary~2.4]{lewis1996convex}.
    \begin{proposition}
    \label{prop:sufficient_H}
    Let $\TT = \SS_+^d$ and $\HH = \SS^d$. Assume that $\cR$ is of the form \eqref{eq:spectral_R}. Then the corresponding lift $r$ satisfies
Assumption~\ref{assm:lift}(iii). Moreover, if $(b,A) \in \mathbb{R}^d\times \HH$ with orthogonal eigendecomposition $A = U \mathrm{diag}(\lambda(A))U^\top$, then the unique minimizer of \eqref{eq:affine_coeff_prox} is the pair $(b^\star, A^\star) \in \RR^d \times \TT$ given by
\begin{equation}
\label{eq:easy_proximal}
   \begin{aligned}
        &b^\star
= 
\argmin_{b' \in \mathbb{R}^d}
\left\{
R_1(b') + 
\frac{1}{2h}
\|b' - b\|_2^2\right\},\\
&A^\star = U\mathrm{diag}(a^\star)U^\top, \ a^\star
= 
\argmin_{a \in  \mathbb{R}^{d}_{+}}
\left\{
R_2(a) + 
\frac{1}{2h}\|a - \lambda(A)\|_2^2\right\}.
\end{aligned}
\end{equation}
\end{proposition}
Table~\ref{tab:R_prox} lists examples of regularizers
$\cR$ satisfying the assumptions of
Proposition~\ref{prop:sufficient_H}.
For simplicity, we restrict attention to cases where
$R_1$ is constant, in which case $b^\star=b$.
For each example, the proximal problem associated
with $R_2$ is numerically tractable; see, e.g.,
\citep[Section~6]{parikh2014proximal}.

\begin{table}[t]
\centering
\caption{Examples of functions $\cR$ satisfying Assumption~\ref{assm:lift}\emph{(iii)} when $\TT = \SS_+^d$ and $\HH = \SS^d$, with the corresponding functions $R_2$ in Proposition \ref{prop:sufficient_H} ($R_1$ is constant). The Schatten $p$-norm and the operator norm are respectively given by $\|S\|_{\mathrm{S}_p} \coloneq \|\lambda(S)\|_p$ and $\|S\|_\mathrm{op} = \lambda_{\max}(S)$, for all $S \in \mathbb{S}_+^d$, and $C_{\underline{\sigma}, \overline{\sigma}}\coloneq \{\mu \in \mathsf{BW}: \underline{\sigma}^2I_d \preceq \Sigma_\mu \preceq \overline{\sigma}^2I_d\}$, for $0 \leq \underline{\sigma} \leq \overline{\sigma}$.}
\label{tab:R_prox}
\vspace{1em}
\resizebox{0.35\textwidth}{!}{
\begin{tabular}{l|l}
\toprule
$\mathcal R(\mu)$
& $R_2(\lambda),\ \lambda\in \mathbb{R}_+^d$ \\
\midrule

$-\operatorname{Ent}(\mu)$
&
$-\sum_{i=1}^d\log \lambda_i$ %

\\[2mm]

$\operatorname{tr}[(\Sigma_\mu^{1/2})^{-q}],\ q>0$
&
$\displaystyle\sum_i \lambda_i^{-q}$

\\[2mm]

$\operatorname{tr}[(\Sigma_\mu^{1/2})^{p}], \ 1 \leq p\neq 2$
&
$\|\lambda\|_p^p$

\\[2mm]

$\|\Sigma_\mu^{1/2}\|_{\mathrm{S}_p}, \ p \geq 1$
&
$\|\lambda\|_p$

\\[2mm]

$\|\Sigma_\mu^{1/2}\|_{\mathrm{op}}$
&
$\max_{i=1,\dots,d}\lambda_i$

\\[2mm]

$\iota_{ C_{\underline\sigma,\overline\sigma}}(\mu)$
&
$\iota_{[\underline\sigma,\overline\sigma]^d}
(\lambda)$

\\

\bottomrule
\end{tabular}
}
\end{table}

\begin{remark}
Gaussian VI was previously considered through the
forward-backward Gaussian VI (FB-GVI) method of
\cite{diao2023forward}, which achieves an
$O(k^{-1})$ objective-gap rate for deterministic
gradients and convex potentials, in contrast to
our $O(k^{-2})$ rate. Even when setting the momentum weights $\beta_k = 0$,
the methods differ in that FB-GVI performs an
intrinsic $\BW$-gradient step followed by a JKO step,
whereas \texttt{BW-FISTA} performs both steps in the
lifted Hilbert space. Consequently, 
the two methods generally produce different
iterates.
\end{remark}

\begin{remark}[Lower triangular parameterization]
\label{rem:lower-triangular}
The proximal problem in \eqref{eq:affine_coeff_prox} for Gaussian VI also admits an efficient implementation with $\TT=\LL_+^d$ and $\HH=\LL^d$. Indeed, if $A \in \LL^d$, then the proximal problem for the linear coefficient is
$$
\min_{L\in\mathbb L_+^d}
\left\{
-\log\det L+\frac1{2h}\|L-A\|_\mathrm{F}^2
\right\},
$$
with the unique solution is $L^\star = A + \Delta A$, where $(\Delta A)_{ij} = \frac{1}{2}((A_{ii}^2 + 4h)^{1/2} - A_{ii})$ for $i = j$ and $(\Delta A)_{ij} = 0$ otherwise \citep[Theorem~13]{domke2020provable}. The proximal step can thus be implemented more efficiently than the eigendecomposition approach of Proposition~\ref{prop:sufficient_H}, where $\TT = \SS_+^d$ and $\HH = \SS^d$. However, the symmetric parametrization $\TT=\SS_+^d$ has the
advantage that the maps in $\cT$ are cyclically monotone.
Hence Assumption~\ref{assm:extras}(ii) holds, yielding
$R=W_2(\mu_0,\mu^\star)$ in the convergence rates of Theorems~\ref{thm:w-fista_rate} and \ref{thm:w-fista_rate_stochastic}.
Gaussian VI with proximal steps over $\LL_+^d$ was
previously studied by \cite{Kim_2023_VI} in a
nonaccelerated stochastic gradient setting.
\end{remark}

    \subsection{Finite-Dimensional Wasserstein FISTA}
    \label{sec:finite-W-FISTA}

    The affine coefficient representation of \texttt{BW-FISTA} 
in Section~\ref{sec:BW-FISTA} extends naturally to
general finite-dimensional Hilbert spaces of transport
maps. Moreover, finite-dimensional parametrizations
can be used to approximate infinite-dimensional
spaces of transport maps, as we
demonstrate in the subsequent subsections.

\textbf{Finite-dimensional $\cH$.} Suppose Assumption~\ref{assm:lift} holds with
$N \defeq \dim(\cH)<\infty$, and let
$\Phi=\{\phi_j\}_{j=1}^N$ be a basis of $\cH$.
For $\theta\in\RR^N$, define
\[
    T_\theta\defeq\sum_{j=1}^N\theta_j\phi_j,
    \quad
    \Theta\defeq\{\theta\in\RR^N:T_\theta\in\cT\}.
\]
Here, $\Theta$ is closed and convex, since $\cT$ is closed and convex and $\theta \mapsto T_\theta$ is linear and continuous. Let $M\in\RR^{N\times N}$ denote the Gram matrix of
$\Phi$, with entries
$M_{ij}=\langle\phi_i,\phi_j\rangle_{L^2(\rho)}$.
Since $\Phi$ is linearly independent, $M$ is positive
definite. We have the linear isometric isomorphism
\[
    \RR^N\to\cH,\ \theta\mapsto T_\theta,
\]
where $\RR^N$ is equipped with the inner product
$\langle\theta,\eta\rangle_M=\theta^\top M\eta$.

Under the above assumptions, Algorithm~\ref{alg:W-FISTA} reduces to a version of the Euclidean FISTA over the parameter $\theta$. To see this, let $g$ and $r$ be as in Assumption~\ref{assm:lift}
and define the functions
\begin{align*}
    &G:\RR^N\to\RR,\quad
    G(\theta)\coloneq g(T_\theta),\\
    &\bar{R}:\RR^N\to\RR\cup\{+\infty\},\quad
    \bar{R}(\theta)\coloneq r(T_\theta).
\end{align*}
Since $\theta\mapsto T_\theta$ is a linear isometric
isomorphism, $G$ is differentiable at $\theta$ if and
only if $g$ is differentiable at $T_\theta$. In this case, writing
$\nabla_\cH g(T_\theta)=\sum_{j=1}^N c_j(\theta)\phi_j$,
the chain rule yields
\[
\begin{aligned}
    \frac{\partial G(\theta)}{\partial\theta_i} =
    \left\langle
        \nabla_\cH g(T_\theta),
        \frac{\partial T_\theta}{\partial\theta_i}
    \right\rangle_{L^2(\rho)} =
    \langle\nabla_\cH g(T_\theta),\phi_i\rangle_{L^2(\rho)}
    =
    \sum_{j=1}^N M_{ij}c_j(\theta), \quad \forall i \in [N].
\end{aligned}
\]
Thus, $$\nabla G(\theta) = \left(\langle \nabla_\cH g(T_\theta), \phi_i \rangle_{L^2(\rho)}\right)_{i=1}^N, \quad c(\theta) = M^{-1}\nabla G(\theta),$$ for all $\theta \in \RR^N$, where $\nabla G$ denotes the gradient with respect
to the standard Euclidean inner product. Similarly, for the stochastic gradient setting, let $\widehat{c}: \RR^N \times \Xi \to \RR^N$ be such that $\cO_\cH(T_\theta,\xi) = \sum_{j=1}^N \widehat{c}_j(\theta, \xi) \phi_j$ for all $(\theta, \xi) \in \RR^N \times \Xi$ and define $$\cO(\theta, \xi): \RR^N \times \Xi \to \RR^N, \ \cO_i(\theta, \xi) \coloneq \langle \cO_\cH(T_\theta, \xi), \phi_i \rangle_{L^2(\rho)}.$$ It follows that $\cO(\theta, \xi) = M\widehat{c}(\theta, \xi)$, or equivalently, $$\widehat{c}(\theta, \xi) = M^{-1}\cO(\theta, \xi).$$

We further have by the isometric isomorphism that
$\bar{R}$ is proper, convex, and lower semicontinuous if
and only if $r$ is. Additionally, for $T_z\in\cH$,
the proximal problem
\[
    \min_{T\in\cT}
    r(T)+\frac{1}{2h}\|T-T_z\|_{L^2(\rho)}^2
\]
has unique minimizer $T_{\theta^\star}\in\cT$,
where $\theta^\star$ is the unique solution to
the strongly convex problem
\begin{equation}
    \label{eq:finite_coeff_prox}
    \min_{\theta\in\Theta}
    \bar{R}(\theta)+\frac{1}{2h}\|\theta-z\|_M^2.
\end{equation}
    
By the above observations, \texttt{W-FISTA} in
Algorithm~\ref{alg:W-FISTA} reduces to its finite-dimensional
counterpart in Algorithm~\ref{alg:finite-W-FISTA}, where $\theta_k$ denotes the coefficient vector of $T_k$. The convergence guarantees of
Theorems~\ref{thm:w-fista_rate} and
\ref{thm:w-fista_rate_stochastic} remain under their
respective assumptions.

\textbf{Finite-dimensional restrictions.}
When $\cH$ is infinite-dimensional, we may approximate
by the finite-dimensional methods above by selecting a basis $\Phi$ for a finite-dimensional
linear subspace $\cH_\Phi\subseteq\cH$. Define the restrictions of $\cT$ and $\mathsf{Q}$, respectively, by
\[
    \cT_\Phi\coloneq\cT\cap\cH_\Phi,
    \quad
    \sQ_\Phi\coloneq(\cT_\Phi)_\sharp\rho,
\]
and consider the restrictions $g|_{\cH_\Phi}$ and $r|_{\cH_\Phi}$. Provided that $\mathrm{Id} \in \cH_\Phi$ and $r|_{\cH_\Phi}$ is proper, the restricted
setting satisfies Assumption~\ref{assm:lift} with the
same strong convexity and smoothness constants $\ell$ and $L$, respectively.
Moreover, the gradient and stochastic gradient oracle
for the restricted smooth lift are given by
\[
\begin{aligned}
    &\nabla_{\cH_\Phi}(g|_{\cH_\Phi})(T)
    =
    \proj_{\cH_\Phi}[\nabla_\cH g(T)],\\
    &\cO_{\cH_\Phi}(T,\xi)
    \coloneq
    \proj_{\cH_\Phi}[\cO_\cH(T,\xi)],
\end{aligned}
\]
respectively. The latter retains the noise parameters
$\sigma$ and $\zeta_\cH|_{\cH_\Phi}$ by linearity
and nonexpansiveness of the orthogonal projection.

    \begin{algorithm}[t]
    \caption{Finite-dimensional \texttt{W-FISTA}}
    \label{alg:finite-W-FISTA}
    \begin{algorithmic}[1]
    \Statex \textbf{Parameters:} Step size $h>0$, momentum tuning constant $q \in [0,1]$
    \State Initialize $\theta_{-1}=\theta_0 \in \RR^N$ such that $T_{\theta_{-1}} = T_{\theta_{0}} = \mathrm{Id}$, $\, t_0 = 0$
    \For{$k=0,1,2,\ldots$}
        \State $T_{\theta_k} \gets \sum_{j=1}^N(\theta_k)_j\phi_j$
        \State $t_{k+1}
        =
        \frac{1}{2}(1-qt_k^2+\sqrt{(1-qt_k^2)^2+4t_k^2}), \ \beta_k=\frac{t_k-1}{t_{k+1}+qt_k^2} $
        \State $y_k \gets \theta_k+\beta_k(\theta_k-\theta_{k-1})$
        \If{deterministic gradients}
        \State $z_k \gets y_k - h M^{-1}\nabla G(y_k)$
        \ElsIf{stochastic gradients}
        \State $z_k \gets y_k-hM^{-1}\cO(y_k,\xi_k)$, where $\xi_k \sim \PP$ is drawn independently
        \EndIf
        \State
        $\theta_{k+1}
        =
        \argmin_{\theta\in \Theta}
        \left\{
        \bar{R}(\theta)
        +
        \frac{1}{2h}
        \|\theta-z_k\|_{M}^2
        \right\}
        $ \label{line:finite-dim-prox}
    \EndFor
    \end{algorithmic}
    \end{algorithm}

    \subsection{Mean-Field Wasserstein FISTA}
    \label{sec:MF-FISTA}

    Recall the special case of Problem~\eqref{eq:composite_problem}
    obtained by restricting to mean-field distributions and taking 
    $\sQ=\MF$. We refer to \texttt{MF-FISTA} as
    Algorithm~\ref{alg:W-FISTA} with the standard choices
    of transport spaces in Proposition~\ref{prop:standard_choices},
    namely,
    \[
        \rho=\gamma^{\otimes d}=\cN(0,I_d),\qquad
        \cH=L^2(\rho;\RR^d)\cap\cM_{\mathrm{sep}},
\qquad
\cT=\cH\cap\cM_{\mathrm{cm}}(\rho),
    \]
    where $\gamma\coloneq\cN(0,1)$.

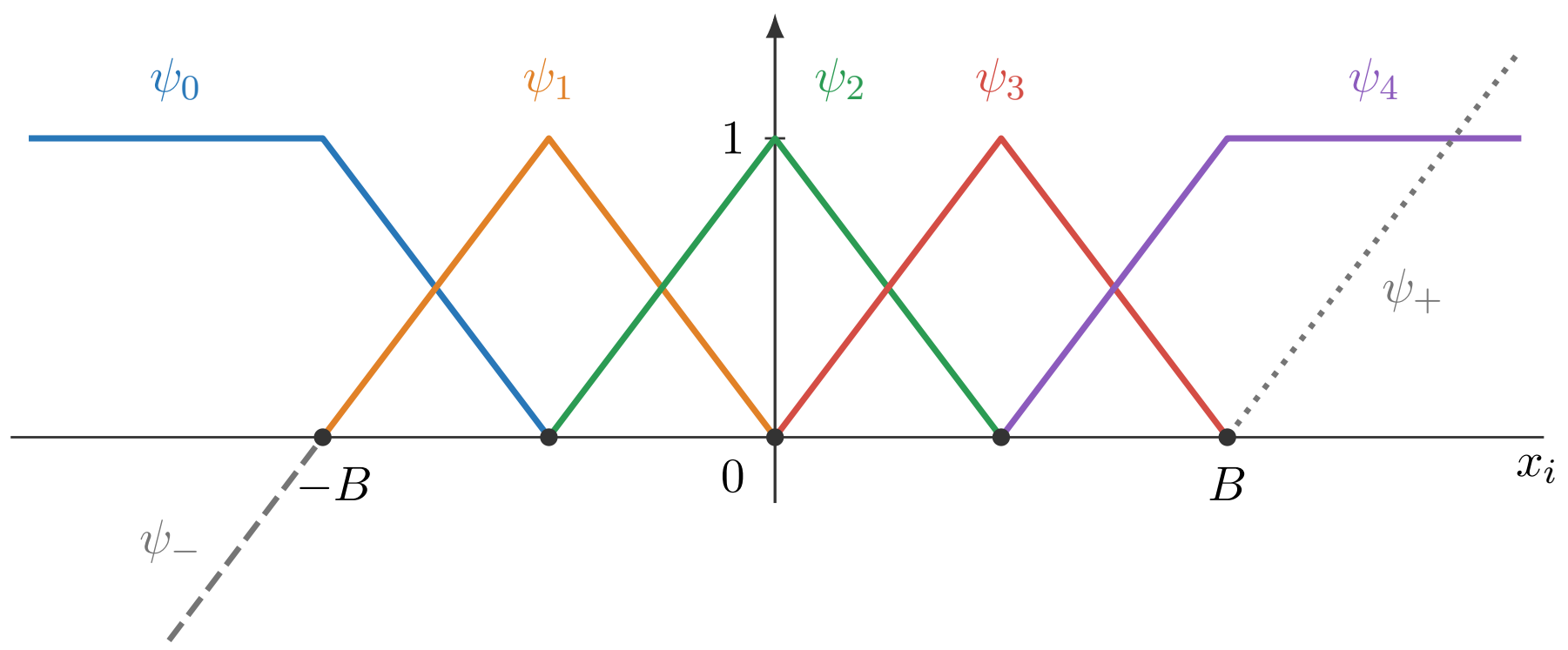
\begin{figure}[t]
    \centering
    \resizebox{0.75\linewidth}{!}{%
        \begin{tikzpicture}[
  x=1.8cm, y=2.38cm,
  basis/.style={line width=1.4pt, line cap=butt, line join=round},
  tail/.style={line width=1.3pt, dash pattern=on 1.3pt off 2.15pt,
               line cap=butt},
  tail_1/.style={line width=1.3pt, dash pattern=on 5pt off 2.15pt,
               line cap=butt},
  axis/.style={line width=0.65pt, -{Latex[length=2.0mm,width=1.5mm]}},
  every node/.style={inner sep=1.5pt}
]
  \definecolor{basisBlue}{RGB}{40,119,184}
  \definecolor{basisOrange}{RGB}{225,129,39}
  \definecolor{basisGreen}{RGB}{43,155,83}
  \definecolor{basisRed}{RGB}{212,77,69}
  \definecolor{basisPurple}{RGB}{140,91,189}
  \definecolor{basisGray}{gray}{0.46}
  \definecolor{basisAxis}{gray}{0.20}

  \draw[basisAxis, line width=0.55pt] (-3.38,0) -- (3.40,0);
  \draw[axis, basisAxis] (0,-0.22) -- (0,1.42);
  \draw[basisAxis, line width=0.55pt] (-0.045,1) -- (0.045,1);

\draw[basis, basisBlue]
  (-3.30,1) -- (-2,1) -- (-1,0);

\draw[basis, basisOrange]
  (-2,0) -- (-1,1) -- (0,0);

\draw[basis, basisGreen]
  (-1,0) -- (0,1) -- (1,0);

\draw[basis, basisRed]
  (0,0) -- (1,1) -- (2,0);

\draw[basis, basisPurple]
  (1,0) -- (2,1) -- (3.30,1);

  \draw[tail_1, basisGray] (-2.68,-0.68) -- (-2,0);
  \draw[tail, basisGray] (2,0) -- (3.30,1.30);

  \foreach \s in {-2,-1,0,1,2}
    \fill[basisAxis] (\s,0) circle[radius=2pt];

  \node[basisBlue,   anchor=south] at (-2.65,1.10) {$\psi_0$};
  \node[basisOrange, anchor=south] at (-1,1.10)    {$\psi_1$};
  \node[basisGreen,  anchor=south west] at (0.14,1.10) {$\psi_2$};
  \node[basisRed,    anchor=south] at (1,1.10)     {$\psi_3$};
  \node[basisPurple, anchor=south] at (2.65,1.10)  {$\psi_4$};
  \node[basisGray, anchor=east] at (-2.5,-0.35) {$\psi_-$};
  \node[basisGray, anchor=west] at (2.65,0.49)   {$\psi_+$};

  \node[anchor=east] at (-0.10,1) {$1$};
  \node[anchor=north] at (-1.95,-0.075) {$-B$};
  \node[anchor=north east] at (-0.10,-0.055) {$0$};
  \node[anchor=north] at (2,-0.075)  {$B$};
  \node[anchor=north] at (3.37,-0.035) {$x_i$};
\end{tikzpicture}
    }
    \caption{The fixed-coordinate nodal basis $\Psi$ for $J = 4$. Only the nonzero values of each basis function are illustrated. The black dots on the horizontal axis represent the discretization points.}
    \label{fig:mf_basis}
\end{figure}

     The coordinate-wise separability of maps in $\cT$ and $\cH$ lends itself to finite-dimensional approximations using piecewise affine functions in each coordinate; see \citep{jiang2025algorithms} for a similar treatment. We consider a nodal basis representation that, in the VI setting, leads to a numerically tractable proximal problem. More precisely, for $B > 0$ and $J \in \NN_+$, we let $a_0 < \dots < a_J$ be the discretization of the interval $[-B,B]$ with discretization width $\delta \coloneq 2B/J$, that is, $a_j \coloneq -B + j\delta$ for $j \in \{0,\dots, J\}$. Define the set of basis functions by $$\Phi \coloneq \{x \mapsto \psi_\ell(x_i) e_i:i\in [d], \, \ell \in\{-,0,\dots,J,+\}\} \subseteq \cH,$$
     where $e_i$ is the $i$-th standard basis vector in $\RR^d$, and $\psi_\ell$ for $\ell \in \{ -,0,\dots,J,+\}$ are defined, respectively, by $\psi_{-}(x) \coloneq (x + B)\mathbf{1}_{\{x\leq -B\}}$, $\psi_{+}(x) \coloneq (x - B)\mathbf{1}_{\{x\geq B\}}$,
\begin{align*}
&\psi_0(x) \coloneq
\begin{cases}
1,&x<a_0,\\
\dfrac{a_1-x}{\delta},&x \in [a_0, a_1],\\
0,&x>a_1.
\end{cases} \quad \psi_J(x) \coloneq
\begin{cases}
0,&x<a_{J-1},\\
\dfrac{x-a_{J-1}}{\delta},&x \in [a_{J-1}, a_J],\\
1,&x>a_J.
\end{cases} \\
&\qquad \qquad \psi_j(x) \coloneq \begin{cases}
        \dfrac{x - a_{j-1}}{\delta}, &x\in[a_{j-1}, a_j],\\
        \dfrac{a_{j+1} - x}{\delta}, &x\in(a_{j}, a_{j+1}],\\
        0, & \text{otherwise},
    \end{cases} \quad \forall j \in \{1,\dots,J-1\}.
\end{align*}
    Here, we suppress the the dependence of $\Phi$ on $B$ and $J$ to avoid clutter. Denote by $\Psi \coloneq \{\psi_\ell:\ell \in\{-,0,\dots,J,+\}\}$ the fixed-coordinate nodal basis, which is illustrated in Figure~\ref{fig:mf_basis}. The Gram matrix of $\Phi$ is $I_d \otimes M_\Psi$, where $M_\Psi$ is the Gram matrix of $\Psi$. Due to the local support of the basis functions $\Psi$, $M_\Psi$ is tridiagonal.

    Define the $d(J+3)$-dimensional subspace $\cH_{\Phi} \coloneq \mathrm{span}(\Phi)$ of $\cH$ and the restrictions of $\cT$ and $\MF$ by $\cT_\Phi = \cT \cap \cH_\Phi$ and $\MF_\Phi = (\cT_\Phi)_\sharp\rho$, respectively. Every $T_\theta \in \cH_\Phi$ has the representation $$T_{\theta}(x) = \sum_{i=1}^d \theta^{i\top} \psi(x_i) e_i,
    $$ 
    where we write $\theta = (\theta^{1\top}, \dots, \theta^{d\top})^\top$ with $\theta^i = (s_-^i,q_0^i,\dots,q_J^i,s_+^i)^\top \in \RR^{J+3}$ and define the vector of basis functions as $\psi(x_i)\coloneq
(\psi_-(x_i),\psi_0(x_i),\ldots,\psi_J(x_i),\psi_+(x_i))^\top$. The admissible transport maps are given by
    $$\cT_{\Phi} = \left\{T_{\theta} \in \cH_{\Phi}: q_0^i \leq \dots \leq q_J^i, \ s_-^i \geq 0, \ s_+^i \geq 0, \ \forall i\in [d]\right\}.$$ In other words, every $T_{\theta} \in \cH_{\Phi}$ is coordinate-wise separable with piecewise linear components. The $i$-th component is parametrized by its values $q_0^i,\ldots,q_J^i$ at $a_0,\ldots,a_J$ and its slopes $s_-^i$ and $s_+^i$ on $(-\infty,-B]$ and $[B,\infty)$, respectively.
    Every $T_{\theta} \in \cT_{\Phi}$ is, in addition, nondecreasing.

    \textbf{Approximation guarantees for mean-field VI.} We refer to Problem \eqref{eq:VI_in_sec} with $\sQ = \mathsf{MF}_\Phi$ as approximate mean-field VI.
    By refining the basis $\Phi$, we obtain increasingly accurate approximations of the spaces $\cT$ and $\cH$ through their restrictions $\cT_\Phi$ and $\cH_\Phi$, respectively. This allows us to approximate the mean-field VI minimizer to arbitrary accuracy, under the following assumption.

    \begin{assumption}
    \label{assm:mean_field_V}
    The potential $V\in C^2(\RR^d)$ is
    $\ell_V$-strongly convex and $L_V$-smooth,
    with $0<\ell_V\leq L_V<\infty$.
    \end{assumption}

Under Assumption~\ref{assm:mean_field_V}, mean-field VI
admits a unique minimizer
\cite[Proposition~5.2]{jiang2025algorithms}.
Approximate mean-field VI likewise admits a unique
minimizer by strong convexity of the lifted objective
\cite[Corollary~11.17]{bauschke2017convex}. To that end, let
$$\pi^\star \coloneq\argmin_{\mu \in \mathsf{MF}} \mathrm{KL}(\mu \, \|\, \pi), \quad \pi_\Phi^\star \coloneq\argmin_{\mu \in \mathsf{MF}_\Phi} \mathrm{KL}(\mu \, \|\, \pi).$$
be the respective minimizers. For fixed $B$, $J$ and provided that the discretization grids coincide, the class $\cT_\Phi$ contains the class of piecewise linear transport maps
constructed in \citet[Appendix~C.3]{jiang2025algorithms} This yields the following extension of their
approximation result in Theorem~5.9. %

\begin{proposition}
\label{prop:MF_approximation}
Let Assumption~\ref{assm:mean_field_V} hold,
and set $\kappa\coloneq L_V/\ell_V$.
For every $\varepsilon \in (0,1]$,
there exist
\[
B=O\left(
\sqrt{\log\frac{\kappa^3 d}{\varepsilon}}
\right),
\quad
J=\widetilde O\left(
\kappa^2\sqrt{\frac{d}{\varepsilon}}
\right),
\]
such that the respective mean-field VI minimizers
$\pi^\star$ and $\pi_\Phi^\star$ satisfy
\[
\frac{\ell_V}{2} W_2(\pi_\Phi^\star,\pi^\star)^2
\leq
\cF(\pi_\Phi^\star)-\cF(\pi^\star)
\lesssim\varepsilon.
\]
\end{proposition}

\textbf{Applicability of finite-dimensional \texttt{W-FISTA}.}
Solving approximate mean-field VI to $\varepsilon$ accuracy can be achieved at accelerated rates due to Theorems~\ref{thm:w-fista_rate} and \ref{thm:w-fista_rate_stochastic} via finite-dimensional \texttt{W-FISTA} in Algorithm~\ref{alg:finite-W-FISTA}. Indeed, mean-field VI satisfies
    Assumption~\ref{assm:lift} by
    Proposition~\ref{prop:vi_summary}. Consequently, since $\Id\in\cH_\Phi$ and
    $r(\Id)=-\Ent(\rho)<\infty$,
   {approximate mean-field VI} also satisfies
    Assumption~\ref{assm:lift} by the restriction argument in Section~\ref{sec:finite-W-FISTA} with spaces $\cT_\Phi$ and $\cH_\Phi$, lifts $g|_{\cH_{\Phi}}$ and $r|_{\cH_{\Phi}}$ and constants
    $\ell=\ell_V$ and $L=L_V$. We refer to Algorithm~\ref{alg:finite-W-FISTA} applied to approximate mean-field VI with spaces $\cT_\Phi$ and $\cH_\Phi$ as finite-dimensional \texttt{MF-FISTA}. We describe below how to compute the deterministic and stochastic gradients as well as the proximal operator associated with finite-dimensional \texttt{MF-FISTA}.

First, we treat the gradients. By Proposition~\ref{prop:vi_summary}, deterministic and stochastic $\cH$-gradient oracles for
mean-field VI are given in Table~\ref{tab:vi_summary}.
The corresponding $\cH_\Phi$-gradient oracles for {approximate mean-field VI} are obtained by orthogonally projecting $\cH$-gradient counterparts onto $\cH_\Phi$. Since projection onto $\cH_\Phi$ amounts to projection onto $\operatorname{span}(\Psi)$ in each coordinate, the Gram-matrix formula for orthogonal projection yields, for all $T\in\cH_\Phi$ and $i\in[d]$,
\begin{equation}
\label{eq:grad_finite_mf_vi}
\begin{aligned}
    [\nabla_{\cH_\Phi}g(T)(x)]_i
    &=
    [\mathrm{proj}_{\cH_\Phi}
    [\nabla_\cH g(T)](x)]_i\\
    &=
    \psi(x_i)^\top M_\Psi^{-1}
    \EE_{X\sim\rho}
    [\psi(X_i)\partial_i V(T(X))].
\end{aligned}
\end{equation}
Similarly, taking $(\Xi,\PP)=(\RR^d,\rho)$,
the stochastic gradient oracle is given by
\begin{equation}
\label{eq:stochastic_grad_finite_mf_vi}
\begin{aligned}
    [\cO_{\cH_\Phi}(T,\xi)(x)]_i
    &\coloneq
    [\mathrm{proj}_{\cH_\Phi}
    [\cO_\cH(T,\xi)](x)]_i\\
    &=
    \psi(x_i)^\top M_\Psi^{-1}
    \EE_{X_i\sim\gamma}
    \left[
        \psi(X_i)
        \partial_i V(T_i(X_i),T_{-i}(\xi_{-i}))
    \right].
\end{aligned}
\end{equation}
By linearity and nonexpansiveness of the orthogonal
projection, $\cO_{\cH_\Phi}$ inherits the noise
parameters of $\cO_\cH$, namely, $\sigma_\Phi^2=L_V^2$ and $\zeta_{\cH_\Phi}(T)
    =
    \|T-\EE_\rho[T(X)]\|_{L^2(\rho)}.$

\textbf{Outer iteration complexity.} Combining Proposition~\ref{prop:MF_approximation} with
the convergence guarantees of finite-dimensional
\texttt{MF-FISTA}, we recover the accuracy dependence
of the oracle complexity bounds in
Theorems~\ref{thm:w-fista_rate}
and~\ref{thm:w-fista_rate_stochastic}
for mean-field VI. Indeed, with $\mu_k$ denoting the measure valued iterates of finite-dimensional
\texttt{MF-FISTA}, the optimality gap decomposes as
\[
    \cF(\mu_k)-\cF(\pi^\star)
    =
    \underbrace{\cF(\mu_k)-\cF(\pi_\Phi^\star)}
    _{\text{optimization error}}
    +
    \underbrace{\cF(\pi_\Phi^\star)-\cF(\pi^\star)}
    _{\text{approximation error}}.
\]
To attain an optimality gap of at most $\varepsilon$ for mean field VI,
we 
\begin{itemize}
    \item select $B$ and $J$ according to
Proposition~\ref{prop:MF_approximation} so that the
approximation error is at most $\varepsilon/2$,
\item solve {approximate mean-field VI} by finite-dimensional \texttt{MF-FISTA} to accuracy $\varepsilon/2$.
\end{itemize}
Importantly, the constants in Theorems~\ref{thm:w-fista_rate} and \ref{thm:w-fista_rate_stochastic} of finite-dimensional \texttt{MF-FISTA} can be made independent of $B$ and $J$. This follows since $\sigma$, $\ell$, and $L$
and the noise scaling function $\zeta$ associated with approximate mean-field VI are inherited
from mean-field VI, while the initial radius $R$ remains uniformly bounded
in $B$ and $J$, since
\[
\begin{aligned}
    R^2
    &\leq \|T_\Phi^\star-T_0\|_{L^2(\rho)}^2 \leq \frac{2}{\ell}
    \bigl(\cF(\rho)-\cF(\pi_\Phi^\star)\bigr) \leq \frac{2}{\ell}
    \bigl(\cF(\rho)-\cF(\pi^\star)\bigr),
\end{aligned}
\]
where $T_\Phi^\star$ is a minimizer of approximate mean-field VI.
The first inequality follows from the definition of $R$,
the second from strong convexity of $f=g+r$ and the fact that $T_0=\Id$,
and the last from optimality of $\pi^\star$. Consequently, the outer iteration complexity of
finite-dimensional \texttt{MF-FISTA} for solving mean-field VI
to accuracy $\varepsilon$ is
$O(\log(1/\varepsilon))$ in the deterministic setting and
$O(\varepsilon^{-1}\log(1/\varepsilon))$ in the stochastic
setting by
Theorems~\ref{thm:w-fista_rate}
and~\ref{thm:w-fista_rate_stochastic}, respectively.

\begin{remark}[Comparison to \citealp{jiang2025algorithms}]
\label{rem:MF_deterministic_comparison}
\citet{jiang2025algorithms} apply accelerated projected
gradient methods to the lifted KL objective by
restricting the transport maps to a domain on which
the negative entropy is smooth.
The resulting smoothness bound deteriorates as the
discretization is refined, leading to a polynomial
dependence on $1/\varepsilon$ in their
deterministic iteration complexity \citep[Corollary~5.14]{jiang2025algorithms}.
In contrast,
finite-dimensional \texttt{MF-FISTA} retains the logarithmic dependence on $1/\varepsilon$
in its deterministic iteration complexity
at the expense of solving a numerically tractable proximal
problem at each iteration.
\end{remark}

\begin{remark}[Computation of the gradients]
\label{rem:MF-VI-approx}
A computational challenge in finite-dimensional \texttt{MF-FISTA} is the evaluation of the
expectations appearing in the gradient formulas \eqref{eq:grad_finite_mf_vi} and \eqref{eq:stochastic_grad_finite_mf_vi}. While the $\cH_\Phi$-gradient involves the $d$-dimensional
expectation,
the stochastic $\cH_\Phi$-gradient oracle requires only the
one-dimensional expectation arising from the projection.
Replacing the latter by a sampling-based
stochastic $\cH_\Phi$-gradient oracle introduces additional variance, potentially making the
noise parameter $\sigma$ dependent on the basis
parameters $B$ and $J$. Since these parameters are chosen
as a function of the prescribed accuracy $\varepsilon$,
such a sampling scheme would lead to different end-to-end
complexity bounds than those described above. Nevertheless, the existing stochastic treatment in \cite{jiang2025algorithms}
suggest that the dependence on $B$ and $J$ can be controlled.
We leave the analysis of such sampling schemes
for future work.
\end{remark}

We next derive the proximal problem in \eqref{eq:finite_coeff_prox} for finite-dimensional \texttt{MF-FISTA}. By coordinate-wise separability of the transport maps, the proximal problem
decomposes into $d$ independent subproblems. Moreover, the choice of basis $\Phi$ is particularly
attractive computationally, as each subproblem admits a tridiagonal Hessian and can thus be solved efficiently via Newton's method.

For $T_\theta\in\cT_\Phi$, we have by \citep[Section~5.5.1]{jiang2025algorithms}) that
\begin{align*}
    \bar{R}(\theta)
    &=-\Ent(\rho)
    -\sum_{i=1}^d\EE_{X_i\sim\gamma}
    [\log (T_\theta)_i'(X_i)]\\
    &=-\Ent(\rho)
    -\sum_{i=1}^d p^\top\log(D\theta^i),
\end{align*}
where $p\coloneq
    \bigl(
        \gamma((-\infty,a_0)),
        \gamma((a_0,a_1)),\ldots,
        \gamma((a_{J-1},a_J)),
        \gamma((a_J,\infty))
    \bigr)^\top,$ and $D\in\RR^{(J+2)\times(J+3)}$ is the matrix
satisfying
\[
    D\theta^i=
    \left(
        s_-^i,
        \frac{q_1^i-q_0^i}{\delta},
        \ldots,
        \frac{q_J^i-q_{J-1}^i}{\delta},
        s_+^i
    \right)^\top.
\]
Here, the logarithm is applied componentwise,
with the convention $\log(0)=-\infty$. The proximal
problem in \eqref{eq:finite_coeff_prox} separates into $d$ independent subproblems, 
\begin{equation}
    \label{eq:mf_coeff_prox}
    \min_{\theta:\,D\theta>0}
    \left\{
        -p^\top\log(D\theta)
        +\frac{1}{2h}\|\theta-z^i\|_{M_\Psi}^2
    \right\}, \quad i \in [d],
\end{equation}
where we used that $M_\Phi=I_d\otimes M_\Psi$ to separate the quadratic term $\|\theta - z\|_{M_\Phi}^2$, where $z = (z^{1\top},\dots, z^{d\top})^\top \in \RR^{d(J+3)}$.
Let $F_i$, for $i \in [d]$, denote the objectives in
\eqref{eq:mf_coeff_prox}. It is readily verified that the gradients and Hessians
are
\begin{align*}
    \nabla F_i(\theta)
    &=
    -D^\top
    \left(\frac{p_\ell}{(D\theta)_\ell}\right)_\ell
    +\frac{1}{h}M_\Psi(\theta-z^i),\\
    \nabla^2 F_i(\theta)
    &=
    D^\top\mathrm{diag}
    \left(\frac{p_\ell}{(D\theta)_\ell^2}\right)_\ell D
    +\frac{1}{h}M_\Psi, \quad i \in [d],
\end{align*}
where $\ell$ ranges over the indices for the $J+2$ components of $p$ and $D\theta$. By strong convexity, the unique minimizers $(\theta^\star)^i$, for $i \in [d]$, of \eqref{eq:mf_coeff_prox} are characterized by
$\nabla F_i((\theta^\star)^i)=0$, which can be solved efficiently using Newton's method \cite[Section~9.5]{boyd2004convex}:
Starting from a strictly feasible point
$\theta_0$, we compute the Newton direction
$\Delta_n$ by solving
\begin{equation}
    \label{eq:newton}
    \nabla^2 F_i(\theta_n)\Delta_n
    =
    -\nabla F_i(\theta_n),
\end{equation}
and update
$\theta_{n+1}=\theta_n+\alpha_n\Delta_n$,
where the step size $\alpha_n\in(0,1]$ is chosen to preserve feasibility.

\begin{remark}[Computational cost of the proximal step]
Since $M_\Psi$ is tridiagonal and $D$
is upper bidiagonal,
the Hessian $\nabla^2 F_i(\theta)$ is tridiagonal.
Consequently, each Newton system in \eqref{eq:newton} can be solved in $O(J)$ arithmetic
operations. Moreover, the $d$ subproblems can be solved in parallel.
\end{remark}

\subsection{Radial Wasserstein FISTA}
\label{sec:Rad-FISTA}

\begin{figure}[t]
    \centering
    \resizebox{0.9\linewidth}{!}{%
        \begin{tikzpicture}[
  x=1.80cm, y=2.38cm,
  basis/.style={line width=1.4pt, line cap=butt, line join=round},
  tail/.style={line width=1.3pt, dash pattern=on 1.3pt off 2.15pt,
               line cap=butt},
  axis/.style={line width=0.65pt, -{Latex[length=2.0mm,width=1.5mm]}},
  every node/.style={inner sep=1.5pt}
]
  \definecolor{basisBlue}{RGB}{40,119,184}
  \definecolor{basisOrange}{RGB}{225,129,39}
  \definecolor{basisGreen}{RGB}{43,155,83}
  \definecolor{basisRed}{RGB}{212,77,69}
  \definecolor{basisPurple}{RGB}{140,91,189}
  \definecolor{basisGray}{gray}{0.46}
  \definecolor{basisAxis}{gray}{0.20}

  \draw[axis, basisAxis] (0,0) -- (8.43,0);
  \draw[axis, basisAxis] (0,-0.22) -- (0,1.40);
  \draw[basisAxis, line width=0.55pt] (-0.045,1) -- (0.045,1);

  \draw[basis, basisBlue]   (0,0) -- (2,1) -- (3,0);
  \draw[basis, basisOrange] (2,0) -- (3,1) -- (4,0);
  \draw[basis, basisGreen]  (3,0) -- (4,1) -- (5,0);
  \draw[basis, basisRed]    (4,0) -- (5,1) -- (6,0);
  \draw[basis, basisPurple] (5,0) -- (6,1) -- (8.20,1);

  \draw[tail, basisGray] (6,0) -- (7.50,1.50);

  \foreach \s in {2,3,4,5,6}
    \fill[basisAxis] (\s,0) circle[radius=2pt];

  \node[basisBlue,   anchor=south] at (2,1.12) {$\psi_0$};
  \node[basisOrange, anchor=south] at (3,1.12) {$\psi_1$};
  \node[basisGreen,  anchor=south] at (4,1.12) {$\psi_2$};
  \node[basisRed,    anchor=south] at (5,1.12) {$\psi_3$};
  \node[basisPurple, anchor=south] at (7,1.12) {$\psi_4$};
  \node[basisGray, anchor=west] at (7.00,0.68) {$\psi_+$};

  \node[anchor=east] at (-0.10,1) {$1$};
  \node[anchor=north east] at (-0.10,-0.055) {$0$};
  \node[anchor=north] at (2,-0.09) {$\sqrt{d}-B$};
  \node[anchor=north] at (4,-0.09) {$\sqrt{d}$};
  \node[anchor=north] at (6,-0.09) {$\sqrt{d}+B$};
  \node[anchor=west] at (8.46,-0.025) {$\|x\|_2$};
\end{tikzpicture}
    }
    \caption{The radial nodal basis $\Lambda$ for $J=4$.
    Only the nonzero values of each basis function
    are plotted. The dots on the horizontal axis
    represent the discretization points
    $b_0<\dots<b_4$.}
    \label{fig:radial_nodal_basis}
\end{figure}
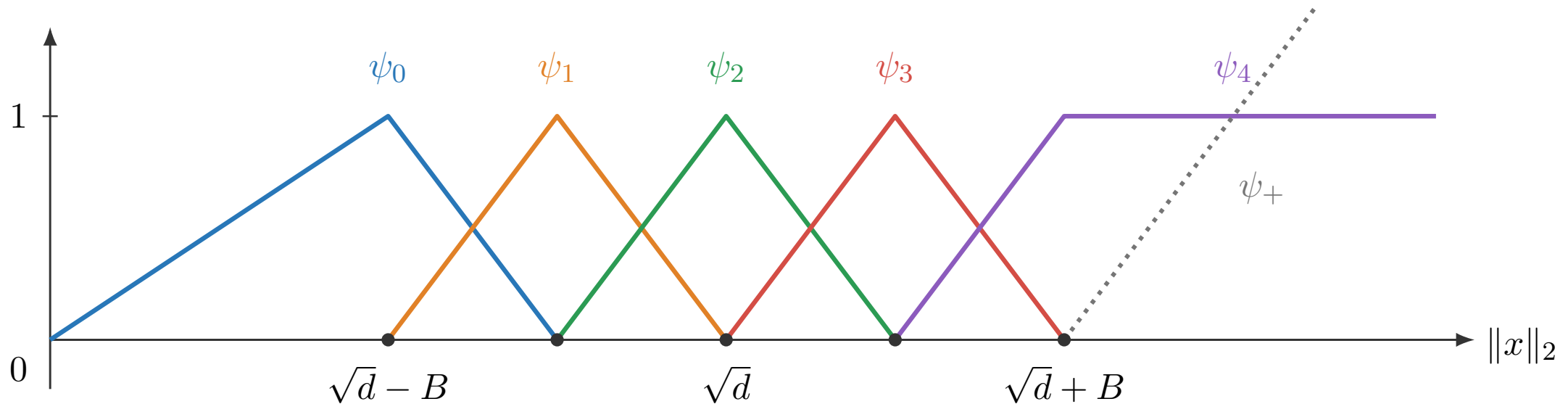

We next consider Problem~\eqref{eq:composite_problem}
with $\sQ=\mathsf{Rad}$ and refer to \texttt{Rad-FISTA} as Algorithm~\ref{alg:W-FISTA} with the
standard choices of transport spaces in
Proposition~\ref{prop:standard_choices}, namely,
\[
    \rho=\cN(0,I_d),\qquad
    \cH=L^2(\rho;\RR^d) \cap \cM_{\mathrm{rad}},\qquad
    \cT=\cH\cap\cM_{\mathrm{cm}}(\rho).
\]
As in the mean-field setting, we consider a nodal
basis representation of the transport maps, now
parametrizing their radial profiles; see \citep{ghafourpour2026variational} for a similar treatment. For $0<B<\sqrt d$ and $J\in\NN_+$, let
$b_j\coloneq\sqrt d-B+j\delta$,
$j\in\{0,\ldots,J\}$, where $\delta\coloneq 2B/J$.
Define the radial nodal basis
$\Lambda \coloneq\{\psi_0,\ldots,\psi_J,\psi_+\}$
on $[0,\infty)$ as in Section~\ref{sec:MF-FISTA},
replacing the grid points $a_j$ by $b_j$,
omitting $\psi_-$, setting
$\psi_+(s)\coloneq(s-b_J)\mathbf{1}_{\{s\geq b_J\}}$,
and replacing $\psi_0$ by
\[
    \psi_0(s)\coloneq
    \begin{cases}
        \dfrac{s}{b_0},&s<b_0,\\[2mm]
        \dfrac{b_1-s}{\delta},&s\in[b_0,b_1],\\[2mm]
        0,&s>b_1.
    \end{cases}
\]
The resulting basis is illustrated in
Figure~\ref{fig:radial_nodal_basis}.
Define the basis of radial transport maps by
\[
    \Phi\coloneq
    \{\phi_\ell:\ell\in\{0,\ldots,J,+\}\}
    \subseteq\cH,
\]
where
$\phi_\ell(x)\coloneq
\psi_\ell(\|x\|_2)x/\|x\|_2$ for $x\neq0$
and $\phi_\ell(0)\coloneq0$. Let $\chi_d\coloneq(\|\cdot\|_2)_\sharp\rho$
denote the chi distribution with $d$ degrees of
freedom. Since $\langle\phi_k,\phi_\ell\rangle_{L^2(\rho)}
    =
    \langle\psi_k,\psi_\ell\rangle_{L^2(\chi_d)}$,
the Gram matrix $M_\Phi$ of $\Phi$ coincides with the Gram
matrix $M_\Lambda$ of $\Lambda$ in $L^2(\chi_d)$.
As in the mean-field setting, $M_\Lambda$ is tridiagonal
due to the local support of the basis functions.

Define the $(J+2)$-dimensional subspace
$\cH_\Phi\coloneq\mathrm{span}(\Phi)$ of $\cH$
and the restrictions of $\cT$ and $\mathsf{Rad}$ by
$\cT_\Phi=\cT\cap\cH_\Phi$ and
$\mathsf{Rad}_\Phi=(\cT_\Phi)_\sharp\rho$,
respectively.
Every $T_\theta\in\cH_\Phi$ has the representation
\[
    T_\theta(x)
    =
    \Psi_\theta(\|x\|_2)\frac{x}{\|x\|_2},
    \quad
    \Psi_\theta(s)=\psi(s)^\top\theta,
\]
for $x\neq0$, with $T_\theta(0)=0$, where we write $\theta = (q_0,\ldots,q_J,s_+)^\top
    \in\RR^{J+2}$ and define $\psi(s) \coloneq 
    (\psi_0(s),\ldots,\psi_J(s),\psi_+(s))^\top$. Here $\Psi_\theta$ is a signed radial profile, and is not necessarily nonnegative. The admissible transport maps are given by
\[
    \cT_\Phi=
    \left\{
        T_\theta\in\cH_\Phi:
        0\leq q_0\leq\cdots\leq q_J,\
        s_+\geq0
    \right\}.
\]
In other words, for every $T_\theta \in \cH_\Phi$, the
piecewise linear radial profile $\Psi_\theta$ 
is parametrized by its function values $q_0,\ldots,q_J$ at $b_0,\ldots,b_J$ and its slope
$s_+$ on $[b_J,\infty)$, respectively.
If $T_\theta\in\cT_\Phi$, the radial profile
$\Psi_\theta$ is, in addition, nonnegative and
nondecreasing.

\textbf{Approximation guarantees for VI.}
We refer to Problem~\eqref{eq:VI_in_sec} with $\sQ=\mathsf{Rad}_\Phi$
as
approximate radial VI. Similarly to the mean-field case, we can refine the basis $\Phi$ to approximate the radial VI minimizer by the approximate radial VI minimizer. Under Assumption~\ref{assm:mean_field_V}, both problems admit unique minimizers
\[
    \pi^\star\coloneq
    \argmin_{\mu\in\mathsf{Rad}}\KL(\mu\Vert\pi),
    \qquad
    \pi_\Phi^\star\coloneq
    \argmin_{\mu\in\mathsf{Rad}_\Phi}\KL(\mu\Vert\pi),
\]
by strong convexity of their lifted objectives. The space $\cT_\Phi$ contains the piecewise affine
radial transport maps of
\citet[Eqs.~(5)--(6)]{ghafourpour2026variational}
when the discretization grids coincide.
Consequently, their approximation guarantees extend
to our setting, as stated in the next proposition.

\begin{proposition}
\label{prop:Rad_approximation}
Let Assumption~\ref{assm:mean_field_V} hold,
suppose additionally that $\nabla V(0)=0$ and
$d\geq3$, and set $\kappa\coloneq L_V/\ell_V$.
Let $\varepsilon\in(0,1]$ be such that
\begin{equation}
    \label{eq:eta}
    \log(d\kappa^3/\varepsilon)\leq cd,
\end{equation}
for a sufficiently small constant $c>0$.
Then there exist
\[
    B=O\left(
        \sqrt{\log\frac{d\kappa^3}{\varepsilon}}
    \right),
    \qquad
    J=\widetilde O\left(
        \frac{\kappa^2}{\sqrt{\varepsilon}}
    \right),
\]
such that the respective radial VI minimizers $\pi^\star$ and $\pi_\Phi^\star$ satisfy
\[
    \frac{\ell_V}{2} W_2^2(\pi_\Phi^\star,\pi^\star)
    \leq
    \cF(\pi_\Phi^\star)-\cF(\pi^\star)
    \lesssim\varepsilon.
\]
\end{proposition}

For fixed $\kappa$, condition~\eqref{eq:eta} holds
for any prescribed $\varepsilon\in(0,1]$ provided
that the dimension $d$ is sufficiently large.

\textbf{Applicability of finite-dimensional \texttt{W-FISTA}.}
As in the mean-field setting, $\Id\in\cH_\Phi$ and
$r(\Id)<\infty$, so approximate radial VI satisfies
Assumption~\ref{assm:lift} by
Proposition~\ref{prop:vi_summary} and the restriction
argument in Section~\ref{sec:finite-W-FISTA}.
We may thus apply Algorithm~\ref{alg:finite-W-FISTA}
to approximate radial VI, which we refer to as finite-dimensional \texttt{Rad-FISTA}.

The deterministic gradient and stochastic gradient
oracle for approximate radial VI are obtained by
orthogonally projecting the corresponding expressions
in Table~\ref{tab:vi_summary} onto $\cH_\Phi$.
Using the Gram-matrix formula for orthogonal
projection, we obtain for $T\in\cH_\Phi$,
\begin{align*}
    \nabla_{\cH_\Phi}g(T)(x)
    &=
    \proj_{\cH_\Phi}[\nabla_\cH g(T)](x)\\
    &=
    \psi(\|x\|_2)^\top M_\Lambda^{-1}
    \EE_{X\sim\rho}\left[
        \psi(\|X\|_2)
        \left\langle
            \nabla V(T(X)),\frac{X}{\|X\|_2}
        \right\rangle
    \right]
    \frac{x}{\|x\|_2}.
\end{align*}
Similarly, taking $(\Xi,\PP)=(\cS^{d-1},\mathrm{Unif}(\cS^{d-1}))$ and letting $\Psi_T$ be the signed radial profile of $T$,
the stochastic gradient oracle is given by
\begin{equation}
\label{eq:rad_oracle_finite_dim}
\begin{aligned}
    \cO_{\cH_\Phi}(T,\xi)(x)
    &\coloneq
    \proj_{\cH_\Phi}[\cO_\cH(T,\xi)](x)\\
    &=
    \psi(\|x\|_2)^\top M_\Lambda^{-1}
    \EE_{s\sim\chi_d}\left[
        \psi(s)\operatorname{sgn}(\Psi_T(s))
        \left\langle
            \nabla V(|\Psi_T(s)|\xi)-\nabla V(0),\xi
        \right\rangle
    \right]
    \frac{x}{\|x\|_2}.
\end{aligned}
\end{equation}
Both expressions are defined to be zero at $x=0$.
As before, $\cO_{\cH_\Phi}$ inherits the noise
parameters of $\cO_\cH$, namely,
$\sigma_\Phi^2=L_V^2$ and
$\zeta_{\cH_\Phi}(T)=\|T\|_{L^2(\rho)}$.

The arguments in Section~\ref{sec:MF-FISTA} apply
equally to radial VI. Thus, provided that the approximation
guarantee in Proposition~\ref{prop:Rad_approximation}
applies at the required accuracy, the outer iteration complexity of
finite-dimensional \texttt{Rad-FISTA} for solving radial VI to accuracy $\varepsilon$ is
$O(\log(1/\varepsilon))$ in the deterministic setting and
$O(\varepsilon^{-1}\log(1/\varepsilon))$ in the stochastic
setting by Theorems~\ref{thm:w-fista_rate}
and~\ref{thm:w-fista_rate_stochastic}, respectively.

\begin{remark}[Comparison to \citealp{ghafourpour2026variational}]
As in Remark~\ref{rem:MF_deterministic_comparison}
for mean-field VI, \citet{ghafourpour2026variational}
restrict the transport maps to a domain on which
the negative entropy is smooth, resulting in a
polynomial dependence on $1/\varepsilon$ in their
deterministic iteration complexity
\citep[Theorem~4.3]{ghafourpour2026variational}.
In contrast, finite-dimensional \texttt{Rad-FISTA}
retains a logarithmic dependence on $1/\varepsilon$
in its deterministic iteration complexity
at the expense of solving a numerically tractable
proximal problem at each iteration.
\end{remark}

\begin{remark}
As discussed in Remark~\ref{rem:MF-VI-approx}
for mean-field VI, approximating the one-dimensional
expectation in \eqref{eq:rad_oracle_finite_dim}
by sampling introduces additional variance.
For their radial parametrization,
\citet{ghafourpour2026variational} derive variance
bounds depending on the number of basis functions.
We leave the corresponding analysis for our basis
with a free-tail slope to future work.
\end{remark}

\textbf{Proximal operator.}
Unlike in mean-field VI, the Jacobian determinant of a radial
transport map contains an additional contribution
from the tangential directions.
For $T_\theta\in\cT_\Phi$ with strictly positive
slopes, we have by \cite[Appendix C.5]{ghafourpour2026variational}, that
\begin{align*}
    \bar{R}(\theta)
    &=
    -\Ent(\rho)
    -\int_0^\infty
    \left[
        \log \Psi_\theta'(s)
        +(d-1)\log\frac{\Psi_\theta(s)}{s}
    \right]\dd\chi_d(s)\\
    &=
    -\Ent(\rho)-p^\top\log(D\theta)
    -(d-1)\int_0^\infty
    \log\frac{\psi(s)^\top\theta}{s}\dd\chi_d(s),
\end{align*}
where $p\coloneq
    \bigl(
        \chi_d((0,b_0)),
        \chi_d((b_0,b_1)),\ldots,
        \chi_d((b_{J-1},b_J)),
        \chi_d((b_J,\infty))
    \bigr)^\top$ and $D\in\RR^{(J+2)\times(J+2)}$ is the matrix
satisfying
\[
    D\theta=
    \left(
        \frac{q_0}{b_0},
        \frac{q_1-q_0}{\delta},\ldots,
        \frac{q_J-q_{J-1}}{\delta},
        s_+
    \right)^\top.
\]
As in mean-field VI, $\bar{R}(\theta)$ is finite if and only if
$D\theta>0$.
Consequently, the proximal problem in
\eqref{eq:finite_coeff_prox} reduces to
\begin{equation}
    \label{eq:rad_coeff_prox}
    \theta^\star
    =
    \argmin_{\theta:\,D\theta>0}  \left\{-p^\top\log(D\theta)
    -(d-1)\int_0^\infty
    \log\frac{\psi(s)^\top\theta}{s}\dd\chi_d(s)
    +\frac{1}{2h}\|\theta-z\|_{M_\Lambda}^2\right\}
\end{equation}
   
Let $F$ denote the objective in \eqref{eq:rad_coeff_prox}. The gradient and Hessian of $F$ are given by
\begin{align*}
    \nabla F(\theta)
    &=
    -D^\top
    \left(\frac{p_\ell}{(D\theta)_\ell}\right)_\ell
    -(d-1)\int_0^\infty
    \frac{\psi(s)}{\psi(s)^\top\theta}\dd\chi_d(s)
    +\frac{1}{h}M_\Lambda(\theta-z),\\
    \nabla^2F(\theta)
    &=
    D^\top\mathrm{diag}
    \left(\frac{p_\ell}{(D\theta)_\ell^2}\right)_\ell D +
    (d-1)\int_0^\infty
    \frac{\psi(s)\psi(s)^\top}
    {(\psi(s)^\top\theta)^2}\dd\chi_d(s)
    +\frac{1}{h}M_\Lambda.
\end{align*}
Again, the unique minimizer of
\eqref{eq:rad_coeff_prox} can be computed using
Newton's method.

\begin{remark}[Computational cost of the proximal step]
Since the Hessian $\nabla^2 F(\theta)$ is again tridiagonal, each Newton system can be solved in $O(J)$ arithmetic operations.
Unlike in mean-field VI, evaluating the gradient and Hessian
of the proximal objective requires one-dimensional
integration, which one can approximate numerically. We leave the subsequent approximation analysis for future work.
\end{remark}

\newpage
\section{Experiments}
\label{sec:experiments}

    We conduct a simple Gaussian VI experiment\footnote{The code is available at \texttt{https://github.com/jupall/wfista}.}, that is Problem  \eqref{eq:VI_in_sec} using $\mathsf{Q}=\mathsf{BW}$, demonstrating the competitiveness of \texttt{BW-FISTA} against the closest work to ours, namely, FB-GVI~\citep{diao2023forward}. The experiment targets the convex and strongly convex regimes separately in the smooth functional. The convex regime is that of the logistic regression; the strongly convex one simply solves the Gaussian VI problem with a quadratic potential. 
    
    \subsection{Quadratic potential}
    For the strongly convex setting, we solve the Gaussian VI problem, with
    \begin{equation}
        V(\theta) = \frac{1}{2}(\theta-c)^\top Q (\theta-c),
    \end{equation}
    where $Q\succ 0$ ensures strong convexity. The exact gradients and Hessians are respectively:
    $$\nabla V(\theta) = Q (\theta - c), \qquad \nabla^2 V(\theta) = Q $$

\begin{figure}[htbp]
        \centering
        \begin{subfigure}[b]{0.48\textwidth}
            \centering
            \includegraphics[width=\textwidth]{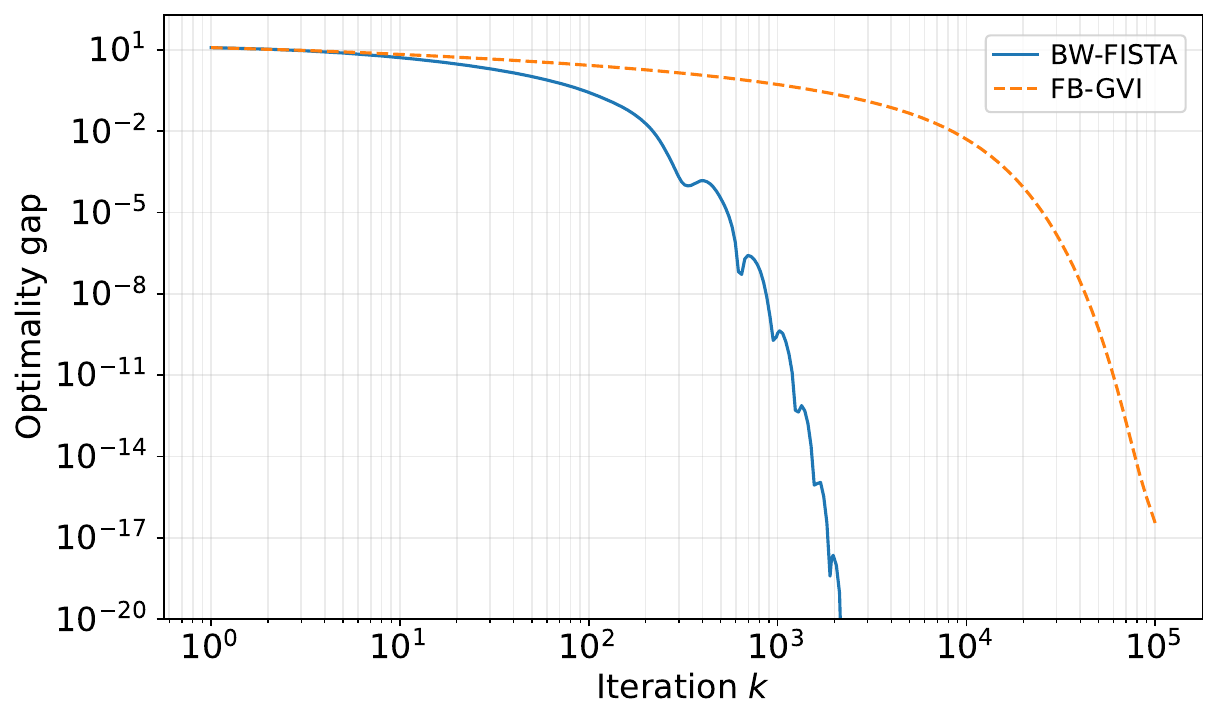}
            \caption{Deterministic oracle}
            \label{fig:optimality_gap_quadratic_deterministic}
        \end{subfigure}
        \hfill
        \begin{subfigure}[b]{0.48\textwidth}
            \centering
            \includegraphics[width=\textwidth]{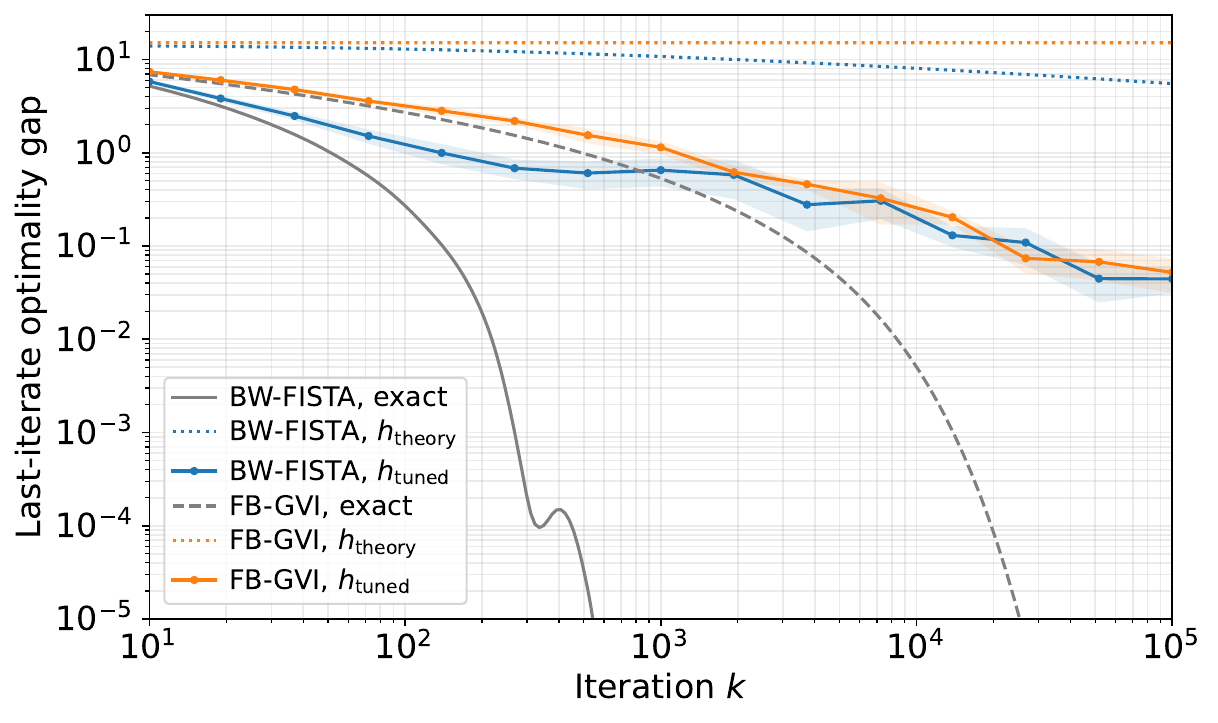}
            \caption{Stochastic oracle, $m=1$}
            \label{fig:optimality_gap_quadratic_m1}
        \end{subfigure}
    
        \vspace{0.5em}
        \begin{subfigure}[b]{0.48\textwidth}
            \centering
            \includegraphics[width=\textwidth]{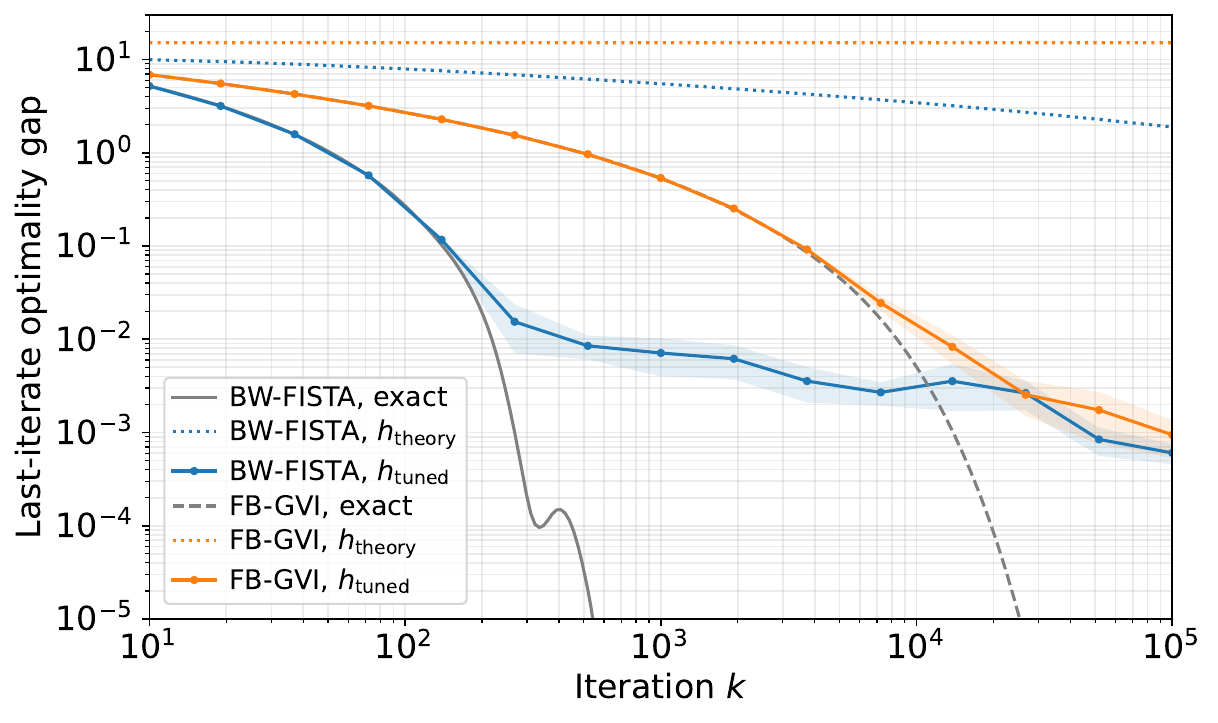}
            \caption{Stochastic oracle, $m=100$}
            \label{fig:optimality_gap_quadratic_m100}
        \end{subfigure}
        \caption{Optimality gap vs.\ iteration for deterministic and stochastic gradient versions of \texttt{BW-FISTA} on the quadratic potential.}
        \label{fig:optimality_gap_quadratic}
    \end{figure}

    We construct $Q = \Gamma \diag(\lambda_1, \ldots, \lambda_d) \Gamma^\top$ with $\Gamma^\top \Gamma = I_d$ and chosen eigenvalues $\lambda_i, i\in[d]$ to control both the smoothness constant $L = \max_{i\in [d]} \lambda_i$ and the strong convexity constant $\alpha = \min_{i\in[d]} \lambda_i$, and by extension the condition number $\kappa = L/\alpha$. In particular, we take $\lambda_j = \kappa^{-(d-j)/(d-1)}, \ \forall j\in [d]$ such that $\alpha=1/\kappa$ and $\beta = 1$. We then impose $d=8$ and, to favor a regime in which acceleration is beneficial, $\kappa=10^4$. The results are presented in Figure \ref{fig:optimality_gap_quadratic}. 
    
    Figure \ref{fig:optimality_gap_quadratic_deterministic} showcases the deterministic oracle. It can be seen that \texttt{BW-FISTA} converges faster than FB-GVI. Figures \ref{fig:optimality_gap_quadratic_m1} and \ref{fig:optimality_gap_quadratic_m100} compare the stochastic gradient version of each method with batch sizes $m=1$ and $m=100$, respectively. We plot rates for step sizes selected as a function of the final iteration count $k$, using both theory-prescribed step sizes and those tuned via a grid search (selecting for lowest median objective gap of the last iterate, over $4$ random seeds). The plotted objective gaps are averaged over 16 additional independent trials. In both cases, \texttt{BW-FISTA} enjoys an initially accelerated convergence rate, as predicted, before entering a noise-dominated regime where the two methods are comparable. As $m$ grows, the stochastic oracle methods align more closely with their exact oracle counterparts.

    \subsection{Logistic regression}

    For the convex case, we solve the logistic regression with a flat prior following the structure of \cite{diao2023forward}.   We observe samples $\{(X_i,Y_i)\}_{i=1}^n\in(\R^d\times\{0,1\})^n$ from the logistic model $Y_i\mid X_i\sim\mathrm{Bern}\bigl(\sigma(\langle\theta^\circ,X_i\rangle)\bigr)$, with $\theta^\circ\sim\mathcal{N}(0,I_d)$, $X_i=\tilde X_i/\lambda_{\max}\bigl(\sum_j\tilde X_j\tilde X_j^\top\bigr)^{1/2}$ the scaled features generated from $\tilde X_1,\dots,\tilde X_n\overset{\text{iid}}{\sim}\mathcal N(0,\diag(1, 1, \ldots, s^2))$, and where $\sigma(t)\defeq e^{t}/(1+e^{t})$. The constant $s >0$ is used to tune the smallest eigenvalue of the Hessian. We consider a uniform improper prior on $\theta\in\R^d$. The posterior on $\theta$ is then $\pi\propto\exp(-V)$, with the potential:
    \begin{equation}\label{eq:exp_potential}
    \begin{aligned}
        V(\theta)
        &=\sum_{i=1}^n\Bigl(\ln\bigl(1+e^{\langle\theta,X_i\rangle}\bigr)-Y_i\langle\theta,X_i\rangle\Bigr) =\sum_{i=1}^n\ln\bigl(1+e^{-(2Y_i-1)\langle\theta,X_i\rangle}\bigr).
    \end{aligned}
    \end{equation}

    The gradient and Hessian of $V$ are
\begin{equation}\label{eq:exp_derivatives}
    \nabla V(\theta)=\sum_{i=1}^n\bigl(\sigma(\langle\theta,X_i\rangle)-Y_i\bigr)\,X_i,\qquad
    \nabla^2V(\theta)=\sum_{i=1}^n\sigma'(\langle\theta,X_i\rangle)\,X_iX_i^\top,
\end{equation}
with $\sigma'=\sigma(1-\sigma)$. Since $0<\sigma'\leq\frac14$, $V$ is convex and $L$-smooth with $L=\frac14\lambda_{\max}\bigl(\sum_iX_iX_i^\top\bigr)=\frac14$. The posterior $\pi$ is a probability measure, and $\mu^\star$ exists only if the data are not linearly separable, which a linear program certifies for our data. Finally, consider $\kappa_{\mathrm{loc}}\defeq L/\lambda_{\min}$ with $\lambda_{\min}$ the smallest eigenvalue of the Hessian at $\mu^\star=\mathcal N(b^\star,\Sigma^\star)$, we choose $s$ to fix $\kappa_{\mathrm{loc}} \approx 10^4$. In this experiment, $d=10$ and $n=20$.  Due to the decoupled structure of the regression, the expectations reduce to one-dimensional Gaussian integrals over the marginals, which we compute by composite Gauss–Legendre quadrature~\citep{davis2007methods} to machine precision \citep{trefethen2008gauss}.

\begin{figure}[htbp]
    \centering
    \begin{subfigure}[b]{0.48\textwidth}
        \centering
        \includegraphics[width=\textwidth]{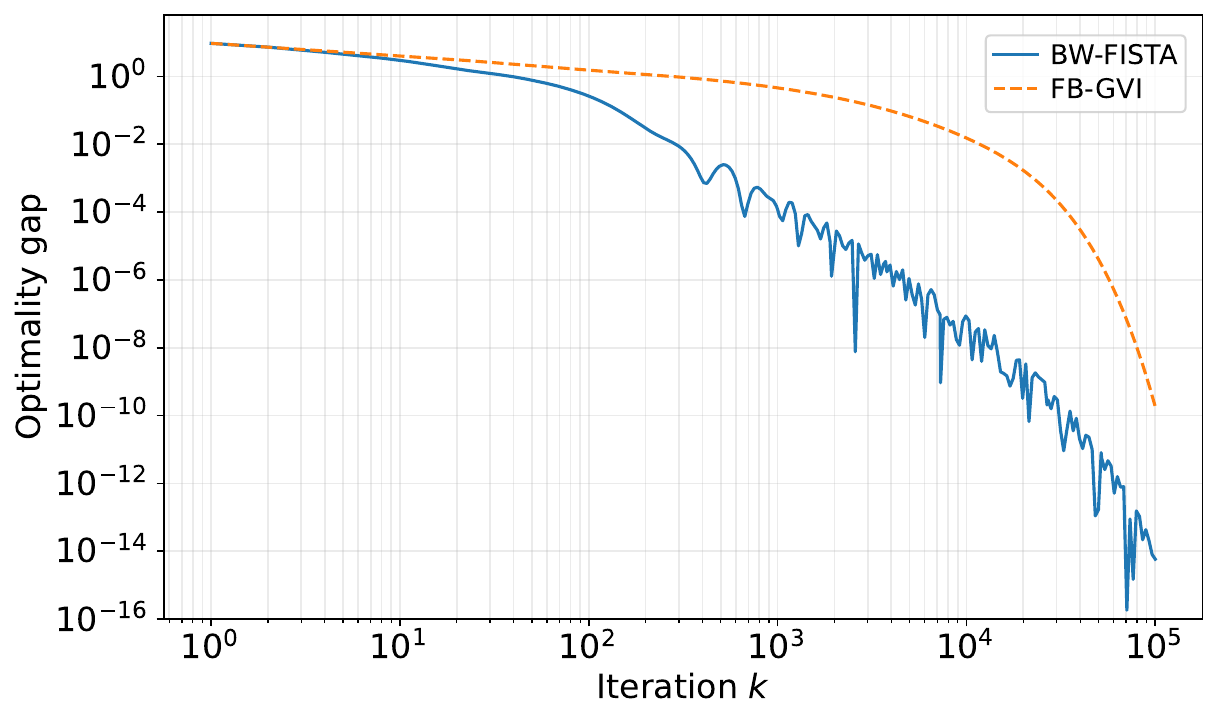}
        \caption{Deterministic oracle}
        \label{fig:optimality_gap_deterministic}
    \end{subfigure}
    \hfill
    \begin{subfigure}[b]{0.48\textwidth}
        \centering
        \includegraphics[width=\textwidth]{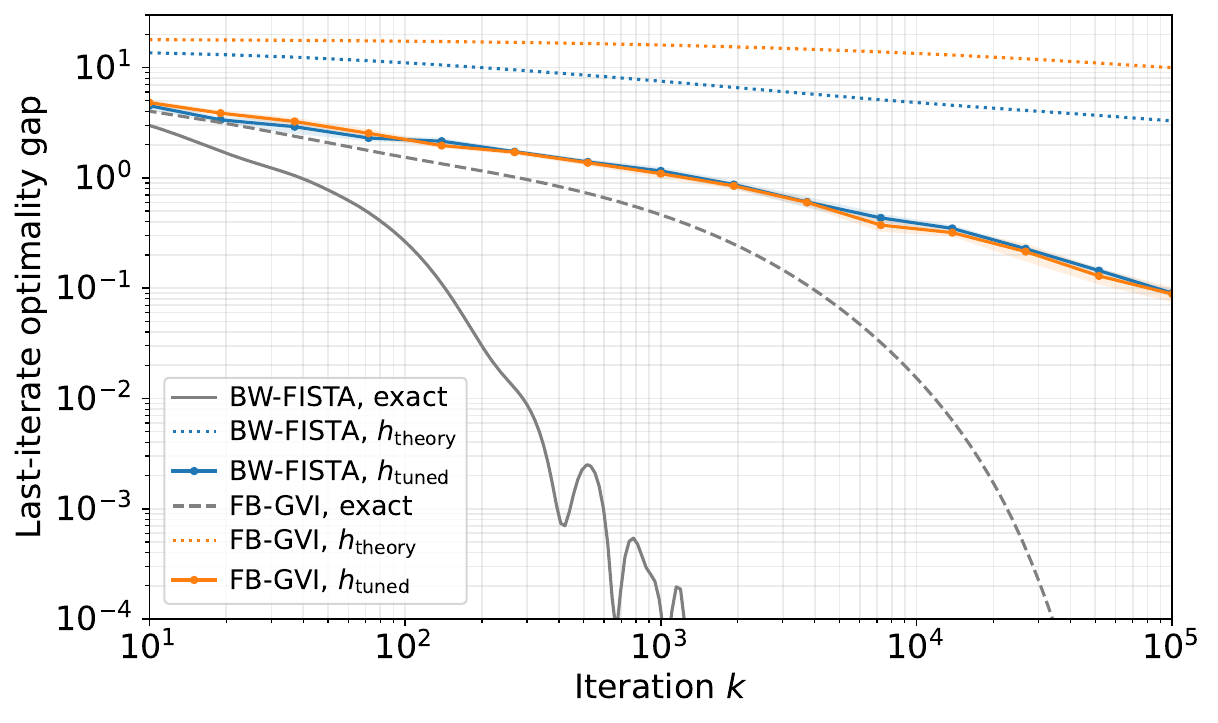}
        \caption{Stochastic oracle, $m=1$}
        \label{fig:optimality_gap_m1}
    \end{subfigure}

    \vspace{0.5em}
    \begin{subfigure}[b]{0.48\textwidth}
        \centering
        \includegraphics[width=\textwidth]{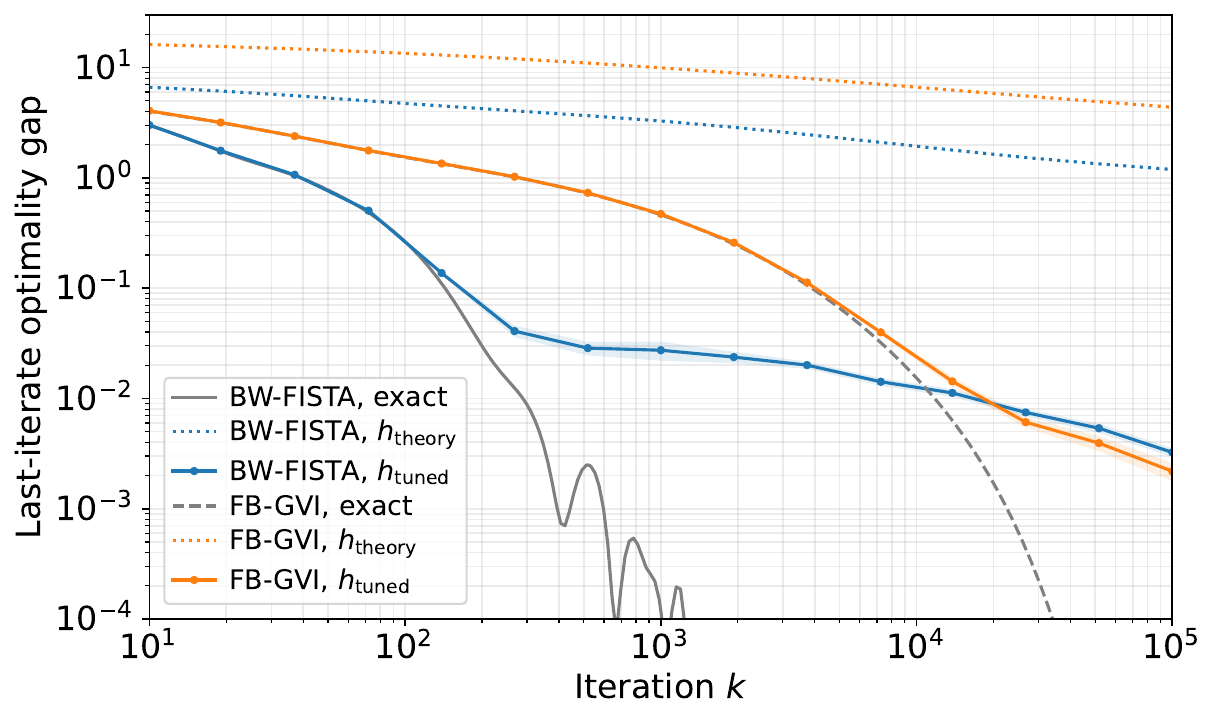}
        \caption{Stochastic oracle, $m=100$}
        \label{fig:optimality_gap_m100}
    \end{subfigure}
    \caption{Optimality gap vs.\ iteration for deterministic and stochastic gradient versions of \texttt{BW-FISTA} on the logistic regression.}
    \label{fig:optimality_gap_logistic}
\end{figure}

Results are presented in Figure \ref{fig:optimality_gap_logistic}. We observe in the deterministic case that strong convexity holds close to the minimizer, as the tuned step-size of FB-GVI infers an exponential convergence close to the minimizer. Nonetheless, the \texttt{BW-FISTA} profits from a faster rate on the finite horizon. Similarly to the previous quadratic experiment, Figures \ref{fig:optimality_gap_m1} and \ref{fig:optimality_gap_m100} compare the stochastic gradient version of each method with batch sizes $m=1$ and $m=100$, respectively. With $m=1$, both methods converge at the same speed. Just like in the quadratic experiment, as $m$ grows, the stochastic oracle methods align more closely with their exact oracle counterparts in the beginning of the experiment before being slowed down by noise. 

\section*{Acknowledgements}

This work was supported as a part of NCCR Automation, a National Centre of Competence in Research, funded by the Swiss National Science Foundation (grant number 51NF40\_225155). The authors acknowledge OpenAI for providing complimentary access to ChatGPT through the ChatGPT for Academic Researchers program.\newpage

\bibliographystyle{plainnat}
\bibliography{ref_light}

\etocdepthtag.toc{appendix}

\appendix

\newpage
\etocsettagdepth{mainmatter}{none}
\etocsettagdepth{appendix}{subsection}
\renewcommand{\contentsname}{Appendix Contents}
\tableofcontents

\section{Proofs for Section~\ref{sec:prelims}}
\label{app:prelims}

\subsection{Proof of Proposition~\ref{prop:pushforward-expressivity}}
\label{prf:pushforward-expressivity}

Fix $\sQ$ and $\rho \in \sQ \cap \cP_{2,\mathrm{ac}}(\R^d)$
as in the statement. First, we note that the inclusion
\begin{equation*}
    (L^2(\rho;\R^d) \cap \cM_\sQ)_\sharp \rho \subseteq \sQ
\end{equation*}
is direct in all three cases. Indeed, the pushforward of a Gaussian
by an affine map is Gaussian, the pushforward of a product measure
by a coordinate-wise separable map is a product measure, and the
pushforward of a radially symmetric measure by a radial map is
radially symmetric. For radial maps, this follows from
$U_\sharp(T_\sharp\rho)
=T_\sharp(U_\sharp\rho)=T_\sharp\rho$. Moreover, the intersection with $L^2(\rho;\R^d)$ ensures
finite second moments.

We next prove that each $\mu\in\sQ$ admits an optimal
transport map from $\rho$ belonging to $\cM_\sQ$, and that
the cyclically monotone maps in $\cM_\sQ$ are exactly those
admitting representatives with the conditions listed in
the statement.

For the first claim, since $\rho$ is absolutely continuous,
Brenier's theorem gives a unique optimal map
$T^{\rho\to\mu}\in\cM_{\mathrm{cm}}(\rho)$ for every
$\mu\in\sQ$.
In each case below, we construct a representative of this
map in $\cM_\sQ$ satisfying the listed conditions.
This proves
\[
    \sQ
    \subseteq
    (\cM_{\mathrm{cm}}(\rho)\cap\cM_\sQ)_\sharp\rho
    \subseteq
    (L^2(\rho;\R^d)\cap\cM_\sQ)_\sharp\rho,
\]
which, together with the reverse inclusion proved above,
establishes the pushforward identity.

For the characterization, these same constructions show
that every cyclically monotone map in $\cM_\sQ$ admits
a representative with the listed conditions.
Indeed, if $T\in\cM_{\mathrm{cm}}(\rho)\cap\cM_\sQ$,
then $\mu\defeq T_\sharp\rho$ belongs to $\sQ$, and
optimality and uniqueness give $T=T^{\rho\to\mu}$
$\rho$-a.e.; see, for example, 
\citealp[Theorems~1.22 and~1.48]{santambrogio2015optimal}.
The construction for this $\mu$ therefore provides the
required representative.
To complete the characterization, each case ends by
checking that every square-integrable map with the listed
conditions is cyclically monotone.\medskip

\textbf{Gaussian measures ($\sQ = \BW$).}
Since $\rho$ is absolutely continuous and Gaussian, we have
$\Sigma_\rho \succ 0$. Fixing $\mu \in \BW$ and defining
\begin{equation*}
    A \defeq
    \Sigma_\rho^{-1/2}
    \bigl(
        \Sigma_\rho^{1/2}
        \Sigma_\mu
        \Sigma_\rho^{1/2}
    \bigr)^{1/2}
    \Sigma_\rho^{-1/2}
    \in \PSD,
\end{equation*}
we have
\begin{equation*}
    T^{\rho \to \mu}(x)
    =
    A(x-m_\rho)+m_\mu
    \qquad \rho\text{-a.e.}
\end{equation*}
This is the standard formula for Gaussian optimal transport
(see, e.g., Section~1.6.3 of
\citealp{panaretos2020invitation}). Thus $T^{\rho \to \mu}$ admits a representative in
$\cM_{\mathrm{aff}}(\PSD)
\subseteq \cM_{\mathrm{aff}}(\Sym)$, as desired.

To complete the characterization, suppose that $T\in L^2(\rho;\R^d)$
satisfies $T(x)=Ax+b$ $\rho$-a.e.\ for some $A\in\PSD$.
Then $T\in\cM_{\mathrm{cm}}(\rho)\cap\cM_{\mathrm{aff}}(\Sym)$,
since it is the $\rho$-a.e.\ gradient of the convex function
$x\mapsto\frac12 x^\top A x+b^\top x$.\medskip

\textbf{Product measures ($\sQ = \MF$).}
For a probability measure $\eta$ on $\R$, define its cumulative
distribution and quantile functions by
\begin{equation*}
    F_\eta(t) \defeq \eta((-\infty,t]),
    \quad t \in \R,
    \qquad
    Q_\eta(u) \defeq
    \inf\{t \in \R:F_\eta(t)\ge u\},
    \quad u \in (0,1).
\end{equation*}
Now fix $\mu \in \MF$ and decompose
\begin{equation*}
    \rho=\rho_1\otimes\cdots\otimes\rho_d,
    \qquad
    \mu=\mu_1\otimes\cdots\otimes\mu_d.
\end{equation*}
For each $i=1,\ldots,d$, $\rho_i$ is absolutely
continuous and hence atomless, and
$\rho_i,\mu_i \in \cP_2(\R)$. Therefore, $T^{\rho_i \to \mu_i} = Q_{\mu_i}\circ F_{\rho_i}$ $\rho_i$-a.e.,
and this map is non-decreasing $\rho_i$-a.e. (see, e.g.,
Theorems~2.5 and~2.9 of \citealp{santambrogio2015optimal}).
The formula is finite on the Borel set
$\{t \in \R:0<F_{\rho_i}(t)<1\}$, which has full
$\rho_i$-measure. Assigning the value $0$ elsewhere gives a
finite Borel representative.

Using these representatives, consider the separable Borel map
\begin{equation*}
    T(x_1,\ldots,x_d)
    \defeq
    \bigl(
        T^{\rho_1 \to \mu_1}(x_1),
        \ldots,
        T^{\rho_d \to \mu_d}(x_d)
    \bigr).
\end{equation*}
The product structure implies that $T_\sharp\rho=\mu$.
In particular, $T \in L^2(\rho;\R^d)$ because
$\mu \in \cP_2(\R^d)$. This certifies that
\begin{align*}
    W_2(\rho,\mu)^2
    &\le
    \E_\rho[\|T-\Id\|^2] =
    \sum_{i=1}^d
    \E_{\rho_i}
    \bigl[
        |T^{\rho_i \to \mu_i}-\Id|^2
    \bigr]
    =
    \sum_{i=1}^d W_2(\rho_i,\mu_i)^2.
\end{align*}
Moreover, for any coupling $\pi \in \Pi(\rho,\mu)$,
\begin{equation*}
    \int\|x-y\|^2\,d\pi(x,y)
    =
    \sum_{i=1}^d
    \int|x_i-y_i|^2\,d\pi(x,y)
    \ge
    \sum_{i=1}^d W_2(\rho_i,\mu_i)^2.
\end{equation*}
Thus, $W_2(\rho,\mu)^2
    =
    \sum_{i=1}^d W_2(\rho_i,\mu_i)^2$
and uniqueness of the optimal map gives
$T^{\rho \to \mu}=T$ $\rho$-a.e. In particular,
$T^{\rho \to \mu}$ admits a representative in
$\cM_{\mathrm{sep}}$, as claimed.

To complete the characterization, suppose that $T\in L^2(\rho;\R^d)$
admits a separable Borel representative whose components
$T_i$ are non-decreasing $\rho_i$-a.e. Set
\begin{equation*}
    \nu_i \defeq (T_i)_\sharp\rho_i,
    \qquad
    \nu \defeq T_\sharp\rho
    =
    \nu_1\otimes\cdots\otimes\nu_d.
\end{equation*}
Each $T_i$ is a non-decreasing transport map from $\rho_i$ to $\nu_i$
on a set of full $\rho_i$-measure, and hence is optimal.
The identity just proved therefore gives
\begin{equation*}
    \E_\rho[\|T-\Id\|^2]
    =
    \sum_{i=1}^d
    \E_{\rho_i}[|T_i-\Id|^2]
    =
    \sum_{i=1}^d W_2(\rho_i,\nu_i)^2
    =
    W_2(\rho,\nu)^2.
\end{equation*}
Thus $T$ is optimal from $\rho$ to $\nu$, and so
$T \in \cM_{\mathrm{cm}}(\rho)\cap\cM_{\mathrm{sep}}$.\medskip

\textbf{Radially symmetric measures ($\sQ = \Rad$).}
Fix $\mu \in \Rad$ and write
\begin{equation*}
    \rho_{\mathrm r} \defeq \|\cdot\|_\sharp\rho,
    \qquad
    \mu_{\mathrm r} \defeq \|\cdot\|_\sharp\mu.
\end{equation*}
These measures belong to $\cP_2(\R)$ and are supported on
$[0,\infty)$. Since every sphere has zero Lebesgue measure,
absolute continuity of $\rho$ implies that $\rho_{\mathrm r}$
is atomless and $\rho_{\mathrm r}(\{0\})=0$.

Using the quantile notation above, define
\begin{equation*}
    q(r)
    \defeq
    Q_{\mu_{\mathrm r}}(F_{\rho_{\mathrm r}}(r))
    \qquad \rho_{\mathrm r}\text{-a.e.}
\end{equation*}
The formula is finite, non-negative, and non-decreasing on the
$\rho_{\mathrm r}$-full Borel set
$\{r>0:0<F_{\rho_{\mathrm r}}(r)<1\}$.
Set $q=0$ outside this set, and define the radial Borel map
\begin{equation*}
    T(x)
    \defeq
    \begin{cases}
        \displaystyle
        q(\|x\|)\frac{x}{\|x\|}, & x\ne0,\\[2mm]
        0, & x=0.
    \end{cases}
\end{equation*}

To check the pushforward identity, write $X=RU\sim\rho$, where
$R\sim\rho_{\mathrm r}$ and $U$ is independent of $R$ and uniform
on $\cS^{d-1}$. The one-dimensional transport formula gives
$q(R)\sim\mu_{\mathrm r}$. Since a radially symmetric measure is
determined by its radius law, $q(R)U\sim\mu$. Thus, $T_\sharp\rho=\mu$, and in particular
$T \in L^2(\rho;\R^d)\cap\cM_{\mathrm{rad}}$.

The reverse triangle inequality implies that every
$\pi \in \Pi(\rho,\mu)$ satisfies
\begin{equation*}
    \int\|x-y\|^2\,d\pi(x,y)
    \ge
    \int(\|x\|-\|y\|)^2\,d\pi(x,y)
    \ge
    W_2(\rho_{\mathrm r},\mu_{\mathrm r})^2.
\end{equation*}
On the other hand, the proposed map attains this lower bound:
\begin{equation*}
    \E_\rho[\|T-\Id\|^2]
    =
    \int(r-q(r))^2\,d\rho_{\mathrm r}(r)
    =
    W_2(\rho_{\mathrm r},\mu_{\mathrm r})^2.
\end{equation*}
Consequently, $W_2(\rho,\mu)^2
    =
    W_2(\rho_{\mathrm r},\mu_{\mathrm r})^2$, and uniqueness gives $T^{\rho \to \mu}=T$.
Thus $T^{\rho \to \mu}$ admits a representative in
$\cM_{\mathrm{rad}}$ with a non-negative, non-decreasing
$\rho_{\mathrm r}$-a.e.\ profile. This is the standard formula for radial OT; see Example~3.5(1), Equation~(3.6), of
\citealp{cuestaalbertos1993optimal}.

To complete the characterization, suppose that $T\in L^2(\rho;\R^d)$
admits a radial Borel representative whose profile
$\Psi_T$ is non-negative and non-decreasing
$\rho_{\mathrm r}$-a.e. Set
\begin{equation*}
    \nu \defeq T_\sharp\rho,
    \qquad
    \nu_{\mathrm r} \defeq \|\cdot\|_\sharp\nu.
\end{equation*}
Since $\Psi_T\ge0$ $\rho_{\mathrm r}$-a.e., $\nu_{\mathrm r}
    =
    (\Psi_T)_\sharp\rho_{\mathrm r}.$
Moreover, $\Psi_T$ is non-decreasing
$\rho_{\mathrm r}$-a.e., so it is an optimal transport from
$\rho_{\mathrm r}$ to $\nu_{\mathrm r}$.
The identity just proved yields
\begin{align*}
    \E_\rho[\|T-\Id\|^2]
    &=
    \int(r-\Psi_T(r))^2\,d\rho_{\mathrm r}(r) =
    W_2(\rho_{\mathrm r},\nu_{\mathrm r})^2
    =
    W_2(\rho,\nu)^2.
\end{align*}
Thus, $T$ is optimal from $\rho$ to $\nu$, and so
$T \in \cM_{\mathrm{cm}}(\rho)\cap\cM_{\mathrm{rad}}$.\qed

\subsection{Proof of Proposition~\ref{prop:main-tangent-spaces}}
\label{prf:main-tangent-spaces}

We first establish a useful helper lemma.

\begin{lemma}[Optimal couplings for the main families]
\label{lem:structured-optimal-couplings}
Let $\sQ \in \{\BW,\MF,\Rad\}$ and $\mu,\nu\in\sQ$. Then there exists $\pi\in\Pi_{\mathrm{opt}}^\sQ(\mu,\nu)$ with the
following structure:
\begin{enumerate}[label=(\roman*), nosep]
    \item If $\sQ=\BW$, then $\pi$ is jointly Gaussian;
    \item If $\sQ=\MF$, writing
    $\mu=\bigotimes_{i=1}^d\mu_i$ and
    $\nu=\bigotimes_{i=1}^d\nu_i$, the coordinate pairs
    $(X_i,Y_i)$ under $(X,Y)\sim\pi$ are independent, with
    $\operatorname{Law}(X_i,Y_i)\in
    \Pi_{\mathrm{opt}}(\mu_i,\nu_i)$;
    \item If $\sQ=\Rad$, then $\pi$ is the law of $(RU,SU)$,
    where $(R,S)$ optimally couples the radius laws
    $\|\cdot\|_\sharp\mu$ and $\|\cdot\|_\sharp\nu$, and
    $U$ is independent of $(R,S)$ and uniform on $\cS^{d-1}$.
\end{enumerate}
\end{lemma}

\begin{proof}
For each case, we construct a suitable optimal coupling and verify that its induced
displacement interpolation remains in $\sQ$.\medskip

\textbf{Gaussian measures.}
Choose any $\pi_0\in\Pi_{\mathrm{opt}}(\mu,\nu)$ and let
$\pi \defeq\cN(m_{\pi_0},\Sigma_{\pi_0})$ be the Gaussian measure on
$\R^{2d}$ with the same mean and covariance.
Since $\mu$ and $\nu$ are Gaussian, $\pi$ has the same marginals
as $\pi_0$. Its transport cost also agrees with that of $\pi_0$,
since the quadratic cost depends only on first and second moments.
Thus $\pi$ is optimal, and its displacement interpolation remains
Gaussian because linear images of jointly Gaussian vectors are
Gaussian.\medskip

\textbf{Product measures.}
Choose independent pairs $(X_i,Y_i)$ with laws
$\pi_i\in\Pi_{\mathrm{opt}}(\mu_i,\nu_i)$, and let $\pi$ be the
joint law of $X \defeq (X_1,\ldots,X_d)$ and $Y \defeq (Y_1,\ldots,Y_d)$.
Then $\pi\in\Pi(\mu,\nu)$ and, for all $\gamma\in\Pi(\mu,\nu)$,
\[
    \int\|x-y\|^2\dd\gamma(x,y)
    =
    \sum_{i=1}^d\int|x_i-y_i|^2\dd\gamma(x,y)
    \ge
    \sum_{i=1}^d W_2(\mu_i,\nu_i)^2
    =
    \E[\|X-Y\|^2].
\]
Thus $\pi$ is optimal. The coordinates of $(1-t)X+tY$ remain
independent, so its law remains in $\MF$ for every $t\in[0,1]$.
\medskip

\textbf{Radially symmetric measures.}
Write $\mu_{\mathrm r} \defeq \|\cdot\|_\sharp\mu$ and
$\nu_{\mathrm r} \defeq \|\cdot\|_\sharp\nu$, and choose $(R,S,U)$
as in the statement. Radial symmetry gives $RU\sim\mu$ and
$SU\sim\nu$, so their joint law $\pi$ belongs to $\Pi(\mu,\nu)$.
For every $\gamma\in\Pi(\mu,\nu)$, the reverse triangle
inequality gives
\begin{align*}
    \int\|x-y\|^2\dd\gamma(x,y)
    &\ge
    \int(\|x\|-\|y\|)^2\dd\gamma(x,y)\\
    &\ge
    W_2(\mu_{\mathrm r},\nu_{\mathrm r})^2
    =
    \E[(R-S)^2]
    =
    \E[\|RU-SU\|^2].
\end{align*}
Thus $\pi$ is optimal. Its displacement interpolation is the
law of $((1-t)R+tS)U$, which is radially symmetric for every
$t\in[0,1]$.
\end{proof}

We now prove Proposition~\ref{prop:main-tangent-spaces}. Let $V_\mu^\sQ$ denote the corresponding space on the
right-hand side of the proposition. We prove $\Tan_\mu\sQ = V_\mu^\sQ$ by establishing the two inclusions
\[
    \Tan_\mu\sQ\subseteq V_\mu^\sQ
    \qquad\text{and}\qquad
    V_\mu^\sQ\subseteq\Tan_\mu\sQ.
\]

For $\sQ=\cP_2(\R^d)$, these follow from
\cite[Proposition~8.5.4 and Corollary~8.5.5]{ambrosio2008gradient};
every optimal barycentric displacement belongs to the closure of smooth
gradients, and positive rescalings of optimal map displacements are dense
in that closure. In what remains, we thus consider the three subfamilies.\medskip

\textbf{The inclusion $\Tan_\mu\sQ\subseteq V_\mu^\sQ$.}
We first note that each proposed space $V_\mu^\sQ$ is a closed linear
subspace. The Gaussian space is finite-dimensional, and the product
space is isometric to $\bigoplus_{i=1}^d L^2(\mu_i)$. The radial space
is isometric to $L^2(\alpha|_{(0,\infty)})$, for the radius law
$\alpha=\|\cdot\|_\sharp\mu$, through the correspondence
\begin{equation}
\label{eq:tangent-radial-profile}
    \psi \in L^2(\alpha|_{(0,\infty)}) \longmapsto v_\psi \in \cM_{\mathrm{rad}} \cap L^2(\mu;\R^d),
    \qquad
    v_\psi(x)
    =
    \begin{cases}
        \displaystyle
        \psi(\|x\|)\frac{x}{\|x\|}, & x\ne0,\\[5pt]
        0, & x=0.
    \end{cases}
\end{equation}
Since $V_\mu^\sQ$ is closed and linear, to prove
$\Tan_\mu\sQ\subseteq V_\mu^\sQ$, it suffices to show that
each admissible barycentric displacement belongs to $V_\mu^\sQ$.

Fix $\nu\in\sQ$ and $\pi\in\Pi_{\mathrm{opt}}^\sQ(\mu,\nu)$. Consider the geodesic midpoint and associated coupling
\[
    \eta \defeq
    \left(\frac{x+y}{2}\right)_\sharp\pi,
    \qquad
    \gamma \defeq
    \left(x,\frac{x+y}{2}\right)_\sharp\pi.
\]
The definition of $\Pi_{\mathrm{opt}}^\sQ(\mu,\nu)$ implies that
$\eta\in\sQ$, while
\citet[Lemma~7.2.1]{ambrosio2008gradient} gives
$\Pi_{\mathrm{opt}}(\mu,\eta)=\{\gamma\}$, even for singular $\mu$.
Lemma~\ref{lem:structured-optimal-couplings} supplies an optimal
coupling of $\mu$ and $\eta$ which satisfies the structural assumptions stated there, and uniqueness implies that this coupling is $\gamma$.
Since $v_\pi=2v_\gamma$, it suffices to show that
$v_\gamma\in V_\mu^\sQ$.\medskip

\textbf{Gaussian measures.}
Set $E \defeq\operatorname{ran}\Sigma_\mu$. If $E=\{0\}$, then $\mu$ is a Dirac mass, $v_\gamma$
is constant, and the conclusion is immediate. Otherwise, since $\gamma$ is jointly Gaussian, there exists a linear map $B:E \to \R^d$ such that
\[
    b_\gamma(x)
    \defeq
    \int z \dd\gamma(z | x)
    =
    m_\eta+B(x-m_\mu).
\]
Since $\gamma \in \Pi_{\mathrm{opt}}(\mu,\eta)$, there exists a proper convex function
$\varphi$ such that $b_\gamma(x)\in\partial\varphi(x)$ $\mu$-a.e. (see \citealp[proof of Theorem~12.4.4]{ambrosio2008gradient}).
Indeed, $\gamma(\cdot| x)$ is concentrated on the closed convex
set $\partial\varphi(x)$, so its mean belongs to the same set.

Define $q:E\to\R\cup\{+\infty\}$ by
$q(\xi)=\varphi(m_\mu+\xi)$.
Since $\mu$ has a strictly positive Gaussian density on $m_\mu+E$,
the convex domain $\dom(q)$ has full Lebesgue measure in $E$ and
therefore equals $E$. Thus $q$ is finite and convex on $E$.
Projecting the subgradient inclusion onto $E$, we obtain
\[
    \nabla_E q(\xi)
    =
    P_E m_\eta+(P_EB)\xi
    \qquad\text{for Lebesgue-a.e.\ }\xi\in E.
\]
Thus, $A\defeq P_EB$ is symmetric, being the linear part of an affine gradient. Relative to $E\oplus E^\perp$, write
\[
    B\defeq
    \begin{pmatrix}
        A\\ C
    \end{pmatrix},
    \qquad
    S\defeq
    \begin{pmatrix}
        A-I_E & C^\top\\
        C & 0
    \end{pmatrix}.
\]
Then, $S$ is symmetric and, for $\mu$-almost all $x \in \R^d$, we have
\begin{align*}
    v_\gamma(x) &= b_\gamma(x) - x\\
    &= m_\eta - m_\mu +B(x-m_\mu)-(x-m_\mu)\\
    &= m_\eta-m_\mu + S(x-m_\mu).
\end{align*}
Thus, $v_\gamma\in \cM_{\mathrm{aff}}(\Sym) \cap L^2(\mu;\R^d) =  V_\mu^{\mathsf{BW}}$.\medskip

\textbf{Product measures.}
Write $\mu=\bigotimes_{i=1}^d\mu_i$,
$\eta=\bigotimes_{i=1}^d\eta_i$, and
$\gamma=\bigotimes_{i=1}^d\gamma_i$, where
$\gamma_i\in\Pi_{\mathrm{opt}}(\mu_i,\eta_i)$.
The product structure gives
\[
    (v_\gamma)_i(x)
    =
    \int_{\R}(z_i-x_i)\,\dd\gamma_i(z_i\mid x_i).
\]
Thus,
$v_\gamma\in\cM_{\mathrm{sep}}\cap L^2(\mu;\R^d)
=V_\mu^{\mathsf{MF}}$.\medskip

\textbf{Radially symmetric measures.}
Write $\gamma=\operatorname{Law}(RU,SU)$ for random variables $(R,S,U)$ as in
Lemma~\ref{lem:structured-optimal-couplings}.
Its barycentric displacement is $v_\psi$ from
\eqref{eq:tangent-radial-profile}, with
\[
    \psi(r)=\mathbb E[S\mid R=r]-r,
    \qquad r>0.
\]
At the origin, the barycentric displacement is zero because
$U$ is independent of $(R,S)$ and has mean zero. Choosing a finite
Borel version of the conditional expectation thus gives
$v_\gamma\in V_\mu^{\mathsf{Rad}}$.\medskip

For each family $\sQ$, $v_\pi=2v_\gamma\in V_\mu^\sQ$. Since $V_\mu^\sQ$
is closed and linear, $ \Tan_\mu\sQ\subseteq V_\mu^\sQ$, as desired.\medskip

\textbf{The inclusion $V_\mu^\sQ\subseteq \Tan_\mu \sQ$.}
We provide a dense class of fields in $V_\mu^\sQ$, each of which is a
positive multiple of an admissible barycentric displacement.
Consider the following classes:
\begin{align*}
    \mathsf{BW}:&\quad
    v(x)=Sx+b,
    \qquad S=S^\top;\\
    \mathsf{MF}:&\quad
    v(x)=(f_1(x_1),\ldots,f_d(x_d)),
    \qquad f_i\in C_c^\infty(\R);\\
    \mathsf{Rad}:&\quad
    v=v_\psi,
    \qquad \psi\in C_c^\infty((0,\infty)).
\end{align*}
The Gaussian class is all of $V_\mu^{\mathsf{BW}}$. In the other two
cases, density follows from regularity of finite Borel measures and
smooth approximation of compactly supported continuous functions.

Every field in these classes is the gradient of a smooth function
$f$ with bounded Hessian. For the radial case, one can take
\[
    f(x)=-\int_{\|x\|}^{\infty}\psi(s)\,\dd s,
\]
which is smooth because it is constant near the origin.
For sufficiently small $\varepsilon>0$,
\[
    T_\varepsilon
    \defeq
    \id+\varepsilon v
    =
    \nabla\left(\frac12\|\cdot\|^2+\varepsilon f\right)
\]
is the gradient of a convex function, for all cases. Thus $\pi_\varepsilon
    \defeq
    (\id,T_\varepsilon)_\sharp\mu$ is an optimal coupling between $\mu$ and $(T_\eps)_\sharp \mu$.
For $t \in [0,1]$, the maps $\id+t\varepsilon v$ are respectively affine,
coordinate-wise separable, or orthogonally equivariant. Hence, their
pushforwards remain in the corresponding family $\sQ$, and so
\[
    \pi_\varepsilon
    \in
    \Pi_{\mathrm{opt}}^\sQ
    \bigl(\mu,(T_\varepsilon)_\sharp\mu\bigr),
    \qquad
    v=\varepsilon^{-1}v_{\pi_\varepsilon}
    \in \Tan_\mu \sQ.
\]
Taking $L^2(\mu)$ closures gives
$V_\mu^\sQ\subseteq\Tan_\mu\sQ$, establishing that $V_\mu^\sQ = \Tan_\mu\sQ$.\qed

\subsection{Proof of Uniqueness for \texorpdfstring{$\sQ$}{Q}-Gradients (Definition~\ref{def:W2-diff})}
\label{app:W2-diff}

Uniqueness of the gradient in the ambient case $\sQ=\cP_2(\R^d)$ is standard; see, e.g., Definition~4.9 and the discussion following equation~(4.5) in \cite{gangbo2011differential}.
We give a proof that applies to an arbitrary family
$\sQ\subseteq\cP_2(\R^d)$.

\begin{proposition}
\label{prop:Q-gradient-uniqueness}
Let $\sQ\subseteq\cP_2(\R^d)$, $\mu\in\sQ$, and
$\cG:\sQ\to\R$. There is at most one
$\xi\in\Tan_\mu\sQ$ such that
\[
    \cG(\nu)-\cG(\mu)
    =
    \langle \xi,v_\pi\rangle_{L^2(\mu)}
    +o\left(W_2(\mu,\nu)\right)
\]
as $\nu\to\mu$ within $\sQ$, uniformly over
$\pi\in\Pi_{\mathrm{opt}}^\sQ(\mu,\nu)$.
\end{proposition}

\begin{proof}
Suppose $\xi_1,\xi_2\in\Tan_\mu\sQ$ both satisfy the
stated expansion, and write $\delta \defeq \xi_1-\xi_2$.
We show that $\delta$ is orthogonal to every generator
of $\Tan_\mu\sQ$.

Fix $\nu\in\sQ\setminus\{\mu\}$ and
$\pi\in\Pi_{\mathrm{opt}}^\sQ(\mu,\nu)$. For $t\in(0,1]$,
define
\[
    \mu_t \defeq ((1-t)x+ty)_\sharp\pi,
    \qquad
    \pi_t \defeq (x,(1-t)x+ty)_\sharp\pi.
\]
By the displacement-interpolation theorem
\cite[Theorem~7.2.2]{ambrosio2008gradient}, we have
\[
    \pi_t\in\Pi_{\mathrm{opt}}(\mu,\mu_t),
    \qquad
    W_2(\mu,\mu_t)=tW_2(\mu,\nu).
\]
Moreover, the interpolation associated with $\pi_t$ is
$s\mapsto\mu_{st}$, which lies in $\sQ$ since $\pi\in\Pi_{\mathrm{opt}}^\sQ(\mu,\nu)$. Thus
$\pi_t\in\Pi_{\mathrm{opt}}^\sQ(\mu,\mu_t)$.
Disintegrating with respect to the first coordinate gives
\[
    v_{\pi_t}(x)
    =
    \int\bigl((1-t)x+ty-x\bigr)\,\dd\pi(y\mid x)
    =
    t\,v_\pi(x)
    \qquad \mu\text{-a.e.}
\]
Applying the two differentiability expansions to
$(\mu_t,\pi_t)$ and subtracting yields
\[
    t\langle\delta,v_\pi\rangle_{L^2(\mu)}
    =
    o\left(t \,W_2(\mu,\nu)\right)
    =
    o(t).
\]
Dividing by $t$ and letting $t\downarrow0$, we obtain
$\langle\delta,v_\pi\rangle_{L^2(\mu)}=0$.
For $\nu=\mu$, the same conclusion holds because every
optimal coupling is concentrated on the diagonal, in which case
$v_\pi \equiv 0$.

By positive scaling and continuity of the inner product,
$\delta$ is therefore orthogonal to every element of
$\Tan_\mu\sQ$. Since $\xi_1,\xi_2\in\Tan_\mu\sQ$, we obtain
\[
    \|\delta\|_{L^2(\mu)}^2
    =
    \langle\delta,\xi_1\rangle_{L^2(\mu)}
    -\langle\delta,\xi_2\rangle_{L^2(\mu)}
    =0,
\]
and so $\xi_1=\xi_2$ in $L^2(\mu;\R^d)$.
\end{proof}

\subsection{Proof of Proposition~\ref{prop:interaction-gradient}}
\label{prf:interaction-gradient}

Set $M\defeq\sup_z\|D^2V(z)\|_{\mathrm{op}} < \infty$ and define
\[
    \psi_\mu(x)
    \defeq
    \sum_{j=1}^n
    \int
    V(x_1,\ldots,x_{j-1},x,x_{j+1},\ldots,x_n)
    \prod_{i\ne j}\dd\mu(x_i).
\]
Differentiation under the integral gives
$w_\mu=\nabla\psi_\mu\in L^2(\mu;\R^d)$ and
$\|D^2\psi_\mu\|_\infty\le nM$.
For sufficiently small $t>0$, the map
$\Id+t w_\mu$ is therefore the gradient of a convex
function and induces an optimal coupling. By the
definition of the tangent cone,
$w_\mu\in\Tan_\mu\cP_2(\R^d)$.

Now fix $\pi\in\Pi_{\mathrm{opt}}(\mu,\nu)$ and draw
independent pairs $(X_j,Y_j)\sim\pi$, $j=1,\ldots,n$.
A Taylor expansion then gives
\[
\begin{aligned}
&\left|
\cG(\nu)-\cG(\mu)
-\sum_{j=1}^n
\E\bigl[
\langle\nabla_jV(X_1,\ldots,X_n),Y_j-X_j\rangle
\bigr]
\right| \le
\frac{M}{2}\sum_{j=1}^n\E\|Y_j-X_j\|^2
=
\frac{nM}{2}W_2(\mu,\nu)^2.
\end{aligned}
\]
Independence identifies the linear term with the inner product
$\langle w_\mu,v_\pi\rangle_{L^2(\mu)}$.
The remainder is uniform over $\pi$, proving the claim.
\qed

\subsection{Proof of Proposition~\ref{prop:Q-grad-projection}}
\label{prf:Q-grad-projection}

Let $H\defeq\Tan_\mu\sQ$, $g\defeq\nabla_{W_2}\cG(\mu)$,
and $\xi\defeq\proj_H g$. Since $H$ is a closed linear subspace
of $L^2(\mu;\R^d)$, this projection is well-defined. For every
$\nu\in\sQ$ and $\pi\in\Pi_{\mathrm{opt}}^\sQ(\mu,\nu)$,
we have $v_\pi\in H$, so
\[
    \langle \xi,v_\pi\rangle_{L^2(\mu)}
    =
    \langle g,v_\pi\rangle_{L^2(\mu)}.
\]
Consequently, $W_2$-differentiability of $\cG$ gives
\[
    \cG(\nu)-\cG(\mu)
    =
    \langle \xi,v_\pi\rangle_{L^2(\mu)}
    +o\left(W_2(\mu,\nu)\right).
\]
The remainder remains uniform over the admissible couplings,
since
$\Pi_{\mathrm{opt}}^\sQ(\mu,\nu)
\subseteq\Pi_{\mathrm{opt}}(\mu,\nu)$.
Thus $\cG|_\sQ$ is $\sQ$-differentiable at $\mu$, with the
asserted gradient.\qed

\subsection{Proof of Proposition~\ref{prop:primary-projections}}
\label{prf:primary-projections}

We treat the three cases separately.\medskip

\textbf{Gaussian measures.}
Let $X \sim \mu$ and $Z \defeq X - m_\mu$, so that $a=\mathbb E[\xi(X)]$ and $C=\mathbb E[\xi(X)Z^\top]$. By Proposition~\ref{prop:main-tangent-spaces}, the orthogonal projection $p$ of $\xi$ onto $\Tan_\mu\BW$ has the form
\[
    p(x)=b+A(x-m_\mu)
\]
for some $b\in\R^d$ and $A\in\Sym$.
Orthogonality of $\xi-p$ to constant fields forces $b=a$. For $B \in \Sym$, orthogonality of $\xi-p$ to the affine field $x\mapsto B(x-m_\mu)$ forces
\[
    0
    =
    \E\bigl[\langle \xi(X)-a-AZ,BZ\rangle\bigr] = \E\bigl[\langle \xi(X)-AZ,BZ\rangle\bigr]
    =
    \langle C-A\Sigma_\mu,B\rangle_{\mathrm F},
\]
where $\langle\cdot,\cdot\rangle_{\mathrm F}$ denotes
the Frobenius inner product. This holding for all $B \in \Sym$ is equivalent to
\[
    A\Sigma_\mu+\Sigma_\mu A=C+C^\top.
\]
Conversely, any symmetric solution $A$ of this equation
makes the field $x \mapsto \xi(x)-a-A(x-m_\mu)$ orthogonal to every field in
$\Tan_\mu\BW$, and therefore yields the orthogonal projection.
Existence of a symmetric solution and uniqueness of its induced
$L^2(\mu)$ vector field thus follow from existence and uniqueness
of the orthogonal projection.\medskip

\textbf{Product measures.}
Let $X \sim \mu$. For each coordinate $i$, choose a Borel function $g_i$ such that
\[
    g_i(X_i)=\mathbb E[\xi_i(X)\,|\, X_i].
\]
Then, by the $L^2$ projection interpretation of the conditional expectation, we have that $g_i\in L^2(\mu_i)$ and, for every
$h_i\in L^2(\mu_i)$,
\[
    \E\bigl[
        (\xi_i(X)-g_i(X_i))h_i(X_i)
    \bigr]=0.
\]
Thus, by Proposition~\ref{prop:main-tangent-spaces}, $g(x)=(g_1(x_1),\ldots,g_d(x_d))$ belongs to
$\Tan_\mu\MF$, and summing these identities over $i$
shows that $\xi-g$ is orthogonal to $\Tan_\mu\MF$.
This proves the claimed formula.\medskip

\textbf{Radially symmetric measures.}
Fix a Borel representative of $\xi$, and write $\dd U$
for normalized Haar measure on $\Orth$. Define
\[
    E\defeq
    \left\{
        x:\int_{\Orth}\|\xi(U^\top x)\|^2\dd U<\infty
    \right\},
    \qquad
    p(x)\defeq
    \begin{cases}
        \displaystyle
        \int_{\Orth} U \xi(U^\top x)\dd U, & x\in E,\\
        0, & x\notin E.
    \end{cases}
\]
Invariance of $\mu$ and Tonelli's theorem give $\mu(E)=1$.
The set $E$ is orthogonally invariant, and Haar invariance
gives $p(Vx)=Vp(x)$ for every $V\in\Orth$ and $x\in\R^d$.
Moreover, Jensen's inequality gives
$\|p\|_{L^2(\mu)}\le\|\xi\|_{L^2(\mu)}$.
We thus have $p\in\Tan_\mu\Rad$.

For every $v\in\Tan_\mu\Rad$, invariance of $\mu$ and
equivariance of a representative of $v$ give
\[
    \langle p,v\rangle_{L^2(\mu)}
    =
    \int_{\Orth}
    \langle\xi,U v(U^\top\,\cdot)\rangle_{L^2(\mu)}
    \dd U
    =
    \langle\xi,v\rangle_{L^2(\mu)}.
\]
Hence, $p$ is precisely the orthogonal projection of $\xi$ onto
$\Tan_\mu\Rad$.
\qed

\subsection{Proof of Corollary~\ref{cor:interaction-energy-Q-grads}}
\label{prf:interaction-energy-Q-grads}
By Propositions~\ref{prop:interaction-gradient}
and~\ref{prop:Q-grad-projection}, $\nabla_\sQ\cG(\mu) = \proj_{\Tan_\mu\sQ}w_\mu$, for $w_\mu$ defined in Proposition~\ref{prop:interaction-gradient}. The mean-field formula therefore follows immediately from
Proposition~\ref{prop:primary-projections}\emph{(ii)}. We now treat the remaining two cases.\medskip

\textbf{Gaussian measures.} Boundedness of the Hessian of $V$
permits differentiation under the integrals defining $w_\mu$,
giving
\[
    Dw_\mu(x)
    =
    \sum_{j=1}^n
    \int
    D^2_{jj}V
    (x_1,\ldots,x_{j-1},x,x_{j+1},\ldots,x_n)
    \prod_{k\ne j}\dd\mu(x_k).
\]
Thus $Dw_\mu$ is bounded and symmetric. Setting $H \defeq \mathbb E_\mu[Dw_\mu]$, Gaussian integration by parts gives
\[
    C\defeq
    \mathbb E_\mu\bigl[w_\mu(X)(X-m_\mu)^\top\bigr]
    =
    H\Sigma_\mu.
\]
This identity also holds when $\Sigma_\mu$ is singular;
writing $X=m_\mu+\Sigma_\mu^{1/2}Z$, with
$Z\sim\cN(0,I_d)$, we can apply Gaussian integration by
parts in $Z$.

Since $H$ is symmetric, it solves $H\Sigma_\mu+\Sigma_\mu H=C+C^\top$. Thus, 
Proposition~\ref{prop:primary-projections}\emph{(i)} with $A=H$ yields
\[
    \nabla_{\BW}\cG(\mu)(x)
    =
    \E_\mu[w_\mu] +H(x-m_\mu),
\]
as claimed.\medskip

\textbf{Radially symmetric measures.}
Proposition~\ref{prop:primary-projections}\emph{(iii)} gives
\[
    q(x)\defeq\nabla_{\Rad}\cG(\mu)(x)
    =
    \mathbb E_{U\sim\Unif(\Orth)}
    \bigl[U^\top w_\mu(Ux)\bigr].
\]
Since $q$ is equivariant, for $r>0$ and $\|e\|=1$
we have $q(re)=\alpha(r)e$, where
\[
\begin{aligned}
    \alpha(r)
    &=
    \langle q(re),e\rangle =
    \mathbb E_{U\sim\Unif(\Orth)}
    \bigl[\langle w_\mu(rUe),Ue\rangle\bigr] =
    \mathbb E_{u\sim\Unif(\cS^{d-1})}
    \bigl[\langle w_\mu(ru),u\rangle\bigr].
\end{aligned}
\]
The last equality uses that $Ue$ is uniform on
$\cS^{d-1}$ whenever $U\sim\Unif(\Orth)$.
Taking $r=\|x\|$ and $e=x/\|x\|$ proves the asserted
formula for $x\ne0$. At the origin, equivariance gives
$q(0)=-q(0)$, and hence $q(0)=0$.\qed

\section{Proofs for Section~\ref{sec:W-FISTA}}
\label{app:W-FISTA}

\subsection{Proof of Proposition~\ref{prop:weak-geodesic-convexity}}
\label{prf:weak-geodesic-convexity}

The claim is direct for $\sQ = \cP_2(\R^d)$, in which case $\Pi_{\mathrm{opt}}^{\sQ}(\mu,\nu) = \Pi_{\mathrm{opt}}(\mu,\nu) \neq \emptyset$.

For $\sQ\in\{\BW,\MF,\Rad\}$,
Lemma~\ref{lem:structured-optimal-couplings} gives
$\Pi_{\mathrm{opt}}^\sQ(\mu,\nu)\neq\emptyset$ for every
$\mu,\nu\in\sQ$. Thus each family is weakly geodesically convex.\qed

\subsection{Proof of Proposition~\ref{prop:gradient_connection}}
\label{prf:gradient_connection}

Let $\mu \coloneq T_\sharp\rho$, and, for every $D \in \cH$, define $\nu \coloneq (T+D)_\sharp\rho$ and $\gamma \coloneq (T, T+D)_\sharp\rho$. Since $\mathsf{Q} = \mathcal{H}_\sharp \rho$ by Assumption~\ref{assm:extras}\emph{(iii)}, we have $\mu, \nu \in \mathsf{Q}$ for all $D \in \cH$. From differentiability of $\cG$, it holds that
$$\begin{aligned}
g(T + D) - g(T) = \mathcal{G}(\nu) - \mathcal{G}(\mu) = \mathbb{E}_{(X,Y)\sim \pi}[\langle \nabla_\sQ \cG(\mu)(X), Y - X \rangle] + o(W_2(\mu, \nu)).
\end{aligned}$$
as $\nu \to \mu$ in $\mathsf{Q}$, uniformly in $\pi \in \Pi_\mathrm{opt}^\mathsf{Q}(\mu, \nu)$.
Note that $\gamma$ is a suboptimal transport plan between $\mu$ and $\nu$, so that $W_2(\mu, \nu) \leq \|D\|_{L^2(\rho)}$. Thus, Fréchet differentiability of $g$ follows if show that there exists $v \in \cH$ such that $$\mathbb{E}_{(X,Y)\sim \pi}[\langle \nabla_\sQ \cG(\mu)(X), Y - X \rangle] = \langle v, D \rangle_{\cH} + o\left(\|D\|_{L^2(\rho)}\right),$$
as $\|D\|_{L^2(\rho)} \to 0$,
in which case $\nabla_\cH g = v$.

Let $\varepsilon >0$. Since $\nabla_\sQ\cG(\mu) \in \mathrm{Tan}_\mu \sQ \subseteq \mathrm{Tan}_\mu\mathcal{P}_2(\R^d)$, let $\varphi \in C_c^\infty(\mathbb{R}^d)$ such that $\|\nabla_\sQ\cG(\mu) - \nabla \varphi\|_{L^2(\mu)} \leq \varepsilon$. We have that
$$\begin{aligned}
    &\ \int\langle\nabla_\sQ \mathcal{G}(\mu)(x), y - x \rangle\dd\pi(x,y)\\
    &= \ \int\langle \nabla\varphi(x), y-x\rangle \dd\pi(x,y) + \int \langle\nabla_\sQ\cG(\mu)(x) - \nabla \varphi(x), y-x\rangle \dd \pi(x,y)\\
    &=\ \int\varphi(y) \dd\nu(y) - \int\varphi(x) \dd\mu(x) + O_\varepsilon(W_2^2(\mu, \nu)) + O(\varepsilon W_2(\mu, \nu))\\
    &=\ \int\langle\nabla_\sQ \mathcal{G}(\mu)(x), y - x \rangle\dd\gamma(x,y) + O_\varepsilon\left(\|D\|_{L^2(\rho)}^2 + W_2^2(\mu, \nu)\right) + O\left(\varepsilon (\|D\|_{L^2(\rho)} + W_2(\mu, \nu))\right)\\
    &=\ \int\langle\nabla_\sQ \mathcal{G}(\mu)(x), y - x \rangle\dd\gamma(x,y) + O_\varepsilon\left(\|D\|_{L^2(\rho)}^2\right) + O\left(\varepsilon\|D\|_{L^2(\rho)}\right)\\
    &=\ \langle \nabla_\mathsf{Q}\mathcal{G}(\mu) \circ T, D \rangle_{L^2(\rho)} + O_\varepsilon\left(\|D\|_{L^2(\rho)}^2\right) + O\left(\varepsilon\|D\|_{L^2(\rho)}\right).
\end{aligned}$$
Here, the constants implicit in $O(\cdot)$ are independent of $\varepsilon$ and $D$, while the constants in $O_\varepsilon(\cdot)$ are independent of $D$ but depend on $\varepsilon$. The first equality follows by linearity of the integral. The second equality follows by
Taylor's theorem  $\langle \nabla\varphi(x), y - x\rangle = \varphi(y) - \varphi(x) + O_\varepsilon(\|y - x\|_2^2)$ and Cauchy-Schwarz inequality $\int \langle\nabla_\sQ\cG(\mu)(x) - \nabla \varphi(x), y-x\rangle \dd \pi(x,y) \leq \|\nabla_\mathsf Q \cG(\mu) - \nabla \varphi\|_{L^2(\mu)}W_2(\mu, \nu)$. The third equality follows by repeating the argument of the second inequality but for the transport plan $\gamma$. The forth equality follows since $W_2(\mu, \nu) \leq \|D\|_{L^2(\rho)}$, by suboptimality of $\gamma$. The fifth inequality follows from the definition of $\gamma$. Thus,
$$
\begin{aligned}
\limsup_{D \neq 0,\, \|D\|\to 0} \frac{\left| \mathbb{E}_{(X,Y)\sim \pi}[\langle \nabla_\sQ \cG(\mu)(X), Y - X \rangle] - \langle \nabla_\sQ\cG(\mu)\circ T, D \rangle_{L^2(\rho)}\right|}{\|D\|_{L^2(\rho)}} \leq C\varepsilon,
\end{aligned}
$$
for all $\epsilon > 0$, where $C \geq 0$ is a constant. Letting $\epsilon \to 0$ gives $$ \mathbb{E}_{(X,Y)\sim \pi}[\langle \nabla_\sQ \cG(\mu)(X), Y - X \rangle] = \langle \nabla_\sQ\cG(\mu)\circ T, D \rangle_{L^2(\rho)} + o\left(\|D\|_{L^2(\rho)}\right).$$
Finally, since $D \in \cH$, we have $$\langle \nabla_\sQ\mathcal{G}(\mu) \circ T, D \rangle_{L^2(\rho)} = \langle \mathrm{proj}_\cH[\nabla_\sQ\cG(\mu)\circ T], D\rangle_{L^2(\rho)}$$ and  $\mathrm{proj}_\cH[\nabla_\sQ\cG(\mu)\circ T] \in \cH$, which proves the claim.
\qed

\subsection{Proof of Proposition~\ref{prop:standard_choices}}
\label{prf:standard_choices}

Assumption~\ref{assm:lift}\emph{(i)} holds by Brenier's theorem for
$\sQ=\cP_2(\R^d)$ and by
Proposition~\ref{prop:pushforward-expressivity} for the three subfamilies,
noting that $\Id\in\cT$ in every case.

The spaces $\cH$ are linear by construction.
Closedness is immediate when $\sQ=\cP_2(\R^d)$, while for the
subfamilies it follows from
Proposition~\ref{prop:main-tangent-spaces}, since
$\cH=\Tan_\rho\sQ$. Moreover, $\cH_\sharp\rho=\sQ$ follows from Brenier's theorem for $\sQ=\cP_2(\RR^d)$ and from
Proposition~\ref{prop:pushforward-expressivity}
for the subfamilies.

We next show that $\cM_{\mathrm{cm}}(\rho)$ is closed and convex
in $L^2(\rho;\R^d)$. Convexity follows from the subgradient
inequality by taking convex combinations of the corresponding
convex potentials. For closedness, suppose
$T_n\in\cM_{\mathrm{cm}}(\rho)$ and $T_n\to T$ in $L^2(\rho;\R^d)$.
Since
\[
    W_2((T_n)_\sharp\rho,T_\sharp\rho)
    \leq \|T_n-T\|_{L^2(\rho)}\to0,
\]
optimality of each $T_n$ gives
\[
    \|T-\Id\|_{L^2(\rho)}
    =
    \lim_{n\to\infty}\|T_n-\Id\|_{L^2(\rho)}
    =
    \lim_{n\to\infty}W_2(\rho,(T_n)_\sharp\rho)
    =
    W_2(\rho,T_\sharp\rho).
\]
Thus $T$ is an optimal transport map from $\rho$ to
$T_\sharp\rho$, and Brenier's theorem implies
$T\in\cM_{\mathrm{cm}}(\rho)$.
Consequently, $\cT=\cH\cap\cM_{\mathrm{cm}}(\rho)$ is closed
and convex.

Lastly, Assumption~\ref{assm:extras}\emph{(i)} follows by
Proposition~\ref{prop:weak-geodesic-convexity}, and
Assumption~\ref{assm:extras}\emph{(ii)} holds by definition of~$\cT$.
\qed

\subsection{Proof of Proposition~\ref{prop:convexity_canonical_lift_on_W2}}
\label{prf:convexity_canonical_lift_on_W2}

Geodesic convexity
of $\cG$ implies convexity along acceleration-free
curves by \cite[Theorem~1.1]{parker2024some} or
\cite[Theorem~5.6]{Cavagnari_2026_Lagrangian}, under the respective assumptions.
For arbitrary $T_0,T_1\in L^2(\rho;\RR^d)$,
let $\mu_i=(T_i)_\sharp\rho$ for $i\in\{0,1\}$.
Then
\[
    \mu_t=((1-t)T_0+tT_1)_\sharp\rho, \quad t\in[0,1],
\]
is an acceleration-free curve induced by the
coupling $(T_0,T_1)_\sharp\rho$. Consequently,
\[
\begin{aligned}
g((1-t)T_0+tT_1)
&=\cG(\mu_t)\\
&\leq(1-t)\cG(\mu_0)+t\cG(\mu_1)\\
&=(1-t)g(T_0)+tg(T_1),
\end{aligned}
\]
which proves convexity of $g$.

It remains to show that $\cG$ is geodesically convex if $g$ is convex. Let $\mu_0, \mu_1 \in \mathcal{P}_2(\RR^d)$ and $\gamma^\star \in \Pi_\mathrm{opt}(\mu_0, \mu_1)$ and consider the Wasserstein geodesic $\mu_t \coloneq ((1-t)x + t y)_\sharp \gamma^\star$, for $t \in [0,1]$. Since $\rho$ is a atomless, that is, $\rho(\{x\}) = 0$ for $x \in \RR^d$, so is $\rho \otimes \delta_0$. Thus, by \cite[Corollary~1.29]{santambrogio2015optimal}, there exists a Borel measurable map $F:\RR^{2d}\to\RR^{2d}$
such that $F_\sharp(\rho\otimes\delta_0)=\gamma^\star$. Define $T_0,T_1$ by $(T_0(x),T_1(x))=F(x,0)$.
Then $(T_0,T_1)_\sharp\rho=\gamma^\star$, and
$T_0,T_1\in L^2(\rho;\RR^d)$ since
$\mu_0,\mu_1\in\cP_2(\RR^d)$. Thus, $\mu_t = ((1-t)T_0 + tT_1)_\sharp\rho$ and
\[
\begin{aligned}
\cG(\mu_t)
&=g((1-t)T_0+tT_1)\\
&\leq (1-t)g(T_0)+tg(T_1)\\
&=(1-t)\cG(\mu_0)+t\cG(\mu_1),
\end{aligned}
\]
where the inequality is due to convexity of $g$. Since $\gamma^\star \in \Pi_\mathrm{opt}(\mu_0, \mu_1)$ is arbitrary, the claim follows. \qed

\subsection{Proof of Proposition~\ref{prop:interaction_energy_lift}}

The result follows from the convexity properties of interaction energies established in \cite[Proposition~9.3.5]{ambrosio2008gradient}, by considering suitable quadratic modifications of the interaction potential $V$. We provide a direct proof below to make the constants explicit.

Define the canonical lift $\tilde g:L^2(\rho;\RR^d)\to\RR$
by $\tilde g(S)=\cG(S_\sharp\rho)$.
Applying Proposition~\ref{prop:gradient_connection}
with $\sQ=\cP_2(\RR^d)$ and the full Hilbert space
$L^2(\rho;\RR^d)$ yields
\[
\nabla_{L^2(\rho)}\tilde g(T)
=
\nabla_{W_2}\cG(T_\sharp\rho)\circ T = w_{T_\sharp \rho} \circ T.
\]
Since $g=\tilde g|_\cH$ and $\cH$ is a closed linear subspace, the gradient of $g$ is the orthogonal
projection of the ambient gradient $\nabla_{L^2(\rho)}\tilde g$ onto $\cH$, that is,
\[
\nabla_\cH g(T)
=
\mathrm{proj}_\cH[\nabla_{L^2(\rho)}\tilde g(T)]
=
\mathrm{proj}_\cH
\left[w_{T_\sharp \rho} \circ T\right].
\]
Note that $$g(T) = \int V(T(x_1),\dots,T(x_n))\dd\rho^{\otimes n}(x_1,\dots, x_n) = \EE_{\rho^{\otimes n}}[V(\mathbf{T})], $$
where we define $\mathbf{T} \coloneq (T(X_1),\dots,T(X_n))$ for all $T \in \cH$. Moreover,
\[
\begin{aligned}
\langle\nabla_\cH g(T),S-T\rangle_{L^2(\rho)} &= \langle \mathrm{proj}_\cH [w_{T_\sharp\rho} \circ T],S-T\rangle_{L^2(\rho)}\\
&= \langle w_{T_\sharp\rho} \circ T,S-T\rangle_{L^2(\rho)}\\
&=\EE_{\rho^{\otimes n}}\left[
\sum_{j=1}^n
\left\langle
\partial_jV(T(X_1),\ldots,T(X_n)),
(S-T)(X_j)
\right\rangle
\right]\\
&=\EE_{\rho^{\otimes n}}\left[
\langle \nabla V(\mathbf{T}),
\mathbf{S} - \mathbf{T}\rangle
\right]
\end{aligned}
\]
where the second equality follows since $S-T \in \cH$, and the third equality follows by the definition of $w_\mu$ and Fubini's theorem. Here, the assumptions of Fubini's theorem follow from the Cauchy--Schwarz inequality, using that $\nabla V$ is Lipschitz continuous and
$S,T\in L^2(\rho;\RR^d)$. Thus,
\begin{align*}
&g(S) - g(T) - \langle\nabla_\cH g(T),S-T\rangle_{L^2(\rho)} = \EE[V(\mathbf{S}) - V(\mathbf{T}) - \langle \nabla V(\mathbf{T}),
\mathbf{S} - \mathbf{T} \rangle],
\end{align*}
and by $\ell_V$-strong convexity and $L_V$-smoothness of $V$,
$$\frac{\ell_V}{2}\EE_{\rho^{\otimes n}}\|\mathbf{S} - \mathbf{T}\|_2^2 \leq g(S) - g(T) - \langle\nabla_\cH g(T),S-T\rangle_{L^2(\rho)} \leq \frac{L_V}{2}\EE_{\rho^{\otimes n}}\|\mathbf{S} - \mathbf{T}\|_2^2.$$
Finally,
$$\EE_{\rho^{\otimes n}}\|\mathbf{S} - \mathbf{T}\|_2^2 = \EE_{\rho^{\otimes n}}\left[\sum_{i=1}^n \|S(X_i) - T(X_i)\|_2^2\right] = \sum_{i=1}^n \|S - T\|_{L^2(\rho)}^2  = n \|S - T\|_{L^2(\rho)}^2,$$
which proves the claim. \qed

\subsection{Proof of Proposition~\ref{prop:nonsmooth_lift_regular}}

The proof builds on the arguments in
\cite[Proposition~2.15 and Remarks~2.16--2.17]{craig2016exponential}
and
\cite[Proposition~2.1 and Lemmas~B.4--B.5]{dumont2026learning},
extending them to restricted classes of transport maps
$\cT\subseteq\cH\subseteq L^2(\rho;\RR^d)$ and measures
$\sQ\subseteq\cP_2(\RR^d)$.

Let $r$ be defined as in Assumption~\ref{assm:lift}\emph{(iii)}.

\textbf{Properness:} By Assumption~\ref{assm:lift}\emph{(i)}, $\cT_\sharp\rho=\sQ$,
so $\cR$ is proper if and only if $r$ is proper.

By Assumption~\ref{assm:lift}\emph{(i)}, every $\mu\in\sQ$
admits a representation $\mu=T_\sharp\rho$ with
$T\in\cT$. Moreover, by
Assumption~\ref{assm:extras}\emph{(ii)} and Brenier's theorem,
this representation is unique, with
$T=T^{\rho\to\mu}$. We use this repeatedly below.

\textbf{Convexity:} We show that $r$ is convex if and only if $\cR$
is generalized geodesically convex with pivot $\rho$.
Indeed, suppose first that $\cR$ has this property.
Let $T,S\in\operatorname{dom}(r)$ and $t\in[0,1]$,
and set $\mu_t\coloneq((1-t)T+tS)_\sharp\rho$. Since $\mathrm{dom}(r) \subseteq \cT$, $T,S\in\cT$ are
the optimal transport maps from $\rho$ to
$T_\sharp\rho$ and $S_\sharp\rho$, respectively.
Thus $(\mu_t)_{t\in[0,1]}$ is the generalized
geodesic between these measures with pivot $\rho$. Consequently,
\[
\begin{aligned}
r((1-t)T+tS)
&=\cR(\mu_t)\\
&\leq (1-t)\cR(T_\sharp\rho)
      +t\cR(S_\sharp\rho)\\
&=(1-t)r(T)+tr(S).
\end{aligned}
\]
where the first equality follows since $(1-t)T+tS\in\cT$ by convexity of $\cT$, the inequality is due to generalized geodesic convexity of $\cR$, and the second equality is due to $T,S \in \cT$. This proves convexity of $r$ on $\cH$. Conversely, suppose that $r$ is convex, and let
$\mu,\nu\in\operatorname{dom}(\cR)$.
Since $\mu, \nu \in \sQ$, by Assumption~\ref{assm:lift}\emph{(i)} there exists $T,S \in \cT$ such that $T_\sharp\rho = \mu$ and $S_\sharp\rho = \nu$, and thus,  $T=T^{\rho\to\mu}$ and $S=T^{\rho\to\nu}$. For the generalized geodesic
$\mu_t=((1-t)T+tS)_\sharp\rho$, convexity of $r$
therefore gives
\[
\begin{aligned}
\cR(\mu_t)
&=r((1-t)T+tS)\\
&\leq (1-t)r(T)+tr(S)\\
&=(1-t)\cR(\mu)+t\cR(\nu),
\end{aligned}
\]
where the first equality follows since $(1-t)T+tS\in\cT$ by convexity of $\cT$, the inequality inequality is due to convexity of $r$, and the second equality is due to $T,S \in \cT$. Thus, $\cR$ is generalized geodesically
convex with pivot $\rho$.

\textbf{Lower semicontinuity:} We show that $r$ is lower semicontinuous w.r.t. the $L^2(\rho)$-topology if and only if $\cR$ is lower semicontinuos w.r.t. the $W_2$-topology. We first show that the map
\[
\varphi:\cT\to\sQ,\quad  \varphi(T) \coloneq T_\sharp\rho,
\]
is a homeomorphism, that is, a continuous bijection with continuous inverse. Indeed, $\varphi$ is a bijection with inverse
$\varphi^{-1}(\mu)=T^{\rho\to\mu}$, continuity of $\varphi$ follows from
\begin{equation}
\label{eq:W_ineq}
W_2(T_\sharp\rho,S_\sharp\rho)
\leq \|T-S\|_{L^2(\rho)},
\end{equation}
while continuity of $\varphi^{-1}$ follows from
stability of optimal transport maps
\cite[Lemma~2.4.6]{panaretos2020invitation}.
It follows that
$\cR$ lower semicontinuous on $\sQ$ if
and only if $r|_\cT=\cR\circ\varphi$ is lower
semicontinuous on $\cT$, by applying \citep[Lemma 1.28]{bauschke2017convex} to the compositions $\cR \circ \varphi$ and $r|_\cT \circ \varphi^{-1}$, respectively.

Since $\cT$ is closed in $\cH$, we have that $r|_\mathcal{T}$ is lower semicontinuous on $\cT$ if and only if $r$ is lower semicontinuous on $\cH$. %
Indeed, since $r$ is the
extension of $r|_\cT$ by $+\infty$, the sublevel
sets of $r$ are closed in $\cH$ if and only if
those of $r|_\cT$ are closed in $\cT$.
Hence, by the characterization of lower
semicontinuity through closed sublevel sets
\cite[Lemma~1.24]{bauschke2017convex},
$r$ is lower semicontinuous on $\cH$ if and
only if $r|_\cT$ is lower semicontinuous on~$\cT$.

Finally, negative entropy $-\mathrm{Ent}$ is lower semicontinuous w.r.t. the $W_2$-topology
and convex along generalized geodesics with an
absolutely continuous pivot
\cite[Remark~9.3.8 and Proposition~9.3.9]
{ambrosio2008gradient}. For the restriction $-\mathrm{Ent}|_\sQ$, lower semicontinuity follows directly, while generalized geodesic convexity follows since the generalized geodesic between $\mu, \nu \in \sQ$ remains in $\sQ$ by convexity of $\cT$. Moreover, properness of $-\mathrm{Ent}|_\sQ$ is the extra assumption in the statement of Proposition~\ref{prop:nonsmooth_lift_regular}. \qed

\subsection{Proof of Proposition~\ref{prop:stochastic_gradient_connection}}

We have that $\cO_\cH$ is an unbiased estimator for $\nabla_\cH g(T)$, since \begin{align*}
        \EE_{\xi \sim \PP}[\cO_\cH(T, \xi)] &= \EE_{\xi \sim \PP}[\mathrm{proj}_\cH\left[\cO_\mathsf{Q}(T_\sharp \rho, \eta(T,\xi)) \circ T\right]] \\
        &= \mathrm{proj}_\cH\left[\EE_{\xi \sim \PP}[\cO_\mathsf{Q}(T_\sharp \rho, \eta(T,\xi))]\circ T\right] \\
        &= \mathrm{proj}_\cH\left[\EE_{\eta \sim \PP_{T_\sharp\rho}}[\cO_\mathsf{Q}(T_\sharp \rho, \eta)]\circ T\right] \\
        &= \mathrm{proj}_\cH\left[\nabla_\mathsf{Q}\cG(T_\sharp \rho)\circ T\right]\\
        &=\nabla_\cH g(T).
    \end{align*}
    The first equality follows by definition, the second equality since $\mathrm{proj}_\cH: L^2(\rho;\RR^d) \to \cH$ and $L^2(T_\sharp\rho;\RR^d)\ni v \mapsto v \circ T \in L^2(\rho;\RR^d)$ are bounded linear operators between Hilbert spaces and thus commute with the expectation, the third equality from the assumption on $\eta$, the fourth inequality from unbiasedness of $\cO_\mathsf{Q}$ in Definition~\ref{def:oracle-Q}, and the fifth equality from Proposition~\ref{prop:gradient_connection}.
    
For the variance bound, we have
    \begin{align*}
        &\E_{\xi\sim\PP}
\left[\left\|\cO_\cH(T,\xi)-\nabla_{\cH}g(T)\right\|_{L^2(\rho)}^2\right]\\ &= \E_{\xi\sim\PP}
\left[\left\| \mathrm{proj}_\cH\left[\cO_\mathsf{Q}(T_\sharp \rho, \eta(T,\xi)) \circ T\right] -\mathrm{proj}_\cH\left[\nabla_\mathsf{Q}\cG(T_\sharp \rho)\circ T\right]\right\|_{L^2(\rho)}^2\right]\\
&= \E_{\xi\sim\PP}
\left[\left\| \mathrm{proj}_\cH\left[\cO_\mathsf{Q}(T_\sharp \rho, \eta(T,\xi)) \circ T -\nabla_\mathsf{Q}\cG(T_\sharp \rho)\circ T\right]\right\|_{L^2(\rho)}^2\right]\\
&\leq \E_{\xi\sim\PP}
\left[\left\| \cO_\mathsf{Q}(T_\sharp \rho, \eta(T,\xi)) \circ T -\nabla_\mathsf{Q}\cG(T_\sharp \rho)\circ T\right\|_{L^2(\rho)}^2\right]\\
&=\E_{\xi\sim\PP}
\left[\left\| \cO_\mathsf{Q}(T_\sharp \rho, \eta(T,\xi)) -\nabla_\mathsf{Q}\cG(T_\sharp \rho)\right\|_{L^2(T_\sharp\rho)}^2\right]\\
&=\E_{\eta\sim\PP_{T_\sharp \rho}}
\left[\left\| \cO_\mathsf{Q}(T_\sharp \rho, \eta) -\nabla_\mathsf{Q}\cG(T_\sharp \rho)\right\|_{L^2(T_\sharp\rho)}^2\right]\\
&\leq \sigma^2\zeta_\mathsf{Q}(T_\sharp \rho)^2 = \sigma^2\zeta_\cH(T)^2.
    \end{align*}
The first equality follows from the definition of $\cO_\cH$ and Proposition~\ref{prop:gradient_connection}, the second equality from linearity of the projection operator, and the first inequality since projections are non-expansive. The third and fourth equalities follow from the
pushforward identities for $T$ and $\eta(T,\cdot)$,
respectively. The second inequality follows since $\cO_\mathsf{Q}$ satisfies the variance bound in Definition~\ref{def:oracle-Q}, and the last equality from the definition of $\zeta_\cH$.

Finally, for any $T,S\in\cH$, we have
\begin{align*}
    |\zeta_\cH(T)-\zeta_\cH(S)| = |\zeta_\sQ(T_\sharp\rho)-\zeta_\sQ(S_\sharp\rho)| \leq W_2(T_\sharp\rho,S_\sharp\rho) \leq\|T-S\|_{L^2(\rho)},
\end{align*}
where the equality is the definition, the first inequality follows from the
$1$-Lipschitz continuity of $\zeta_\sQ$, and the second inequality from \eqref{eq:W_ineq}.
Thus, $\zeta_\cH$ is $1$-Lipschitz, completing
the proof.
\qed

\subsection{Proof of Proposition~\ref{prop:stochastic-Q-examples}}

Since $V\in C^2((\RR^d)^n)$ is $L_V$-smooth,
it is readily verified that $w_\mu$ is continuously
differentiable and satisfies
\[
    \|w_\mu(x)-w_\mu(y)\|_2
    \leq nL_V\|x-y\|_2,
    \quad
    \|Dw_\mu(x)\|_{\mathrm{op}}\leq nL_V.
\]
We first establish unbiasedness and the variance bounds
for each family.

\textbf{Gaussian case ($\sQ=\mathsf{BW}$).} %
Let $\mu = \cN(m_\mu, \Sigma_\mu) \in \BW$. Since the Jacobian $Dw_\mu$ is symmetric, $\cO_\mathsf{BW}(\mu, \eta) \in \mathrm{Tan}_\mu \mathsf{BW}$. 

Unbiasedness follows, since
\[
\E_{\eta\sim \mu}[\cO_\mathsf{BW}(\mu,\eta)(x)]
=
\E_\mu[w_\mu(X)]
+
\E_\mu[Dw_\mu(X)](x-m_
\mu) = \nabla_\mathsf{BW}\cG(\mu)(x),
\]
where the second equality follows 
by Corollary~\ref{cor:interaction-energy-Q-grads}.

We show the variance bound. We have
\begin{align*}
&\E_{\eta \sim \mu}
\|\cO_\mathsf{BW}(\mu,\eta)
-\nabla_\mathsf{BW}\cG(\mu)\|_{L^2(\mu)}^2\\
&=\E_\mu\|w_\mu(X)-\E_\mu[w_\mu(X')]\|_2^2 + 
\E_\mu\|
(Dw_\mu(X)-\E_\mu[D w_\mu(X')])\Sigma_\mu^{1/2}
\|_{\mathrm F}^2,
\end{align*}
where we used that the constant and centered linear components are
orthogonal in $L^2(\mu)$. For the first term, we have
\[
\E_\mu\|w_\mu(X)-\E_\mu[w_\mu(X')]\|_2^2
\leq
\E_\mu\|w_\mu(X)-w_\mu(m_\mu)\|_2^2
\leq n^2 L_V^2\zeta_\mathsf{BW}(\mu)^2.
\]
where the first inequality follows since the expectation $\E_\mu[w_\mu(X')]$ minimizes the function $z \mapsto \EE_\mu\|w_\mu(X) - z\|_2^2$ on $\RR^d$ and the second inequality from $nL_V$-Lipschitz continuity of $w_\mu$. Similarly, for the second term,
\begin{align*}
&\E_\mu\|
(Dw_\mu(X)-\E_\mu[D w_\mu(X)])\Sigma_\mu^{1/2}
\|_{\mathrm F}^2\leq
\E_\mu\|Dw_\mu(X)\Sigma_\mu^{1/2}\|_{\mathrm F}^2\leq
n^2L_V^2\operatorname{tr}(\Sigma_\mu)
=
n^2 L_V^2\zeta_\mathsf{BW}(\mu)^2,
\end{align*}
where the first inequality follows from the variance identity, the second inequality from the operator norm bound on $Dw_\mu$, and the equality from $\zeta_\mathsf{BW}(\mu)^2 = \EE_\mu \| X  - m_\mu\|_2^2 = \mathrm{tr}(\Sigma_\mu)$. Combining the two bounds yields
\[
\E_{\eta\sim \mu}\bigl[
\|\cO_\mathsf{BW}(\mu,\eta)
-\nabla_\mathsf{BW}\cG(\mu)\|_{L^2(\mu)}^2
\bigr]
\leq 2n^2L_V^2\zeta_\mathsf{BW}(\mu)^2.
\]

\textbf{Mean-field case ($\sQ=\mathsf{MF}$).} Let $\mu = \otimes_{i=1}^d \mu_i \in \MF$. The oracle is coordinate-wise separable, that is, $\cO_\mathsf{MF}(\mu, \eta) \in \mathrm{Tan}_\mu \MF$.

Letting $\mu_{-i}\coloneq\bigotimes_{j\neq i}\mu_j$ denote the product of all marginals except the $i$-th, unbiasedness follows since
\[
\E_{\eta\sim \mu}[
\cO_\mathsf{MF}(\mu,\eta)(x)]_i
=
\E_{X_{-i}\sim\mu_{-i}}
[w_{\mu,i}(x_i,X_{-i})] = [\nabla_\mathsf{MF}\cG(\mu)(x)]_i.
\]
where the last equality follows by Corollary~\ref{cor:interaction-energy-Q-grads}.

For the variance bound, we have
\begin{align*}
    &\E_{\eta\sim\mu}
    \left[
        \|\cO_\MF(\mu,\eta)
        -\nabla_\MF\cG(\mu)\|_{L^2(\mu)}^2
    \right]\\
    &=
    \sum_{i=1}^d
    \E_{X_i\sim\mu_i}
    \left[
        \E_{X_{-i}\sim\mu_{-i}}
        \left[
            \left|
                w_{\mu,i}(X_i,X_{-i})
                -
                \E_{X_{-i}'\sim\mu_{-i}}
                [w_{\mu,i}(X_i,X_{-i}')]
            \right|^2
        \right]
    \right]\\
    &=
    \sum_{i=1}^d
    \E_\mu
    \left[
        \operatorname{Var}_\mu
        [w_{\mu,i}(X)\mid X_i]
    \right]\\
    &\leq
    \sum_{i=1}^d
    \operatorname{Var}_\mu[w_{\mu,i}(X)]\\
    &=
    \E_\mu\|w_\mu(X)-\E_\mu[w_\mu(X')]\|_2^2\\
    &\leq
    n^2L_V^2\zeta_\MF(\mu)^2,
\end{align*}
where the first equality follows from the product structure of $\mu$ and Tonelli's theorem, the second equality from the definition of conditional variance, the first inequality from the law of total variance, the third equality since $\mathrm{Var}_\mu [w_{\mu,i}(X)] = \EE_\mu[|w_{\mu, i}(X) - \EE_\mu[w_{\mu, i}(X')]|^2]$, and the second inequality from the same arguments as in the Gaussian case.

\textbf{Radial case ($\sQ=\mathsf{Rad}$).}
Let $\PP \coloneq \mathrm{Unif}(\cS^{d-1})$. Clearly, $\cO_\Rad(\mu, \eta) \in \mathrm{Tan}_\mu\Rad$.

We show unbiasedness. If $x = 0$, then $\EE_{\eta \sim \PP}[\cO_\mathsf{Rad}(\mu, \eta)(x)]=0=\nabla_\mathsf{Rad}\cG(\mu)(x)$ by definition. If $x \neq 0$, then
\begin{align*}
\E_{\eta \sim \PP}[\cO_\mathsf{Rad}(\mu,\eta)(x)]=
\E_{\eta \sim \PP}\left[
\langle w_\mu(\|x\|_2\eta),\eta \rangle
\right]\frac{x}{\|x\|_2} = \nabla_\mathsf{Rad}\cG(\mu)(x),
\end{align*}
where the first equality uses $\EE_\PP[\eta] = 0$, and the second equality follows by Corollary~\ref{cor:interaction-energy-Q-grads}.

We show the variance bound. For all $x \in \RR^d$,
\begin{equation}
\label{eq:rad_oracle_norm_bound}
    \|\cO_\mathsf{Rad}(\mu,\eta)(x)\|_2 \leq nL_V\|x\|_2.
\end{equation}
Indeed, both sides are $0$ if $x = 0$, and for every $x\neq0$,
\begin{align*}
\|\cO_\mathsf{Rad}(\mu,\eta)(x)\|_2
&\leq
\|w_\mu(\|x\|_2\eta)-w_\mu(0)\|_2\leq  nL_V\|x\|_2,
\end{align*}
where the first inequality follows from Cauchy-Schwarz and $\|\eta\|_2 = 1$, and the second inequality from $nL_V$-Lipschitz continuity of $w_\mu$.
The variance bound now follows since
\begin{align*}
\E_{\eta \sim \PP}\bigl[
\|\cO_\mathsf{Rad}(\mu,\eta)
-\nabla_\mathsf{Rad}\cG(\mu)\|_{L^2(\mu)}^2
\bigr]&\leq
\E_{\eta\sim \PP}\|\cO_\mathsf{Rad}(\mu,\eta)\|_{L^2(\mu)}^2 \\
&\leq n^2L_V^2\E_\mu\|X\|_2^2\\
&=n^2L_V^2\zeta_\mathsf{Rad}(\mu)^2,
\end{align*}
where the first inequality follows from the variance identity, the second inequality from \eqref{eq:rad_oracle_norm_bound}, and the equality from $m_\mu = 0$ for any $\mu \in \Rad$.

\textbf{Lipschitz continuity of $\zeta_\sQ$.} For every $\sQ \in \{\mathsf{BW}, \mathsf{MF}, \mathsf{Rad}\}$ and $\mu, \nu \in \sQ$, $\zeta_\mathsf{Q}$ is $1$-Lipschitz. Indeed,
\begin{equation}
\label{eq:zeta_Q_bounds}
\begin{aligned}
\zeta_\sQ(\mu)
&= \left(\EE_\mu\|X - m_\mu\|_2^2\right)^{1/2} \\
&= \inf_{a\in\RR^d}W_2(\mu,\delta_a) \\
&\leq W_2(\mu, \nu) + \inf_{a\in\RR^d}W_2(\nu,\delta_a)\\
&=W_2(\mu, \nu) + \zeta_\mathsf{Q}(\nu),
\end{aligned}
\end{equation}
where the first equality is the definition, the second equality follows from $W_2(\mu,\delta_a) = (\EE_\mu\|X - a\|_2^2)^{1/2}$ and the fact that the latter is minimized at $a = m_\mu$, the inequality follows from the triangle inequality, and the third equality was already established. Exchanging $\mu$ and $\nu$ in \eqref{eq:zeta_Q_bounds}, we obtain
\[
|\zeta_\sQ(\mu) - \zeta_\mathsf{Q}(\nu)| \leq W_2(\mu, \nu),
\]
which proves the claim.

\textbf{Measurability.} Finally, for every $\sQ \in \{\mathsf{BW}, \mathsf{MF}, \mathsf{Rad}\}$ and $\mu\in\sQ$, the map $\eta\mapsto\cO_\sQ(\mu,\eta)$ is continuous as a map into $L^2(\mu;\RR^d)$, and hence measurable.
Continuity follows from continuity of $w_\mu$ and $Dw_\mu$
for $\sQ = \BW$, Lipschitz continuity of $w_\mu$
for $\sQ =\MF$, and dominated convergence using
the pointwise norm bound \eqref{eq:rad_oracle_norm_bound} for $\sQ =\Rad$.

Thus, all three oracles satisfy
Definition~\ref{def:oracle-Q} with the stated
noise parameters.\qed

\subsection{Proof of Theorem~\ref{thm:w-fista_rate}}
\label{app:deterministic_wfista}

The argument below directly applies FISTA analysis from \cite{chambolle2016introduction}, combined with some small tweaks to account for $q=1$ and a slightly sharper rate when $\ell > 0$.

Fix any $T^\star\in\cT$ such that
$(T^\star)_\sharp\rho=\mu^\star$. By Assumption~\ref{assm:lift},
$T^\star$ minimizes $g+r$ over $\cH$, and
\[
\Delta_k\defeq\cF(\mu_k)-\cF(\mu^\star)
=(g+r)(T_k)-(g+r)(T^\star),
\qquad k\geq1.
\]
Since $r=+\infty$ outside $\cT$, the proximal step may equivalently be taken
over all of $\cH$.

Since $Y_0=T_0$, the proximal descent
inequality (4.37) of \citet{chambolle2016introduction} implies that
\begin{equation}
\label{eq:FISTA-det-first-step}
2h\Delta_1+\|T_1-T^\star\|_{L^2(\rho)}^2
\leq
(1-q)\|T_0-T^\star\|_{L^2(\rho)}^2.
\end{equation}
Indeed, while stated for finite-dimensional Euclidean space, their argument only uses its Hilbert space structure and applies unchanged in $\cH$.
If $q=1$, the right-hand side vanishes and hence $T_1=T^\star$.
Our Step~\ref{step:weight-recurrence} recurrence then gives $t_k=1$ and
$\beta_k=0$ for every $k\geq1$, so all subsequent iterates
equal $T^\star$. Thus, the rate and complexity claims trivially hold
in this case.

Suppose now that $q<1$.
The recurrence at Step~\ref{step:weight-recurrence} implies
\[
(1-qt_{k+1})(t_{k+1}+qt_k^2)
=(1-q)t_{k+1},
\]
and thus, since $q<1$,  that
\[
\beta_k
=\frac{t_k-1}{t_{k+1}+qt_k^2}
=\frac{t_k-1}{t_{k+1}}
  \frac{1-qt_{k+1}}{1-q}.
\]
Hence, Algorithm~\ref{alg:W-FISTA} is exactly the FISTA iteration in
\citet[Algorithm~5]{chambolle2016introduction}, with strong
convexity parameter $\ell$ in the smooth term and zero in
the nonsmooth term. Starting from the first-step inequality
\eqref{eq:FISTA-det-first-step}, their energy recursion leading to
\citep[(B.7)]{chambolle2016introduction} yields
\[
\Delta_k
\leq
\gamma_k\frac{\|T_0-T^\star\|_{L^2(\rho)}^2}{2h},
\qquad
\gamma_k\defeq
\frac{1}{t_k^2}\prod_{j=1}^k(1-qt_j).
\]
Taking the infimum over all $T^\star$ such that $(T^\star)_\sharp \rho = \mu^\star$ gives $\Delta_k\leq R^2\gamma_k/(2h)$.

If $q=0$, the Step~\ref{step:weight-recurrence} recurrence gives
$t_{k+1}\geq t_k+1/2$. Since $t_1=1$, we have
$t_k\geq(k+1)/2$ and therefore
$\gamma_k=t_k^{-2}\leq 4/(k+1)^2$, which implies
\eqref{eq:chambolle_bound}.

If instead $0<q<1$, the Step~\ref{step:weight-recurrence} recurrence
implies that
$\gamma_k=(1-qt_k^2)/t_k^2=t_k^{-2}-q>0$ and that
$\gamma_{k+1}=(1-t_{k+1}^{-1})\gamma_k$.
Setting $s_k\defeq\gamma_k^{-1/2}$, the second identity gives
$s_{k+1}^2-s_k^2=s_{k+1}^2/t_{k+1}>0$,
and so $s_{k+1}\geq s_k$.
Factoring the difference of squares and using
$t_{k+1}^{-2}=q+s_{k+1}^{-2}$, we thus obtain
\[
\begin{aligned}
s_{k+1}-s_k
&=\frac{s_{k+1}^2}{t_{k+1}(s_{k+1}+s_k)} =\frac{s_{k+1}\sqrt{1+qs_{k+1}^2}}{s_{k+1}+s_k}
\geq
\frac12\sqrt{1+qs_{k+1}^2}.
\end{aligned}
\]
Moreover, we have
\[
s_1=\frac1{\sqrt{1-q}}
\geq\frac{\sinh(\sqrt q)}{\sqrt q},
\]
where the inequality follows from $\tanh(\sqrt q)\leq\sqrt q$.
The inductive argument in the proof of
\citet[Lemma~2.2.4, (2.2.11)--(2.2.12)]{nesterov2018lectures},
started at $k=1$, therefore gives
\[
s_k\geq
\frac1{\sqrt q}
\sinh\left(\frac{k+1}{2}\sqrt q\right).
\]
Taking $r_k\defeq(k+1)\sqrt q/2$, we obtain
\begin{equation}
\label{eq:WDfista_gamma}
\gamma_k
\leq\frac{q}{\sinh^2(r_k)}
\leq\frac{16}{(k+1)^2}\exp(-r_k),
\end{equation}
where the last inequality uses
$\sinh x=2\sinh(x/2)\cosh(x/2)\geq(x/2)e^{x/2}$.
Combining this with our earlier bound $\Delta_k\leq\gamma_kR^2/(2h)$ proves
\eqref{eq:chambolle_bound}. 

For the iteration complexity, take $h=1/L$.
Let $A$ denote the applicable branch expression in
\eqref{eq:iteration-complexity-det}, and choose
$N_{\mathrm{det}}\defeq\max\{1,\lceil5A\rceil\}$.
Then $N_{\mathrm{det}}\asymp1+A$.
Suppose $k\ge N_{\mathrm{det}}$. If $\eps>\ell R^2/4$, which we note includes the case when $\ell=0$, then
$k\geq5\sqrt{LR^2/\eps}$ and \eqref{eq:chambolle_bound} imply that
\[
\Delta_k
\leq\frac{8LR^2}{(k+1)^2}
\leq\frac{8}{25}\eps
\leq\eps.
\]
Otherwise, set $a\defeq\ell R^2/\eps\geq4$.
If $k\geq5\sqrt{L/\ell}\log a$, then
$k+1\geq\sqrt{L/\ell}$, so \eqref{eq:chambolle_bound} yields
\[
\begin{aligned}
\Delta_k
&\leq
\frac{8LR^2}{(k+1)^2}
\exp\left(-\frac{k+1}{2}\sqrt{\frac{\ell}{L}}\right) \leq
8\ell R^2a^{-5/2}
\leq
\ell R^2a^{-1}
=\eps,
\end{aligned}
\]
where the last inequality uses $a^{3/2}\geq8$.
This proves \eqref{eq:iteration-complexity-det}.

Finally, under Assumption~\ref{assm:extras}(ii), every
$T^\star\in\cT$ satisfying $(T^\star)_\sharp\rho=\mu^\star$
is an optimal transport map from $\rho$ to $\mu^\star$.
Since $T_0=\Id$, it follows that
\[
\|T^\star-T_0\|_{L^2(\rho)}^2
=
\int\|T^\star(x)-x\|_2^2\,\dd\rho(x)
=
W_2^2(\rho,\mu^\star).
\]
Taking the infimum over all such $T^\star$ gives
$R=W_2(\rho,\mu^\star)$.
\qed

\subsection{Proof of Proposition~
\ref{prop:Q_gradient_oracle_complexity}}
\label{prf:Q_gradient_oracle_complexity}

The upper bound follows directly from
Theorem~\ref{thm:w-fista_rate} and
Proposition~\ref{prop:gradient_connection}, since each $\cH$-gradient evaluation for $g$ can be computed using one $\sQ$-gradient evaluation for $\cG$.

For the lower bound, let $R>0$ and fix $(\sQ,\rho,\cH,\cT)$ satisfying Assumption~\ref{assm:lift}(i), Assumption~\ref{assm:extras}(i--ii), and $\cH_\sharp\rho=\sQ$, with $\sQ$ closed under translations (for example, any of the non-radial settings in Proposition~\ref{prop:standard_choices} suffice). We now construct objectives which satisfy the remaining assumptions.

For $a\in\R^d$, write $\mu_a=(\Id+a)_\sharp\rho$.
Translation closure gives $\mu_a\in\sQ$, and
$T^{\rho\to\mu_a}=\Id+a$. Hence $\Id+a\in\cT$ for every
$a$, and $\cH$ contains all constant maps. Define
\[
\mathsf C_{\rho,R}
\defeq\{\mu_a:\|a\|\le R\},
\qquad
\cR\defeq\iota_{\mathsf C_{\rho,R}}.
\]
By uniqueness of the optimal transport maps, the corresponding
lift is exactly $r=\iota_{\{\Id+a:\|a\|\le R\}}$, which is proper, convex, and lower semicontinuous. Thus Assumption~\ref{assm:lift}\emph{(iii)} holds.

Next, for any $\ell$-strongly convex, $L$-smooth function
$f:\R^d\to\R$, define
\[
\begin{aligned}
\cG_f(\mu)
&\defeq f(m_\mu-m_\rho)
+\frac{\ell}{2}\operatorname{tr}\Sigma_\mu,\\
g_f(T)
&\defeq f(\E_\rho T-m_\rho)
+\frac{\ell}{2}\|T-\E_\rho T\|_{L^2(\rho)}^2.
\end{aligned}
\]
Then $g_f(T)=\cG_f(T_\sharp\rho)$ for every $T\in\cH$.
Since $\cH$ contains the constant maps, its orthogonal
decomposition into constant and mean-zero maps shows that
$g_f$ is $\ell$-strongly convex and $L$-smooth.
This verifies the remaining conditions in
Assumptions~\ref{assm:lift} and~\ref{assm:extras}. Ultimately, we will show that each deterministic method in question will suffer suboptimality at least $\eps$ for some choice of $f$.

A Taylor expansion of the mean and second-moment terms gives
the ambient Wasserstein gradient
\[
\xi_\mu(x)
=\nabla f(m_\mu-m_\rho)+\ell(x-m_\mu).
\]
This field also belongs to $\Tan_\mu\sQ$.
Indeed, linearity of $\cH$ and $\cH_\sharp\rho=\sQ$ imply
that $\sQ$ is closed under dilations, and it is closed under
translations by hypothesis. Thus, for every $s>0$, the map
\[
\Id+s\xi_\mu
=(1+s\ell)\Id
+s\bigl(\nabla f(m_\mu-m_\rho)-\ell m_\mu\bigr)
\]
is cyclically monotone, and the displacement interpolation
from $\mu$ to $(\Id + s\xi_\mu)_\sharp \mu$ remains in $\sQ$.
Thus, $\cG_f$ is $\sQ$-differentiable with $\nabla_\sQ\cG_f(\mu) = \xi_\mu$.
In particular, a query to $\nabla_\sQ\cG_f$ at any
$\mu\in\sQ$ can be simulated using one Euclidean gradient
query at $m_\mu-m_\rho$.

Writing $\cF_f=\cG_f+\cR$, we have
\[
\cF_f(\mu_a)
=f(a)+\frac{\ell}{2}\operatorname{tr}\Sigma_\rho,
\qquad
W_2(\rho,\mu_a)=\|a\|.
\]
Thus, minimizing $\cF_f$ over $\sQ$ is equivalent to
minimizing $f$ over the closed ball of radius $R$.
Every output with finite objective gap is of the form $\mu_a$
and yields a Euclidean output $a$ with the same
gap. The reference measure, lift spaces, and
regularizer are all fixed independently of $f$.

We first show that at least one query is necessary for all $\eps>0$. A deterministic procedure making no queries has an
output independent of $f$. An output outside
$\mathsf C_{\rho,R}$ has infinite suboptimality; otherwise,
write it as $\mu_a$ and choose a unit vector $u$ with
$\langle u,a\rangle\le0$. Taking
\[
f(z)=\frac{\ell}{2}\|z\|^2-A\langle u,z\rangle,
\]
for $A>\ell R$ sufficiently large makes $Ru$ the
constrained minimizer and gives
\[
f(a)-f(Ru)
\ge AR-\frac{\ell R^2}{2}
>\eps.
\]
This proves the desired lower bound whenever $\eps>LR^2/16$.
Indeed, $\ell\le L/4$ ensures that the first accuracy regime
applies, and
\[
N_{\mathrm{det}}(\eps;\ell,L,R)
=
\left\lceil5\sqrt{\frac{LR^2}{\eps}}\right\rceil
\le20.
\]
Selecting $c\le1/20$ ensures that any procedure using fewer than $cN_{\mathrm{det}}$ queries
must therefore make no queries, and the construction above
gives a valid instance with suboptimality greater than $\eps$.

It remains to consider the case where $0<\eps\le LR^2/16$. Setting
\[
\lambda\defeq\max\{\ell,4\eps/R^2\}\in(0,L/4],
\]
we consider $f$ ranging over all $\lambda$-strongly convex, $L$-smooth functions
with an unconstrained minimizer in the ball of radius $R$.
They are admissible since $\lambda\ge\ell$, and their
constrained and unconstrained minima coincide.
By \citet[Corollary~3]{drori2022oracle}, every deterministic
procedure using $k$ queries has worst-case suboptimality
at least
\[
\lambda R^2
\left(1-\sqrt{\lambda/L}\right)^{2k}
\ge
\lambda R^2\exp\left(-4k\sqrt{\lambda/L}\right),
\]
provided $d\ge2k+1$; here we used $\lambda/L\le1/4$.
Inverting this bound gives a necessary query budget of order
\[
\sqrt{\frac L\lambda}
\log\frac{\lambda R^2}{\eps}
\asymp N_{\mathrm{det}}(\eps;\ell,L,R),
\]
as desired. Indeed, $\lambda=4\eps/R^2$ in the first accuracy regime
and $\lambda=\ell$ in the second.
Moreover, the dimension requirement holds for every
$k<N_{\mathrm{det}}$ by the assumed bound on $d$.
\qed

\subsection{Proof of Theorem~\ref{thm:w-fista_rate_stochastic}}
\label{app:stochastic_wfista}

Throughout, norms and inner products are those of $\cH$.
Write $f=g+r$ and $f^\star=\cF(\mu^\star)$.
The set $\argmin_\cH f$ is nonempty, closed, and convex.
Its closest point $T^\star$ to $T_0$ therefore exists and
satisfies
\[
f(T^\star)=f^\star,
\qquad
\|T_0-T^\star\|\le R.
\]
Define the gradient error and variance by
\[
e_k \defeq \cO_\cH(Y_k,\xi_k)-\nabla_\cH g(Y_k),
\qquad
v_k \defeq \E\|e_k\|^2.
\]
Let $\cA_k$ denote the $\sigma$-algebra generated by $\xi_0,\ldots,\xi_{k-1}$.
Joint measurability of the oracle and continuity of the
proximal map imply that $T_k$ and $Y_k$ are
$\cA_k$-measurable.
The linear growth of $\nabla g$ and $\zeta_\cH$, the
oracle variance bound, and nonexpansiveness of the
proximal map show inductively that the iterates and
gradient errors have finite second moments.
Independence and oracle unbiasedness therefore give
$\E[e_k\mid\cA_k]=0$, while the variance bound gives
\begin{equation}
\label{eq:WSfista_conditional_noise}
\begin{aligned}
v_k &\le \sigma^2\E[\zeta_\cH(Y_k)^2]
\le \sigma^2\E\bigl[
(\zeta_\cH(T_0)+R+\|Y_k-T^\star\|)^2
\bigr]\\
&\le
2\sigma^2\bigl[
(\zeta_\cH(T_0)+R)^2+\E\|Y_k-T^\star\|^2
\bigr].
\end{aligned}
\end{equation}

We now establish \eqref{eq:WSfista_rate} in two steps. First, we extend to stochastic gradients the argument
yielding $\Delta_k\le\gamma_kR^2/(2h)$ in the proof of
Theorem~\ref{thm:w-fista_rate}.
We will track an energy combining the weighted objective
gap with the squared norm of an auxiliary vector $u_k$,
defined below.
Conditional unbiasedness will eliminate a linear error
term after taking expectations, leaving a bound on the
expected energy in terms of its initial value and the
accumulated error variances.

Second, we use the energy estimate to control
$\E\|Y_k-T^\star\|^2$ in
\eqref{eq:WSfista_conditional_noise}.
This will relate each error variance to the preceding
variances. The step-size restrictions in the statement then allow us to absorb their accumulated sum, giving
$\max_{0\le j<k}v_j\le2\bar\sigma^2$.
Substituting this uniform bound into the energy estimate
will prove \eqref{eq:WSfista_rate}. We then choose the
step size to deduce the stated iteration complexity.\medskip

\textbf{Energy estimate.}
We adapt the energy argument of
\citet[Appendix~B, proof of Theorem~4.10]{chambolle2016introduction},
retaining the terms arising from gradient errors. A similar stochastic energy estimate also appears in
\citet[Theorem~12]{kulunchakov2020estimate}.
We give the calculation in our notation and then combine it
with the oracle assumption to bound variances uniformly
over the first $k$ steps.

To start, proximal optimality and the smoothness and strong convexity
of $g$ give, for every $T\in\operatorname{dom}r$,
\[
\begin{aligned}
f(T_{k+1})
\le{}&
f(T)
+\frac{(1-q)\|T-Y_k\|^2-\|T-T_{k+1}\|^2}{2h}\\
&+\left(\frac L2-\frac1{2h}\right)
  \|T_{k+1}-Y_k\|^2
+\langle e_k,T-T_{k+1}\rangle.
\end{aligned}
\]
We split the gradient error term and apply Young's inequality to obtain
\[
\begin{aligned}
\langle e_k,T-T_{k+1}\rangle
&=
\langle e_k,T-Y_k\rangle
+\langle e_k,Y_k-T_{k+1}\rangle\\
&\le
\langle e_k,T-Y_k\rangle
+h\|e_k\|^2
+\frac1{4h}\|T_{k+1}-Y_k\|^2.
\end{aligned}
\]
Combining these inequalities, the coefficient of
$\|T_{k+1}-Y_k\|^2$ becomes $L/2 - 1/(4h) \leq 0$, since $h \leq 1/(2L)$.
Discarding this nonpositive term yields the descent inequality
\begin{equation}
\label{eq:WSfista_prox}
f(T_{k+1})
\le
f(T)
+
\frac{(1-q)\|T-Y_k\|^2-\|T-T_{k+1}\|^2}{2h}
+
\langle e_k,T-Y_k\rangle
+
h\|e_k\|^2.
\end{equation}

Next, for $j\ge0$, introduce the comparison point
\[
\bar T_j
\defeq
\left(1-\frac1{t_{j+1}}\right)T_j
+\frac1{t_{j+1}}T^\star.
\]
The scaled displacements from $\bar T_j$ before and after
the proximal step are
\[
\begin{aligned}
p_j
&\defeq t_{j+1}(Y_j-\bar T_j)
=t_{j+1}Y_j-(t_{j+1}-1)T_j-T^\star,\\
u_{j+1}
&\defeq t_{j+1}(T_{j+1}-\bar T_j)
=t_{j+1}T_{j+1}-(t_{j+1}-1)T_j-T^\star.
\end{aligned}
\]
Set $u_0 \defeq T_0-T^\star$ and define the energy
\[
\mathcal E_0=\frac{\|u_0\|^2}{2h},
\qquad
\mathcal E_j
=t_j^2(f(T_j)-f^\star)+\frac{\|u_j\|^2}{2h},
\quad j\ge1.
\]
Recall from the deterministic proof, as a consequence of the weight recurrence, that
\[
t_{j+1}(t_{j+1}-1)
=(1-qt_{j+1})t_j^2,
\qquad
1\le t_j\le\min\{j,q^{-1/2}\},
\quad j\ge1.
\]
Throughout, inverse powers of $q$ are interpreted as
$+\infty$ when $q=0$.
The momentum formula gives
\[
p_j=c_ju_j+(1-c_j)(T_j-T^\star),
\qquad
c_j\defeq\frac{1-qt_{j+1}}{1-q}\in[0,1],
\]
where the bounds on $c_j$ follow from $q\le1/2$ and
the bounds on $t_{j+1}$.

For $j\ge1$, write $t \defeq t_{j+1}$.
Since $f=g+r$ is $\ell$-strongly convex, we have
\[
t^2(f(\bar T_j)-f^\star)
\le
t(t-1)(f(T_j)-f^\star)
-\frac{\ell}{2}(t-1)\|T_j-T^\star\|^2.
\]
Applying \eqref{eq:WSfista_prox} at $T=\bar T_j$,
multiplying by $t^2$, and using the definitions of
$p_j$ and $u_{j+1}$ therefore gives
\[
\begin{aligned}
\mathcal E_{j+1}
&\le
t(t-1)(f(T_j)-f^\star)
+\frac{(1-q)\|p_j\|^2-q(t-1)\|T_j-T^\star\|^2}{2h} 
-t\langle e_j,p_j\rangle
+ht^2\|e_j\|^2.
\end{aligned}
\]
By convexity of the squared norm, we obtain
\[
\begin{aligned}
(1-q)\|p_j\|^2
&\le
(1-q)c_j\|u_j\|^2
+(1-q)(1-c_j)\|T_j-T^\star\|^2\\
&=
(1-qt)\|u_j\|^2
+q(t-1)\|T_j-T^\star\|^2.
\end{aligned}
\]
The terms involving $\|T_j-T^\star\|^2$ cancel.
Using $t(t-1)=(1-qt)t_j^2$, we obtain
\begin{equation}\label{eq:WSfista_energy}
\mathcal E_{j+1}
\le
(1-qt_{j+1})\mathcal E_j
-t_{j+1}\langle e_j,p_j\rangle
+ht_{j+1}^2\|e_j\|^2.
\end{equation}
For $j=0$, the same inequality follows directly from
\eqref{eq:WSfista_prox} with $T=T^\star$, since
$t_1=1$, $Y_0=T_0$, and $p_0=u_0$.
The separate definition of $\mathcal E_0$ avoids any
finiteness assumption on $f(T_0)$.

To analyze the recursive inequality \eqref{eq:WSfista_energy}, we introduce
\[
\Pi_0 \defeq 1,
\qquad
\Pi_k \defeq \prod_{j=1}^k(1-qt_j),
\qquad
\gamma_k \defeq \frac{\Pi_k}{t_k^2},
\quad k\ge1.
\]
Since $p_j$ is $\mathcal A_j$-measurable, conditional
unbiasedness gives
\[
\E\langle e_j,p_j\rangle
=
\E\bigl\langle\E[e_j\,|\,\mathcal A_j],p_j\bigr\rangle
=0.
\]
Taking expectations in \eqref{eq:WSfista_energy} and
iterating therefore yields, for every $k\ge1$,
\begin{equation}\label{eq:WSfista_unrolled}
\E\mathcal E_k
\le
\frac{\Pi_kR^2}{2h}
+
h\sum_{j=1}^k
\frac{\Pi_k}{\Pi_j}t_j^2v_{j-1}.
\end{equation}
The same calculation as in the deterministic proof gives
$\gamma_{j+1}=(1-t_{j+1}^{-1})\gamma_j$ for $j\ge1$. Since $t_{j+1}^{-1}\ge\sqrt q$, we obtain
\begin{equation}\label{eq:WSfista_gamma_ratio}
\frac{\gamma_k}{\gamma_j}
\le(1-\sqrt q)^{k-j},
\qquad 1\le j\le k.
\end{equation}
Moreover, the bound \eqref{eq:WDfista_gamma} established in the proof of
Theorem~\ref{thm:w-fista_rate} gives
\begin{equation}\label{eq:WSfista_gamma}
\gamma_k
\le
\frac{16}{(k+1)^2}
\exp\left(-\frac{k+1}{2}\sqrt q\right).
\end{equation}

\textbf{Variance control.}
Fix $k\ge1$ and set
\[
V_k\defeq\max_{0\le j<k}v_j,
\qquad
M_{k-1}\defeq\max_{0\le j<k}\E\|u_j\|^2.
\]

We first relate the distances in
\eqref{eq:WSfista_conditional_noise} to the norm term
in the energy. The definitions give
\[
T_j-T^\star
=\frac{u_j}{t_j}
+\left(1-\frac1{t_j}\right)(T_{j-1}-T^\star),
\qquad 1\le j<k,
\]
and
\[
Y_j-T^\star
=\frac{c_j}{t_{j+1}}u_j
+\left(1-\frac{c_j}{t_{j+1}}\right)(T_j-T^\star),
\qquad 0\le j<k.
\]
Both are convex combinations. Since $u_0=T_0-T^\star$,
induction and convexity of the squared norm imply that
\[
\E\|T_j-T^\star\|^2\le M_{k-1},
\qquad
\E\|Y_j-T^\star\|^2\le M_{k-1},
\qquad 0\le j<k.
\]

We next bound the sum in \eqref{eq:WSfista_unrolled}
at each intermediate iteration.
For $1\le j<k$, using $\Pi_j/\Pi_i\le1$ and $t_i\le i$
gives
\[
\sum_{i=1}^j\frac{\Pi_j}{\Pi_i}t_i^2
\le
\sum_{i=1}^{k-1}i^2
\le\frac{k^3}{3}.
\]
When $q>0$, \eqref{eq:WSfista_gamma_ratio} and
$t_j^2\le q^{-1}$ also give
\[
\begin{aligned}
\sum_{i=1}^j\frac{\Pi_j}{\Pi_i}t_i^2
&=
t_j^2\sum_{i=1}^j\frac{\gamma_j}{\gamma_i} \le
q^{-1}\sum_{i=1}^j(1-\sqrt q)^{j-i}
\le q^{-3/2}.
\end{aligned}
\]
Thus, the sum is at most $K_k\defeq\min\{k^3/3,q^{-3/2}\}$.

Since $\Pi_j\le1$ and
$\mathcal E_j\ge\|u_j\|^2/(2h)$,
\eqref{eq:WSfista_unrolled} and the bound
$\|u_0\|^2\le R^2$ give
\[
\begin{aligned}
M_{k-1}
=\max_{0\le j<k}\E\|u_j\|^2
&\le2h\max_{0\le j<k}\E\mathcal E_j\\
&\le
R^2+2h^2\max_{0\le j<k}
\sum_{i=1}^j\frac{\Pi_j}{\Pi_i}t_i^2v_{i-1}\\
&\le
R^2+2h^2K_k\max_{0\le i<k}v_i
\\
&=R^2+2h^2K_kV_k.
\end{aligned}
\]
Combining this with \eqref{eq:WSfista_conditional_noise}
and the definition of $\bar\sigma$ yields
\[
\begin{aligned}
V_k
&\le
2\sigma^2\bigl[
(\zeta_\cH(T_0)+R)^2+M_{k-1}
\bigr] \le
\bar\sigma^2+4\sigma^2h^2K_kV_k.
\end{aligned}
\]

If $\sigma=0$, then $V_k=0$.
Otherwise, under the first step-size restriction in the statement,
\[
4\sigma^2h^2K_k
\le\frac43\sigma^2h^2k^3
\le\frac13.
\]
Under the second restriction, $\ell>0$ and
\[
4\sigma^2h^2K_k
\le4\sigma^2h^2q^{-3/2}
=\frac{4\sigma^2\sqrt h}{\ell^{3/2}}
\le\frac12.
\]
Thus, under either restriction, the coefficient of $V_k$
is at most $1/2$, and so $V_k\le2\bar\sigma^2$.

Finally, \eqref{eq:WSfista_gamma_ratio} gives
\[
\sum_{j=1}^k\frac{\gamma_k}{\gamma_j}
\le
\sum_{j=1}^k(1-\sqrt q)^{k-j}
\le\min\{k,q^{-1/2}\}.
\]
Dividing \eqref{eq:WSfista_unrolled} by $t_k^2$,
discarding the nonnegative norm term, and using
$V_k\le2\bar\sigma^2$, we therefore obtain
\[
\begin{aligned}
\E[f(T_k)-f^\star]
&\le
\frac{\gamma_kR^2}{2h}
+
h\sum_{j=1}^k\frac{\gamma_k}{\gamma_j}v_{j-1}\\
&\le
\frac{\gamma_kR^2}{2h}
+
2\bar\sigma^2h\min\{k,q^{-1/2}\}.
\end{aligned}
\]
Combining this with \eqref{eq:WSfista_gamma} and the
lift identities proves \eqref{eq:WSfista_rate}.

\textbf{Iteration complexity.}
Fix $\eps>0$. Take $N_{\mathrm{det}}$ to be the explicit
sufficient choice constructed in the proof of
Theorem~\ref{thm:w-fista_rate}, with numerical constant $5$,
and define
\[
N_{\mathrm{stoc}}(\eps;\ell,R,\bar\sigma)
\defeq
1+
\begin{cases}
\displaystyle
\left\lceil
1024\,\frac{\bar\sigma^2R^2}{\eps^2}
\right\rceil,
& \displaystyle \eps>\frac{\ell R^2}{4},
\\[4mm]
\displaystyle
\left\lceil
1024\,\frac{\bar\sigma^2}{\ell\eps}
\log\frac{\ell R^2}{\eps}
\right\rceil,
& \displaystyle 0<\eps\le\frac{\ell R^2}{4}.
\end{cases}
\]
This choice satisfies \eqref{eq:iteration-complexity-stoc}.
Fix any integer
\[
k\ge
N_{\mathrm{det}}(\eps;\ell,L,R)
+
N_{\mathrm{stoc}}(\eps;\ell,R,\bar\sigma).
\]
We choose a constant step size for these $k$ iterations and
show that
\[
\Delta_k
\defeq
\E[f(T_k)-f^\star]
\le\eps.
\]
Recall that $\bar\sigma\ge2\sigma R$. Expressions containing
$\bar\sigma$ in the denominator are interpreted as $+\infty$
when $\bar\sigma=0$. First suppose that $\eps>\ell R^2/4$, and choose
\[
h
=
\min\left\{
\frac1{2L},
\frac{R}{\bar\sigma k^{3/2}}
\right\}.
\]
When $\sigma>0$, the inequality $\bar\sigma\ge2\sigma R$
ensures that $h \leq 1/(2\sigma k^{3/2})$, so the first step-size condition holds.
Dropping the exponential factor in
\eqref{eq:WSfista_rate} and using
$\min\{k,q^{-1/2}\}\le k$, one can show
\[
\begin{aligned}
\Delta_k
&\le
\frac{8R^2}{hk^2}
+2\bar\sigma^2hk \le
\frac{16LR^2}{k^2}
+\frac{10\bar\sigma R}{\sqrt{k}} \le
\left(
\frac{16}{25}+\frac{10}{32}
\right)\eps
<\eps.
\end{aligned}
\]

Now suppose that
$0<\eps\le\ell R^2/4$, so in particular $\ell>0$.
Set $a\defeq \ell R^2/\eps \ge4$ 
and choose
\[
h
=
\min\left\{
\frac1{2L},
\frac{\ell\eps^2}{16\bar\sigma^4}
\right\}.
\]
When $\sigma>0$, the bounds
$\bar\sigma\ge2\sigma R$ and
$\eps\le\ell R^2/4$ imply
\[
\frac{\ell\eps^2}{16\bar\sigma^4}
\le
\frac{\ell^3}{4096\sigma^4}
\le
\frac{\ell^3}{64\sigma^4},
\]
so the second step-size condition in the theorem holds.
The stochastic term in \eqref{eq:WSfista_rate} satisfies
\[
2\bar\sigma^2h\min\{k,(h\ell)^{-1/2}\}
\le
2\bar\sigma^2\sqrt{\frac{h}{\ell}}
\le
\frac{\eps}{2}.
\]

For the deterministic term, set $z\!\defeq\!(k\!+\!1)\sqrt{\ell h}$. The lower bound on $k$ and the definition of $h$ give
\[
z
\ge
\min\left\{
\frac5{\sqrt2},256
\right\}\log a
\ge
3\log a.
\]
Indeed, the first quantity comes from
$k\ge5\sqrt{L/\ell}\log a$ when
$h=1/(2L)$, while the second uses $k\ge
1024\,\bar\sigma^2/(\,\ell\eps)\log a$ when $h=\ell\eps^2/(16\bar\sigma^4)$. Hence, the deterministic term in
\eqref{eq:WSfista_rate} is at most
\[
\frac{8\ell R^2}{z^2}e^{-z/2}
\le
\frac{8\eps}{9\sqrt a\,(\log a)^2}
\le
\frac{\eps}{2},
\]
where the last inequality follows from $a\ge4$.

Thus $\Delta_k\le\eps$ in both accuracy regimes.
Since the construction applies to every integer
$k\ge N_{\mathrm{det}}+N_{\mathrm{stoc}}$, with the step size
chosen as a function of $k$, this proves the stated iteration
complexity.
\qed

\subsection{Proof of Proposition~\ref{prop:vi_summary}}

We verify that Assumptions~\ref{assm:lift} and \ref{assm:extras} hold for each $\sQ \in \{\BW, \MF, \Rad\}$. By Proposition~\ref{prop:weak-geodesic-convexity}, each such $\sQ$ satisfies Assumption~\ref{assm:extras}\emph{(i)}. The standard choices of $\cT$ and $\cH$ in Proposition~\ref{prop:standard_choices} satisfy Assumptions~\ref{assm:lift}\emph{(i)} and \ref{assm:extras}\emph{(ii)}. Since $g$ is the canonical lift, we have by Proposition~\ref{prop:interaction-gradient}, Corollary~\ref{cor:interaction-energy-Q-grads}, and Proposition~\ref{prop:gradient_connection} that $g$ is Fréchet differentiable with gradient $\nabla_\cH g(T) =\mathrm{proj}_\cH[\nabla_\sQ \cG(T_\sharp \rho) \circ T]$, where the $\sQ$-gradients are given in Corollary~\ref{cor:interaction-energy-Q-grads}. By Proposition~\ref{prop:interaction_energy_lift}, $g$ is $\ell_V$-strongly convex and $L_V$-smooth. Hence, $g$ satisfies Assumptions~\ref{assm:lift}\emph{(ii)} and \ref{assm:extras}\emph{(iii)}. Since $\rho=\cN(0,I_d)\in\sQ$ and
$-\Ent(\rho)<+\infty$, Proposition~
\ref{prop:nonsmooth_lift_regular} implies that
$r$ is proper, convex, and lower semicontinuous, and thus satisfies Assumption~\ref{assm:lift}\emph{(iii)}. Altogether, $\rho$, $\cT$, $\cH$, $g$ and $r$ satisfy Assumptions~\ref{assm:lift} and \ref{assm:extras} for each $\sQ \in \{\BW, \MF, \Rad\}$.

We verify the gradient formulas of Table~\ref{tab:vi_summary}. To apply Proposition~\ref{prop:stochastic_gradient_connection},
we take $\eta(T,\xi)=T(\xi)$ with $\xi\sim\rho$
for $\sQ\in\{\BW,\MF\}$, and $\eta(T,\xi)=\xi$
with $\xi\sim\Unif(\cS^{d-1})$ for $\sQ=\Rad$. By  Proposition~\ref{prop:gradient_connection} and Proposition~\ref{prop:stochastic_gradient_connection}, we have
\[
\begin{aligned}
    &\nabla_\cH g(T)
    =
    \mathrm{proj}_\cH[\nabla_\sQ\cG(\mu)\circ T],\\
    &\cO_\cH(T,\xi)
    =
    \mathrm{proj}_\cH[\cO_\sQ(T_\sharp\rho, \eta(T,\xi))\circ T],
\end{aligned}
\]
where $\nabla_\mathsf{Q}\cG$ and $\cO_\sQ$ are as in  Corollary~\ref{cor:interaction-energy-Q-grads} and Proposition~\ref{prop:stochastic-Q-examples}, respectively. These are the deterministic and stochastic $\cH$-gradients of $g$ (up to verifying measurability of $\cO_\cH$). We evaluate the $\cH$-gradients below for the different cases of $\sQ$. In each case, if $T \in \cH$, let $\mu = T_\sharp\rho$.

\textbf{Gaussian case ($\sQ=\BW$).}
Let $T \in \cH$. Then, $T(x)=m+Ax$, where $A=A^\top$. The $\BW$-gradient is
$\nabla_\BW\cG(\mu)(y)=\EE_\rho[\nabla V(T(X))]+\EE_\rho[\nabla^2V(T(X))](y-m)$.
Thus,
\[
\begin{aligned}
    \nabla_\cH g(T)(x)
    &=
    \mathrm{proj}_\cH[\EE_\rho[\nabla V(T(X))]+\EE_\rho[\nabla^2V(T(X))]A(\cdot)](x)\\
    &=
    \EE_\rho[\nabla V(T(X))]+\mathrm{sym}(\EE_\rho[\nabla^2V(T(X))]A)x,\\
    \cO_\cH(T,\xi)(x)
    &=
    \mathrm{proj}_\cH\left[
        \nabla V(T(\xi))
        +\nabla^2V(T(\xi))A(\cdot)
    \right](x)\\
    &=
    \nabla V(T(\xi))
    +\mathrm{sym}(\nabla^2V(T(\xi))A)x.
\end{aligned}
\]
where $\mathrm{sym}(B) \coloneq (B + B^\top)/2$ for $B \in \RR^{d \times d}$ is the symmetrization operator.

\textbf{Mean-field case ($\sQ=\MF$).}
Let $T \in \cH$. Since coordinate-wise separability is preserved under composition, both $\nabla_\MF\cG(\mu)\circ T$ and $\cO_\MF(\mu,T(\xi))\circ T$ belong to $\cH$. Thus, the projection
onto $\cH$ leaves them unchanged,
\[
\begin{aligned}
    \nabla_\cH g(T)(x)
    &=
    \mathrm{proj}_\cH[\nabla_\MF\cG(\mu)\circ T](x)\\
    &=
    \left(
        \EE_\rho\left[
            \partial_iV\bigl(
                T_i(x_i),T_{-i}(X_{-i})
            \bigr)
        \right]
    \right)_{i=1}^d,\\
    \cO_\cH(T,\xi)(x)
    &=
    \mathrm{proj}_\cH[\cO_\MF(\mu,T(\xi))\circ T](x)\\
    &=
    \left(
        \partial_iV\bigl(
            T_i(x_i),T_{-i}(\xi_{-i})
        \bigr)
    \right)_{i=1}^d.
\end{aligned}
\]

\textbf{Radial case ($\sQ=\Rad$).}
Let $T \in \cH$, and write $T(x)=\Psi_T(\|x\|_2) x /\|x\|_2$ when $x \neq 0$ and $T(x) = 0$ otherwise. Since compositions of radial maps are radial,
both $\nabla_\Rad\cG(\mu)\circ T$ and
$\cO_\Rad(\mu,\xi)\circ T$ belong to $\cH$.
Thus, the projection onto $\cH$ leaves them
unchanged. For $x$ such that $T(x)\neq0$, we obtain
\[
\begin{aligned}
    \nabla_\cH g(T)(x)
    &=
    \mathrm{proj}_\cH[\nabla_\Rad\cG(\mu)\circ T](x)\\
    &=
    \EE_{\mathrm{Unif}(\cS^{d-1})}\left[\left\langle
        \nabla V(\|T(x)\|_2 X),X
    \right\rangle\right]
    \frac{T(x)}{\|T(x)\|_2},\\
    \cO_\cH(T,\xi)(x)
    &=
    \mathrm{proj}_\cH[\cO_\Rad(\mu,\xi)\circ T](x)\\
    &=
    \left\langle
        \nabla V(\|T(x)\|_2\xi)-\nabla V(0),\xi
    \right\rangle
    \frac{T(x)}{\|T(x)\|_2},
\end{aligned}
\]
while both expressions are defined to be zero whenever $T(x)=0$.

\textbf{Variance bounds.} The stochastic $\cH$-gradient retains the variance constants of the $\sQ$-gradients. The noise scaling function follows from
\[
    \zeta_\cH(T)^2
    =\zeta_\sQ(T_\sharp\rho)^2
    =\|T-\EE_\rho[T(X)]\|_{L^2(\rho)}^2.
\]

\textbf{Measurability.} Finally, we verify joint measurability of the
oracles $\cO_\cH$. The Gaussian and radial
oracles are jointly continuous and hence jointly measurable, by continuity of
$\nabla^2 V$ and Lipschitz continuity of $\nabla V$,
respectively.
For the MF oracle, the coordinate maps admit
representatives for which the evaluation maps
$(T_i,x_i)\mapsto T_i(x_i)$ are jointly Borel
measurable, by an application of Lévy's upward
convergence theorem
\citep[Theorem~9.24]{kallenberg2021foundations}.
Joint measurability of $\cO_\cH$ then follows
by composing the jointly measurable map
$(T,\xi)\mapsto (T,T(\xi))$ with the map
\[
    (T,z)\mapsto
    \left(\partial_iV(T_i(\cdot_i),z_{-i})\right)_{i=1}^d
    \in\cH,
\]
which is continuous by Lipschitz continuity
of $\nabla V$.
\qed

\section{Proofs for Section~\ref{sec:tractable}}

\subsection{Proof of Proposition~\ref{prop:sufficient_g}}

Let $\cG(\mu) = \int V\dd\mu^{\otimes n}$. Under the stated assumptions on $V$,
Proposition~\ref{prop:interaction_energy_lift}
ensures that the canonical lift
$g(T)=\cG(T_\sharp\rho)$ on $\cH$ is differentiable,
$n\ell_V$-strongly convex and $nL_V$-smooth,
and hence satisfies Assumption~\ref{assm:lift}\emph{(ii)}. By Proposition~\ref{prop:gradient_connection}
and Corollary~\ref{cor:interaction-energy-Q-grads},
setting $\mu=(T_{b,A})_\sharp\rho$, we have
$$\nabla_\cH g(T_{b,A})(x) = \mathrm{proj}_\cH[\nabla_\mathsf{BW}\cG(\mu) \circ T_{b,A}](x) = \EE_\mu[w_\mu(X)] + \mathrm{proj}_\HH [\EE_\mu[Dw_\mu(X)]A] x.$$
The formula for $\nabla G(b,A)$ follows by identifying the affine coefficients of $\nabla_\cH g(T_{b,A})$.

Let $\cG(\mu) = G_1(m_\mu, \Sigma_\mu)$. Consider the canonical lift $g(T_{b,A}) = G_1(b, AA^\top)$ then the coefficient function is $
G(b,A)=G_1(b,AA^\top).$ The $\ell_1$-strong convexity and $L_1$-smoothness follow by assumption. It suffices to show the form of the gradient of $G$.
Consider the map
\[
G_0:\RR^d\times\HH\to\RR^d\times\SS_+^d, \quad
G_0(b,A)=(b,AA^\top).
\]
The coefficient map of $g$ is $G = G_1 \circ G_0$. Note that $G_0$ is Fréchet differentiable as a map into
$\RR^d\times\SS^d$, with Fréchet derivative $DG_0(b,A)[c,C]=(c,CA^\top+AC^\top)$ for all $(c,C) \in \RR^d \times \HH$. Since $G_1$ is Fr\'echet differentiable relative to
$\RR^d\times\SS_+^d$ by assumption, we have that the Fréchet derivative is representable as $DG_1(x,X)[d,D] = \langle \nabla_1G_1(x,X), d\rangle + \langle \nabla_2 G_1(x,X), D\rangle_\mathrm{F}$ for all $(x,X), (d,D) \in \RR^d \times \SS^d$. The relative chain rule
\cite[Lemma~2.4]{fischer2007extending}
implies that $G=G_1\circ G_0$ is Fr\'echet differentiable
on $\RR^d\times\HH$. In particular,
\[
\begin{aligned}
DG(b,A)[c,C] &= (DG_1(G_0(b,A)) \circ DG_0(b,A))[c,C]\\
&=
\langle\nabla_1G_1(b,AA^\top),c\rangle +
\langle\nabla_2 G_1(b,AA^\top),
CA^\top+AC^\top\rangle_{\mathrm F}\\
&=
\langle\nabla_1 G_1(b,AA^\top),c\rangle
+
\langle2\nabla_2 G_1(b,AA^\top)A,
C\rangle_{\mathrm F}\\
&=
\langle\nabla_1 G_1(b,AA^\top),c\rangle
+
\langle2 \mathrm{proj}_\HH[\nabla_2 G_1(b,AA^\top)A],
C\rangle_{\mathrm F},
\end{aligned}
\]
where the first equality follows from the relative chain rule and the second equality from evaluating the composition. The third equality uses that
\begin{align*}
\langle \nabla_2 G_1(b,AA^\top), CA^\top \rangle_\mathrm{F} = \langle \nabla_2 G_1(b,AA^\top), AC^\top \rangle_\mathrm{F} = \langle \nabla_2 G_1(b,AA^\top)A, C \rangle_\mathrm{F},
\end{align*}
since $\nabla_2 G_1(b,AA^\top)\in\SS^d$ and due to properties of the trace. The fourth equality follows since $C \in \HH$. We have thus identified the Fréchet gradient of $G$ as
\[
\nabla G(b,A)
=
\left(
\nabla_1G_1(b,AA^\top),
\mathrm{proj}_\HH
\left[2\nabla_2 G_1(b,AA^\top)A\right]
\right).
\]

Let $\cG(\mu) = G_2(m_\mu, \Sigma_\mu^{1/2})$.
Here $\TT=\SS_+^d$ and $\HH=\SS^d$. Assumption~\ref{assm:lift}\emph{(i)} holds by the
matrix square root construction and since
$I_d\in\TT$. For $A\in\TT$, $(AA^\top)^{1/2}=A$. Therefore, the lift $g(T_{b,A}) = G_2(b,A)$ agrees with $\cG((T_{b,A})_\sharp\rho) = G_2(b, (AA^\top)^{1/2})$ on $\RR^d \times \TT$. Additionally, $G_2$ is $\ell_2$-strongly convex, $L_2$-smooth and differentiable, and thus, $g$ has these properties with coefficient gradient $\nabla G(b,A)
=
\nabla G_2(b,A).$

Let $\cG(\mu) = \frac{1}{2}W_2^2(\mu, \nu)$ and
$\TT=\SS_+^d\Sigma_\nu^{1/2}$ and $\HH=\SS^d\Sigma_\nu^{1/2}$. Then $\TT$ is closed and convex and $\HH$ is a linear subspace. We have $I_d = \Sigma_\nu^{-1/2}\Sigma_\nu^{1/2} \in \TT$. Moreover, $\TT \ni A \mapsto AA^\top \in \SS^d_+$ is surjective. Indeed, for $\Sigma \in \SS^d_+$, note that $$A = \Sigma_\nu^{-1/2}(\Sigma_\nu^{1/2}\Sigma\Sigma_\nu^{1/2})^{1/2}\Sigma_\nu^{-1/2} \Sigma_\nu^{1/2} = \Sigma_\nu^{-1/2}(\Sigma_\nu^{1/2}\Sigma\Sigma_\nu^{1/2})^{1/2} \in \TT,$$ satisfies $AA^\top = \Sigma_\nu^{-1/2}(\Sigma_\nu^{1/2}\Sigma\Sigma_\nu^{1/2})\Sigma_\nu^{-1/2} = \Sigma$. Thus, $\cT$ satisfies Assumption~\ref{assm:lift}\emph{(i)}.
For $A\in\TT$, write $A=S\Sigma_\nu^{1/2}$ with $S\in\SS_+^d$.
The affine map
\[
y\longmapsto b+S(y-m_\nu)
\]
is an optimal transport map from $\nu=\cN(m_\nu,\Sigma_\nu)$ to
$(T_{b,A})_\sharp \rho = \cN(b, S\Sigma_\nu S^\top)$. Consequently, $(b + S(\cdot - m_\nu),\mathrm{id})_\sharp \nu = (T_{b,A}, T_{m_\nu, B})_\sharp \rho$, where $B = \Sigma_\nu^{1/2}$, is an optimal coupling, and
\[
\begin{aligned}
W_2^2((T_{b,A})_\sharp\rho,\nu) = 
\EE_{Z\sim\rho}
\|b + AZ -(m_\nu + \Sigma_\nu^{1/2}Z)\|_2^2 =
\|b-m_\nu\|_2^2+\|A-\Sigma_\nu^{1/2}\|_{\mathrm F}^2.
\end{aligned}
\]
Thus,
\[
G(b,A)
=
\frac12\|b-m_\nu\|_2^2
+
\frac12\|A-\Sigma_\nu^{1/2}\|_{\mathrm F}^2
\]
agrees with
$\cG((T_{b,A})_\sharp\rho)$ on
$\RR^d\times\TT$ and defines a function on
$\RR^d\times\HH$. Clearly, $G$ is the $1$-strongly convex and $1$-smooth quadratic on $\RR^d \times \HH$. Finally, since $\Sigma_\nu^{1/2}\in\HH$,
we have $A-\Sigma_\nu^{1/2}\in\HH$ for every
$A\in\HH$, and hence the gradient is $$\nabla G(b,A)=(b-m_\nu,A-\Sigma_\nu^{1/2}).$$ \qed

\subsection{Proof of Corollary~\ref{cor:BW_variance}}

Fix $T_{b,A}\in\cH$ and set $\mu=(T_{b,A})_\sharp\rho$.
Applying Proposition~\ref{prop:stochastic_gradient_connection}
with $\eta(T,\xi)=T(\xi)$ and using the stochastic $\BW$-gradient oracle
in Proposition~\ref{prop:stochastic-Q-examples}
yields the stated expression for $\cO_\cH$.
Moreover, $\sigma$ is inherited from the stochastic $\BW$-gradient oracle and
\[
    \zeta_\cH(T_{b,A})
    =\zeta_\BW(\mu)
    =\sqrt{\tr(\Sigma_\mu)}
    =\|A\|_{\mathrm F}.
\]
The integrands defining $w_\mu$ and
$Dw_\mu$ under the affine parametrization
are continuous in $(b,A,\xi,x)$, and hence jointly
measurable, by continuity of $\nabla V$ and
$\nabla^2 V$, respectively.
Joint measurability of $w_\mu$ and
$Dw_\mu$, and consequently that of $\cO_\cH$, follows from
\cite[Lemma~1.28(ii)]{kallenberg2021foundations}.
\qed

\subsection{Proof of Proposition~\ref{prop:sufficient_H}}

We extend $R_2$ to $\mathbb R^d$ by $+\infty$ outside $\mathbb R_+^d$, and use the same notation for the extension. We have
\[
    r(T_{b,A})=R_1(b)+R_2(\lambda(A)),
    \quad (b,A)\in\RR^d\times\HH.
\]
Indeed, equality holds if $A\in\TT$ since $(AA^\top)^{1/2}=A$, and both sides are $+\infty$ if $A \in \HH \setminus \TT$.
Since $R_2$ is proper, lower semicontinuous, convex,
and permutation invariant, the spectral function
$\mathbb S^d\ni A\mapsto R_2(\lambda(A))$
is proper, lower semicontinuous, and convex
\citep[Corollary~2.4]{lewis1996convex}. Thus, by the isometric isomorphism between $\RR^d\times\HH$ and $\cH$, the lift $r$ is
proper, convex, and lower semicontinuous. The characterization of $b^\star$ follows by separability, while the characterization of $A^\star$ follows by \cite[Sec.~6.7.2, Eq.~(6.11)]{parikh2014proximal}. \qed
    
\subsection{Proof of Proposition~\ref{prop:MF_approximation}}

For a given $B>0$, $J \in \NN_+$ and coinciding discretization points, the set
$\cT_\Phi$ contains the piecewise affine coordinate-wise separable transport
maps constructed in \cite[Appendix~C.3]{jiang2025algorithms},
since our parametrization allows arbitrary
nonnegative tail slopes. Hence, by \cite[Theorem~5.6 and Corollary C.3]{jiang2025algorithms}, for every $\eta \in (0,1]$, there exists $\widehat \mu_{\eta} \in \mathsf{MF}_\Phi$ such that $$\mathrm{KL}(\widehat \mu_{\eta}\|\pi)
-
\mathrm{KL}(\pi^\star\|\pi)
\lesssim
\kappa^3\sqrt d\,\eta$$ if $B = O(\sqrt{\log(d/\eta^2)})$ and $J = \widetilde{O}(\kappa^{1/2}d^{1/4}\eta^{-1/2})$. Hence, for $\varepsilon \in (0,1]$, let $\eta(\varepsilon) = \varepsilon/(\kappa^3\sqrt{d})$, such that $\mathrm{KL}(\widehat \mu_{\eta(\varepsilon)}\|\pi)
-
\mathrm{KL}(\pi^\star\|\pi)
\lesssim
\varepsilon$ for suitable choices of $B = O(\sqrt{\log(\kappa^6d^2/\varepsilon^2)}) = O(\sqrt{\log(\kappa^3d/\varepsilon)})$ and $J = \widetilde{O}(\kappa^2d^{1/2}\varepsilon^{-1/2})$. Thus,
\begin{align*}
\frac{\ell_V}{2}W_2(\pi^\star_\Phi,\pi^\star)^2 &\leq \mathrm{KL}(\pi^\star_\Phi\|\pi)
-
\mathrm{KL}(\pi^\star\|\pi) \\
&\leq 
\mathrm{KL}(\widehat\mu_{\eta(\varepsilon)}\|\pi)
-
\mathrm{KL}(\pi^\star\|\pi)\\
&\lesssim
\varepsilon.
\end{align*}
The first inequality follows from the $\ell_V$-strong
geodesic convexity of the KL divergence
\cite[Proposition~5.10]{jiang2025algorithms} and weak geodesic convexity of $\MF$ by Proposition~\ref{prop:weak-geodesic-convexity}. The second inequality follows from optimality of
$\pi_\Phi^\star$ over $\mathsf{MF}_\Phi$ and
$\widehat\mu_{\eta(\varepsilon)}\in\mathsf{MF}_\Phi$. \qed

\subsection{Proof of Proposition~\ref{prop:Rad_approximation}}

For given $B\in(0,\sqrt d)$, $J\in\NN_+$, and
coinciding discretization points, the set $\cT_\Phi$
contains the piecewise affine radial transport maps
constructed in \cite[Eqs.~(5)--(6)]{ghafourpour2026variational},
since our parametrization allows an arbitrary
nonnegative tail slope.

Let $T^\star=T^{\rho\to\pi^\star}$.
By \cite[Theorems~3.5 and~4.1 and Appendix~C.1,
Steps~4--5]{ghafourpour2026variational},
for every $\eta\in(0,1]$ satisfying
$\log(d/\eta)\leq c_0d$, where $c_0>0$ is a
sufficiently small constant, there exists
$\widehat T_\eta\in\cT_\Phi$ such that
\[
    \|\widehat T_\eta-T^\star\|_{L^2(\rho)}^2
    \lesssim \frac{\eta^2}{\ell_V},
    \qquad
    \|D(\widehat T_\eta-T^\star)\|_{L^2(\rho)}^2
    \lesssim
    \frac{\kappa\eta\sqrt{\log(d/\eta)}}{\ell_V},
\]
for suitable choices
$B=O(\sqrt{\log(d/\eta)})$ and
$J=\widetilde O(\kappa^{1/2}\eta^{-1/2})$.
The norm in the second inequality is induced by the
Frobenius norm.

Applying the argument in the proof of
\cite[Corollary~C.3]{jiang2025algorithms},
with the first variation vanishing by the
characterization of the radial minimizer in
\cite[Proposition~3.3]{ghafourpour2026variational},
we obtain, for
$\widehat\mu_\eta=(\widehat T_\eta)_\sharp\rho$,
\begin{align*}
\KL(\widehat\mu_\eta\Vert\pi)
-\KL(\pi^\star\Vert\pi)
&\leq
\frac{L_V}{2}
\left(
\|\widehat T_\eta-T^\star\|_{L^2(\rho)}^2
+\kappa\|D(\widehat T_\eta-T^\star)\|_{L^2(\rho)}^2
\right)\\
&\lesssim
\kappa^3\eta\sqrt{\log(d/\eta)}.
\end{align*}

Now set
\[
    \eta(\varepsilon)
    =
    \frac{\varepsilon}
    {\kappa^3\sqrt{\log(d\kappa^3/\varepsilon)}}.
\]
Using \eqref{eq:eta} with $c=2c_0/3$, it follows
that $\log(d/\eta(\varepsilon))\leq c_0d$.
Consequently,
\[
    \KL(\widehat\mu_{\eta(\varepsilon)}\Vert\pi)
    -\KL(\pi^\star\Vert\pi)
    \lesssim
    \kappa^3\eta(\varepsilon)
    \sqrt{\log(d/\eta(\varepsilon))}
    \lesssim\varepsilon
\]
for suitable choices
$B=O(\sqrt{\log(d\kappa^3/\varepsilon)})$ and
$J=\widetilde O(\kappa^2\varepsilon^{-1/2})$.
Moreover,
\begin{align*}
\frac{\ell_V}{2}W_2(\pi^\star_\Phi,\pi^\star)^2
&\leq
\KL(\pi^\star_\Phi\Vert\pi)
-\KL(\pi^\star\Vert\pi)\\
&\leq
\KL(\widehat\mu_{\eta(\varepsilon)}\Vert\pi)
-\KL(\pi^\star\Vert\pi)\\
&\lesssim\varepsilon.
\end{align*}
The first inequality follows from the $\ell_V$-strong
geodesic convexity of the KL divergence
\cite[Proposition~5.10]{jiang2025algorithms},
weak geodesic convexity of $\Rad$
(Proposition~\ref{prop:weak-geodesic-convexity}),
and optimality of $\pi^\star$ over $\Rad$.
The second inequality follows from optimality of
$\pi_\Phi^\star$ over $\mathsf{Rad}_\Phi$ and
$\widehat\mu_{\eta(\varepsilon)}\in\mathsf{Rad}_\Phi$.
\qed

    \end{document}